\documentclass[12pt,psamsfonts,leqno,oneside,letterpaper]{amsart}
\usepackage[dvips,text={6.5truein,9truein},left=1truein,top=1truein]{geometry}
\usepackage{amssymb,amsmath,amscd,enumerate}
\usepackage[pdftex]{graphicx}
\usepackage{url}

\usepackage[colorlinks,linkcolor=blue,citecolor=blue,pdfstartview=FitH]{hyperref}
\input xy
\xyoption{all}
\SelectTips{cm}{12}

\usepackage{color}

\theoremstyle{definition}
\newtheorem{para}{}[section]

\newtheorem{remark}[para]{Remark}
\newtheorem{remarks}[para]{Remarks}
\newtheorem{notation}[para]{Notation}
\newtheorem{convention}[para]{Convention}
\newtheorem{definition}[para]{Definition}
\newtheorem{definitions}[para]{Definitions}

\newcommand\Alternatives{\begin{enumerate}[(i)]}
\newcommand\EndAlternatives{\end{enumerate}}
\newcommand\Conditions{\begin{enumerate}[(1)]}
\newcommand\EndConditions{\end{enumerate}}

\theoremstyle{plain}
\newtheorem{theorem}[para]{Theorem}
\newtheorem{lemma}[para]{Lemma}
\newtheorem{proposition}[para]{Proposition}
\newtheorem{corollary}[para]{Corollary}
\newtheorem{conjecture}[para]{Conjecture}
\newtheorem{claim}[equation]{}

\numberwithin{equation}{para}
\numberwithin{figure}{section}

\newcommand\Number{\begin{para}}
\newcommand\EndNumber{\end{para}}
\newcommand\Definition{\begin{definition}}
\newcommand\EndDefinition{\end{definition}}
\newcommand\Definitions{\begin{definitions}}
\newcommand\EndDefinitions{\end{definitions}}
\newcommand\Theorem{\begin{theorem}}
\newcommand\EndTheorem{\end{theorem}}
\newcommand\Conjecture{\begin{conjecture}}
\newcommand\EndConjecture{\end{conjecture}}
\newcommand\Remark{\begin{remark}}
\newcommand\EndRemark{\end{remark}}
\newcommand\Remarks{\begin{remarks}}
\newcommand\EndRemarks{\end{remarks}}
\newcommand\Convention{\begin{convention}}
\newcommand\EndConvention{\end{convention}}
\newcommand\Notation{\begin{notation}}
\newcommand\EndNotation{\end{notation}}
\newcommand\Lemma{\begin{lemma}}
\newcommand\EndLemma{\end{lemma}}
\newcommand\Proposition{\begin{proposition}}
\newcommand\EndProposition{\end{proposition}}
\newcommand\Corollary{\begin{corollary}}
\newcommand\EndCorollary{\end{corollary}}
\newcommand\Claim{\begin{claim}}
\newcommand\EndClaim{\end{claim}}
\newcommand\Proof{\begin{proof}}
\newcommand\EndProof{\end{proof}}
\newcommand\Equation{\begin{equation}}
\newcommand\EndEquation{\end{equation}}
\newcommand\NoProof{{\hfill$\square$}}
\newcommand\Bullets{\begin{itemize}}
\newcommand\EndBullets{\end{itemize}}

\newcommand\Fix{\mathop{\rm Fix}}
\newcommand\Per{\mathop{\rm Per}}

\newcommand\tF{\widetilde F}

\newcommand\trace{\mathop{\rm trace}}

\newcommand\hGamma{{\widehat\Gamma}}
\newcommand\hM{{\widehat M}}

\newcommand\diam{\mathop{\rm diam}}

\newcommand\tM{\widetilde M}

\newcommand\tGamma{\widetilde \Gamma}

\newcommand\tC{\widetilde C}

\newcommand\tN{\widetilde N}

\newcommand\tPhi{\widetilde \Phi}

\newcommand\tQ{\widetilde Q}
\newcommand\tx{\widetilde x}
\newcommand\ty{\widetilde y}

\newcommand\cald{{\mathcal D}}
\newcommand\cale{{\mathcal E}}
\newcommand\calF{{\mathcal F}}
\newcommand\calL{{\mathcal L}}

\newcommand\calt{{\mathcal T}}

\newcommand\frako{{\mathfrak o}}
\newcommand\frakm{{\mathfrak m}}
\newcommand\frakT{{\mathfrak T}}
\newcommand\calo{{\mathcal O}}
\newcommand\calm{{\mathcal M}}
\newcommand\caln{{\mathcal N}}
\newcommand\calr{{\mathcal R}}

\newcommand\calv{{\mathcal V}}

\newcommand\Aut{\mathop{\rm Aut}}
\newcommand\meas{\mathop{\rm meas}}
\newcommand\St{{\rm St}}

\newtheorem*{TheoremA}{Theorem A}
\newtheorem*{TheoremB}{Theorem B}

\newcommand\ZZ{{\mathbb Z}}
\newcommand\FF{{\mathbb F}}

\newcommand\NN{{\mathbb N}}
\newcommand\CC{{\mathbb C}}
\newcommand\QQ{{\mathbb Q}}
\newcommand\HH{{\mathbb H}}

\newcommand\op{{\mathfrak o}_{(\frakp)}}
\newcommand\cala{{\mathcal A}}

\newcommand\calp{{\mathcal P}}
\newcommand\cals{{\mathcal S}}

\newcommand\frakp{{\mathfrak p}}

\newcommand\dist{\mathop{\rm dist}}

\newcommand\vol{\mathop{\rm vol}}

\newcommand\length{\mathop{{\rm length}}}

\newcommand\isomplus{\mathop{{\rm Isom}_+}}

\newcommand\pizzle{{\rm PSL}_2}
\newcommand\zzle{{\rm SL}_2}
\newcommand\piggle{{\rm PGL}_2}
\newcommand\ggle{{\rm GL}_2}
\newcommand\ok{{\mathcal O}_K}

\newcommand\ooks{\calo_{K,S}}

\begin{document}
\author{Peter B. Shalen}
\address{Department of Mathematics, Statistics, and Computer Science (M/C 249)\\  University of Illinois at Chicago\\
  851 S. Morgan St.\\
  Chicago, IL 60607-7045} \email{shalen@math.uic.edu}
\thanks{Partially supported by NSF {grant DMS-0906155}}

\title{Margulis numbers and number fields}

\begin{abstract}
It is shown that, up to isometry, all but finitely many closed, orientable hyperbolic $3$-manifolds with a given trace field $K$ admit $0.34$ as a Margulis number. This is deduced from a more technical result giving a condition under which $\max(d(P,x\cdot P),d(P,y\cdot P))\ge0.34$ for every $P\in\HH^3$, where $x$ and $y$ lie in $\pizzle(E)$ for some number field $E$, generate a discrete torsion-free group of $\pizzle(\CC)$ and do not commute. Specifically, this is always the case if there is a valuation  $v$  of $E$ such that (1) the residue field $k_v=\frako_v/\frakm_v$ of $v$ has sufficiently large characteristic, (2) $x\in\pizzle(\frako_v)$, and (3) the image of $x$ under the natural homomorphism $\pizzle(\frako_v)\to \pizzle(k_v)$ has order $7$.
\end{abstract}

\maketitle

\section{Introduction}

This paper will illustrate a surprising interaction between the quantitative geometry of hyperbolic $3$-manifolds and their number-theoretical properties. Before stating the main results I will give some conventions and a little background.

\Definitions\label{Margdef}
I will denote by $d$ the  hyperbolic distance on $\HH^n$. If $P$ is a point of $\HH^n$ and $x$ is a hyperbolic isometry, then $d_P(x)$ will denote $d(P,x\cdot P)$.

Let $M$ be a closed, orientable hyperbolic $n$-manifold. We may write
$M=\HH^n/\Gamma$ where $\Gamma\le\isomplus(\HH^n)$ is discrete,
cocompact and torsion-free. The group $\Gamma$ is uniquely determined
by $M$ up to conjugacy in $\isomplus(\HH^n)$. We define a {\it
Margulis number} for $M$, or for $\Gamma$, to be a positive real
number $\mu$ with the following property:
\Claim
If $P$ is a point of $\HH^n$, and if $x$ and $y$ are elements of $\Gamma$ such
that $\max(d_P(x),d_P(y))<\mu$, then $x$ and $y$ commute.
\EndClaim
\EndDefinitions

The Margulis Lemma \cite[Chapter D]{bp} implies that for every
$n\ge2$ there is a positive constant which is a Margulis number for
every closed, orientable hyperbolic $n$-manifold. The largest such number,
$\mu_+(n)$, is called the {\it Margulis constant} for closed, orientable hyperbolic
$n$-manifolds.

Margulis numbers play a central role in the geometry of hyperbolic
manifolds. If $\mu$ is a Margulis number for $M$ then the points of
$M$ where the injectivity radius is less than $\mu/2$ form a disjoint
union of ``tubes'' about closed geodesics whose geometric structure
can be precisely described. Topologically they are open $(n-1)$-ball
bundles over $S^1$. This observation and the Margulis Lemma can be
used to show, for example, that for every $V>0$ there is a finite
collection of compact orientable $3$-manifolds $M_1,\ldots,M_N$, whose
boundary components are tori, such that every closed, orientable
hyperbolic $3$-manifold of volume at most $V$ can be obtained by a
Dehn filling of one of the $M_i$.

Meyerhoff showed in \cite{meyerhoff} that $\mu_+(3)>0.104$. Marc Culler has informed me that according to strong numerical evidence, 
$0.616$ fails to be a Margulis number for the hyperbolic $3$-manifold {\tt m027(-4,1)}, and hence 
$\mu_+(3) < 0.616$.          

The main result of  \cite{finiteness} asserts that up to isometry, all but at most finitely many orientable hyperbolic $3$-manifolds admit $0.29$ as a Margulis number. 

The motivating result of this paper, Theorem A below, asserts that results of the type proved in \cite{finiteness} can be considerably strengthened if one imposes a number-theoretical restriction on a hyperbolic $3$-manifold by specifying its trace field, of which I will now recall the definition.

\Number\label{integral traces}
If $\Gamma$ is a cocompact discrete subgroup of $\zzle(\CC)$, and if $\tGamma$ denotes the preimage of $\Gamma$ under the quotient homomorphism $\zzle(\CC)\to\pizzle(\CC)$, it follows from
\cite[Theorem 3.1.2]{mr} that the traces of the elements of $\tGamma$ generate a finite extension of
$\QQ$, called the {\it trace field} of $\Gamma$. It is also referred to as the trace field of the hyperbolic $3$-manifold $M=\HH^3$.
\EndNumber

\begin{TheoremA}
 Let $K$ be a number field. Then up to isometry, among the closed, orientable hyperbolic $3$-manifolds that have trace field $K$, all but (at most) finitely many admit $0.34$ as a Margulis number.
\end{TheoremA}

This will be proved in the body of the paper as Theorem \ref {motivational research}.

The proof of Theorem A will depend on a more technical result, which I state below as Theorem B. Before giving the statement it will be convenient to make a few more definitions that are used throughout the paper.

\Number\label{lions and tigers and bears}
If $K$ is a field, I will denote by $\Pi_K$ the quotient homomorphism $\ggle(K)\to\piggle(K)$, which by restriction maps $\zzle(K)$ to $\pizzle(K)$.
For a matrix $A$ in $\ggle(K)$, I will write $[A]=\Pi_K(A)$.
\EndNumber

\Number \label{bet you can't eat one}If $\zeta:R\to S$ is a homomorphism of commutative rings, I will denote by $h_\zeta$ the group homomorphism from $\ggle(R)$ to $\ggle(S)$ defined by
$$\begin{pmatrix}a&b\\c&d\end{pmatrix}\mapsto\begin{pmatrix}\zeta(
a)&\zeta( b)\\\zeta( c)&\zeta( d)\end{pmatrix}.$$
I will denote by $\overline h_\zeta: \piggle(R)\to\piggle(S)$
the group homomorphism defined by $\overline h_\zeta\circ\Pi_R= \Pi_S\circ h_\zeta$. 
%By a {\it valuation} of a field $K$ I will mean a surjective homomorphism $v$ from the multiplicative group $K^\times$ to the additive group $\ZZ$ which satisfies the condition that $v(x+y)\ge\min (v(x),v(y))$ whenever $x$, $y$ and $x+y$ belong to $K^\times$. (These are sometimes called normalized discrete rank-$1$ valuations.) 

If $v$ is a valuation of a field $K$, I will denote by $\frako_v$ the corresponding valuation ring. (A general reference for valuations is \cite[Chapter II]{neukirch}.)
%$$\frako=\{0\}\cup\{x\in K:v(x)\ge0\}.$$
The unique maximal ideal of $\frako_v$
%, which may be defined as $\{0\}\cup\{x\in E:v(x)\ge1\}$, 
will be denoted $\frakm_v$.
The residue field of $v$, defined to be the field $\frako_v/\frakm_v$, will be denoted $k_v$. I will denote  the quotient homomorphism $\frako_v\to k_v$ by $\zeta_v$. I will set $h_v =h_{\zeta_v}:\ggle(\frako_v)\to\ggle( k_v)$ and $\overline h_v =\overline h_{\zeta_v}:\piggle(\frako_v)\to\piggle( k_v)$. (In particular, $\overline h_v(\gamma)$ is defined when $x\in \pizzle(\frako_v)\le\piggle(\frako_v)$. This special case will play an important role in the paper.)

\EndNumber

\begin{TheoremB}
There exists a natural number $N$ with the following property. 
Let  $\Gamma$ be any cocompact, discrete, torsion-free subgroup of $\pizzle(\CC) $.
Suppose that $\Gamma\le\pizzle(E)$, where $E$ is a number field.
Let $v$ be a valuation of $E$. Let $x$ and $y$ be non-commuting elements of $\Gamma$. Suppose that either
\Alternatives
\item $x$ does not lie in a $\piggle(E)$-conjugate of $\pizzle(\frako_v)$, or
\item $x\in\pizzle(\frako_v)$, the characteristic of $k_v$ is greater than $N$, and $\overline h_v(x)$ has order $7$ in $\pizzle(k_v)$. 
\EndAlternatives
Then for every point $P\in\HH^3$ we have
$$\max( d_P(x), d_P(y))>0.34.$$
\end{TheoremB}

This will be proved in the body of the paper as Theorem \ref{N lives}.

Like Theorem A, Theorem B represents a surprising interaction between the algebraic (or number-theoretic) aspect  of a hyperbolic manifold and its quantitative geometric aspect.

The proof of Theorem B is somewhat easier if one restricts attentions to groups $\Gamma$ that have ``integral traces'' in the sense that the traces of their elements are all algebraic integers, and I will outline the proof of this restricted version of the theorem first. In the setting of this version fo the theorem, Bass's $\ggle$ Subgroup Theorem \cite{bassgltwo} may be shown to imply that Alternative (i) of the hypothesis cannot hold. The proof in this case consists of an algebraic step and a geometric step. In the algebraic step one uses alternative (ii) of the hypothesis to construct a non-abelian subgroup of $\Theta$ of $\Gamma_1:=\langle x,y\rangle$ which is generated by a word of length $1$ and a word of length $7$ in $x$ and $y$, and whose index $|\Gamma_1:\Theta|$ is bounded below by the characteristic of $k_v$. Thus by making a suitable choice of the natural number $N$ in the statement of Theorem B, one can force $|\Gamma_1:\Theta|$ to be arbitrarily large. In the geometric step one shows that if two suitably short words generate a subgroup of suitably large index then one obtains a good lower bound for $\max( d_P(x), d_P(y))$. 

The algebraic step is carried out in Sections \ref{reversive section}---\ref{two-generator section}. The arguments here involve the $p$-nilpotency of characteristic-$p$ congruence kernels, the canonical involution of a two-generator Kleinian group, and a good deal of delicate finite group theory and linear algebra. The essential result for the later application is Corollary \ref{the daughter of rosie o'grady}.

This corollary, and Proposition \ref{do monkeys have uncles?} on which it depends, involve the assumption that for a certain element $x$ of a group $\Gamma\le\pizzle(\frako_v)$, the order $m$ of $\overline h_v(x)$ is at least $7$. 
This surprising appearance of the number $7$ is explained by Dickson's classification of subgroups of $\zzle(2,q)$, which is used in the paper via in Lemma \ref{basket driving}. In Dickson's classification the groups that appear as exceptional, such as $\cala_5$ and $\cals_4$, contain elements of order less than $7$; thus the assumption that $m\ge7$ forces $\overline h_v(\Gamma)$ to be non-exceptional, which turns out to be important for the group-theoretical arguments.

The geometric step needed to prove the restricted version of Theorem B is based on the techniques used to prove Theorem 4.2 and Corollary 4.3 of \cite {bounds}. The main ingredients are the strong $\log3$ Theorem (in the form proved in \cite{surgery}, improving the earlier version from \cite{accs}) and the theory of algebraic and geometric convergence of Kleinian groups. The refinements of the results of \cite{bounds} that are needed for this step are carried out in Section \ref{better bounds} of the present paper; the essential result is Corollary \ref{ranana}.

Because Corollary \ref{the daughter of rosie o'grady} is valid for any $m\ge7$, the arguments described above would work if the number $7$  in the statement of Theorem B were replaced by any larger integer. However, the resulting geometric estimates become weaker as $m$ increases.

To prove Theorem B in its unrestricted form, one also needs to consider the complementary case in which the group $\Gamma$ does not have integral traces. In this case it is a well-known consequence of Bass's $\ggle$ Subgroup Theorem that $M:=\HH^3/\Gamma$ is a Haken manifold, i.e. it contains an incompressible surface. For the case of a Haken manifold, techniques for estimating Margulis numbers, or more generally giving estimates of the type provided by Theorem B, were developed in \cite{hakenmarg}. While the results of \cite{hakenmarg}, as stated, are not strong enough to prove Theorem B, the general techniques do turn out to apply here. When the hypotheses of Theorem B hold and $\Gamma$ does not have integral traces, it can be shown, using the tree for $\zzle$ of a valued field and $3$-manifold topology,  that $M$ contains an incompressible surface of a restricted type. Using such a surface, the better estimates required for the conclusion of the theorem can be obtained by refining the arguments of \cite{hakenmarg}. These arguments occupy Sections \ref{pizzletree}, \ref{tweetwee}, and \ref{hakenmeasure}.

In Section \ref{bee section}, the above ingredients are combined to give the proof of Theorem B.

The arguments needed to deduce Theorem A from Theorem B occupy sections \ref{character variety section}---\ref{masterpiece section}. A crucial step involves applying Proposition 2.7 of \cite{finitistic}, which depends on a deep number-theoretical result due to Siegel and Mahler. This proposition implies that there is a finite subset 
$W$ of the field $K$ such that for every $x\in \pizzle(K)$ which is represented by an element of $\zzle(\CC)$ whose trace does not lie in $W$, one of the alternatives (i) or (ii) of Theorem B must hold. Using this fact and Theorem B, it is not hard to reduce the proof of Theorem A to showing that, up to conjugacy, there at most finitely pairs $( x, y)$ of elements of $\zzle$ such that $x$ and $y$ have prescribed traces, $\langle x,y\rangle $ is discrete, non-elementary and torsion-free and has trace field $K$, and $\max( d_P([x]), d_P([y]))\le0.34$ for some $P\in\HH^3$.

It turns out that a stronger statement is true, and it is stated in this paper as Proposition \ref{free analogue}. Up to conjugacy, the pairs $( x, y)$ of elements of $\zzle$ for which $\Gamma_1:=\langle x,y\rangle $ is discrete and non-elementary (or even for which $\Gamma_1$ is irreducible) are parametrized by certain points of the character variety of a rank-$2$ free group, which may be identified with $\CC^3$. Those pairs for which $x$ and $y$ have prescribed traces that are algebraic over $\QQ$ are parametrized by points of a line in $\CC^3$, which is in particular an affine curve defined over $\QQ$. Proposition \ref{free analogue} asserts that if $C\subset\CC^3$ is any curve defined over $\overline\QQ$, if $D$ is any positive integer, and if $\alpha$ is any positive number less than $\log3$, then there are at most finitely many points of $C$ parametrizing pairs $(x,y)$ such that (1) $\Gamma_1:=\langle x,y\rangle$ is discrete, non-elementary and torsion-free, (2) $\max(d_P([x]),d_P([y]))\le\alpha$ for some $P\in\HH^3$, and (3) the trace field of $\Gamma_1$ has degree at most $D$.  The proof of Proposition \ref{free analogue} involves a result, Proposition \ref{long reid type prop}, which is a partial generalization of a result due to Long and Reid \cite[Theorem 3.2]{long-reid}, and the proof closely parallels the proof of their result.

I am grateful to Marc Culler, Jason DeBlois, Ben McReynolds, Alan Reid and Steven Smith for helping to make this a better paper. Smith helped me with the finite group theory needed in Sections \ref{Delta section} and \ref{two-generator section}. The transition from Theorem B to Theorem A made possible by Reid's telling me about the result of Siegel and Mahler which are quoted in \cite{finitistic}, and by a suggestion from DeBlois that led to the character variety arguments of Sections \ref{character variety section}---\ref{masterpiece section}.

\section{Preliminaries}

\Number
If $V$ is a vector space over a field, I will denote the group of linear automorphisms of $V$ by ${\rm GL}(V)$. As usual, if $R$ is a commutative ring, $\ggle(R)$ denotes the group of invertible $2\times2$ matrices over $R$, and $\zzle(R)$ is the subgroup consisting of matrices of determinant $1$, while $\piggle(R)$ and $\pizzle(R)$ denote the quotients of $\ggle(R)$ and $\zzle(R)$ by their centers, which consist of scalar matrices. The full ring of $2\times2$ matrices with entries in $R$ will be denoted $\calm_2(R)$.

By a {\it representation} of a group $\Gamma$ in $G$, where $G$ is one of the groups ${\rm GL}(V)$, $\ggle(R)$, $\zzle(R)$, $\piggle(R)$ or $\pizzle(R)$, I will mean simply a homomorphism from $\Gamma$ to $G$.
\EndNumber

\Number

{\it Hyperbolic manifolds} will be understood to be complete and connected.
Up to isometry, every orientable hyperbolic $3$-manifold has the form $\HH^3/\Gamma$, where $\Gamma$ is a discrete, torsion-free subgroup of $\pizzle(\CC)$, uniquely determined by $M$ up to conjugacy. When an
orientable hyperbolic $3$-manifold is written in the form $\HH^3/\Gamma$, it will be understood that $\Gamma$ is a discrete, torsion-free subgroup of $\pizzle(\CC)$ (and is cocompact if and only if $M$ is closed).
\EndNumber

\Number\label{watson}
I will use $1$ as the default notation for the identity element of a group, but I will denote the identity matrix in $\ggle(K)$ by $I$. Thus we have $[I]=1\in\piggle(K)$.

Recall that a group $\Gamma$ is said to be {\it elementary} if it has an abelian subgroup of finite index.

Note that if $\tGamma\le\zzle(\CC)$ is torsion-free then $\Pi_\CC$ maps $\tGamma$ isomorphically onto a subgroup $\Gamma$ of $\pizzle(\CC)$, and that $\Gamma$ is discrete only if $\tGamma$ is. Since $\Gamma$ and $\tGamma$ are isomorphic, $\Gamma$ is torsion-free, or non-elementary, if and only if $\tGamma$ is torsion-free, or non-elementary respectively.
\EndNumber

The following result will be needed at a couple of points in the paper. 

\Proposition\label{just conjugate}
If $\Gamma$ is any cocompact, discrete, torsion-free subgroup of $\pizzle(\CC) $, then $\Gamma$ is conjugate to a subgroup of $\pizzle(E)$ for some algebraic number field $E$.
\EndProposition

\Proof
Since $\Gamma$ is in particular a non-elementary, finitely generated, torsion-free Kleinian group, it follows from \cite[Proposition 4.2]{finitistic} that $\Gamma$ is isomorphic to a Kleinian group $\Gamma_1$ such that $\Gamma_1\le\pizzle(E)$ for some number field $E$. Since $\Gamma$ is cocompact, it follows from the Mostow rigidity theorem that $\Gamma$ is in fact conjugate to $\Gamma_1$ in $\pizzle(\CC)$.
\EndProof

\Number\label{adams spectral nahant} If $\Gamma$ is a group and $p$ is a prime, the {\it mod $p$ commutator subgroup} of $\Gamma$ is defined to be the subgroup of $\Gamma$  generated by all commutators and $p$-th powers.
\EndNumber

\section{Reversive groups}\label{reversive section}
Let $\Gamma$ be a group. I will define a {\it reversive system} for $\Gamma$ to be a pair of elements $(x,y)$ such that (1) $x$ and $y$ generate $\Gamma$, and (2) there is an automorphism $J$ of $\Gamma$ such that $J(x)=x^{-1}$ and $J(y)=y^{-1}$. Note that such an automorphism $J$ must be an involution, and that it is uniquely determined by $(x,y)$. I will call it the {\it reversal} defined by $(x,y)$. The group $\Gamma$ will
be said to be {\it reversive} if it admits a reversive system.

Recall that a subgroup $\Gamma$ of $\zzle(K)$, where $K$ is a field, is termed {\it reducible} if there is a $1$-dimensional subspace of $K^2$ that is invariant under $\Gamma$; or equivalently if $\Gamma$ is conjugate in $\ggle(K)$ to a group of upper-triangular matrices. A subgroup $\Gamma$ of $\pizzle(K)$ is said to be reducible if $\Pi_K^{-1}(\Gamma)\le\zzle(K)$ is reducible.

\Proposition\label{it's so diversible it's reversible}If\, $\Gamma$ is an irreducible subgroup of $\pizzle(\CC)$, and if $x$ and $y$ are elements that generate $\Gamma$, then $(x,y)$ is a reversive system. In particular, every two-generator irreducible subgroup of $\pizzle(\CC)$ is reversive.
\EndProposition

\Proof
Let $X$ and $Y$ be matrices in $\zzle(\CC)$ with $[X]=x$ and $[Y]=y$. Since $\Gamma$ is irreducible, we have $\trace XYX^{-1}Y^{-1}\ne2$ (see \cite[Proposition 1.5.5]{varreps}). Set $C=XY-YX$. We have $\trace C=0$, and $\det C=\det(XYX^{-1}Y^{-1}-I)\ne0$, since $XYX^{-1}Y^{-1}$ does not admit $1$ as an eigenvalue. Hence $C^2=\alpha I$ for some $\alpha\in\CC-\{0\}$. We have 
$$\trace X^{-1}C=\trace(Y-X^{-1}YX)=\trace Y-\trace(X^{-1}YX)=0.$$
Hence $(X^{-1}C)^2= \beta I$ for some $\beta\in\CC-\{0\}$. It follows that $CX^{-1}C^{-1}=\alpha^{-1} CX^{-1}C=\alpha^{-1}\beta X$. Since $\det(CX^{-1}C^{-1})=1=\det X$ we must have $\alpha^{-1}\beta=\pm1$, so that $CX^{-1}C^{-1}=\pm X$. The same argument shows that $CY^{-1}C^{-1}=\pm Y$. It follows that there is an automorphism $J$ of $\Gamma$ defined by $J([Z])=[C^{-1}ZC]$, and that $J(x)=x^{-1}$ and $J(y)=y^{-1}$.
\EndProof

\Lemma\label{and it's gonna stay that way}
Suppose that $(x,y)$ is a reversive system for a group $\Gamma$, and let $J$ be the reversal defined by $(x,y)$. Let $N$ be a normal subgroup of $\Gamma$ such that $J(N)=N$. Then $\Gamma/N$ is a reversive group.
\EndLemma

\Proof
Set $\overline\Gamma=\Gamma/N$, and let $\pi:\Gamma\to\overline\Gamma$ denote the quotient homomorphism. Set $\overline x=\pi(x)$ and $\overline y=\pi(y)$, so that $\overline\Gamma=\langle\overline x,\overline y\rangle$. Since
$J(N)=N$, there is an automorphism $\overline J$ of $\overline\Gamma$ such that $\pi\circ J=\overline J\circ \pi$. It is clear that $J(\overline x)=\overline x^{-1}$ and $J(\overline y)=\overline y^{-1}$, so that $(\overline x,\overline y)$ is a reversive system for $\overline\Gamma$.
\EndProof

\Definitions
Recall that if $p$ is a prime, an {\it elementary abelian $p$-group} is a finite abelian group in which of every non-trivial element has order $p$. An elementary abelian $p$-group $S$ may be regarded as the additive group of an $\FF_p$-vector space, and every automorphism of the group $S$ is a linear automorphism of this vector space. The {\it rank} of $S$ is its dimension as an $\FF_p$-vector space. An automorphism $\phi$ of $S$ will be termed {\it unimodular} if the determinant of $\phi$, regarded as a linear automorphism, has determinant $1$. Likewise, $\phi$  will be termed {\it semi-unimodular} if its determinant is $\pm1$.
\EndDefinitions

\Number\label{harris pat}
If $S$ is an elementary abelian $p$-group for some prime $p$, then every element of the commutator subgroup of $\Aut S$ is unimodular. This is because if $V$ is a finite-dimensional vector space over a field, the commutator subgroup of ${\rm GL}(V)$ consists of matrices of determinant $1$.
\EndNumber

\Number\label{paris hat}
If $S$ is an elementary abelian $p$-group for some prime $p$, then every odd-order semi-unimodular element of $\Aut S$ is unimodular. This is because if $V$ is a finite-dimensional vector space over a field, an element of odd finite order in ${\rm GL}(V)$ cannot have determinant $-1$.
\EndNumber

\Proposition\label{unimodular}
Let $\Delta$ be a reversive finite group, let $p$ be an odd prime, and let $S\triangleleft \Delta$ be an elementary abelian $p$-group. Suppose that $\Delta/S$ is abelian. Then for any $\xi\in \Delta$, the automorphism $s\mapsto \xi s\xi^{-1}$ of $S$ is semi-unimodular.
\EndProposition

\Proof
Let $(x,y)$ be a reversive system for $\Delta$, and let $J$ denote the reversal defined by $(x,y)$. Then $J(x)=x^{-1}$ and $J(y)=y^{-1}$. If we denote by $\Delta_1$ the abelian group $\Delta/[\Delta,\Delta]$, and if we denote by $ x_1,y_1\in\Delta_1$  the images of $x$ and $y$ under the quotient homomorphism, and by  $J_1$ the automorphism of $\Delta_1$ induced by $J$, then we have $J_1(x_1)=-x_1$ and $J_1(y_1)=-y_1$. Hence $J_1(z)=-z$ for every $z\in\Delta_1$. It follows that every subgroup containing $\Delta_1$ is $J$-invariant. Since $\Delta/S$ is abelian, in particular $S$ is $J$-invariant. Let $j$ denote the automorphism $J|S$ of $S$. We may regard $j$ as an element of ${\rm GL}(S)$.

Now regard $S$ as the additive group of a vector space over $\FF_p$, and define a representation $\rho$ of $\Delta$ in ${\rm GL}(S)$ by setting $\rho(\xi)(s)=\xi s\xi^{-1}$ for all $\xi\in\Delta$, $s\in S$. We need to show that $\det\rho(\xi)=\pm1$ for every $\xi\in\Delta$. Since $x$ and $y$ generate $\delta$, it suffices to show this for $\xi=x$ and for $\xi=y$. By symmetry we need only consider $\xi=x$.

Set $\phi=\rho(x)$. For each $s\in S$, we have $xsx^{-1}=\phi(x)$. Since $J(x)=x^{-1}$, it follows that
$$x^{-1}j(s)x=J(xsx^{-1})=J(\phi(s)) =j(\phi(s)).$$
Hence
$$j(s)=xj(\phi(s))x^{-1}=\phi(j(\phi(s))).$$
This shows that $j=\phi\circ j\circ\phi$, so that $j^{-1}\phi j=\phi^{-1}$. Thus $\phi$ and $\phi^{-1}$ are conjugate in ${\rm GL}(S)$, and hence $\det\phi=\pm1$, as required.
\EndProof

\section{Residually $p$-nilpotent groups}\label{resnil section}

\Definition\label{i got st vitus dancers}
Let $p$ be a prime number. A group $\Gamma$ is said to be {\it residually $p$-nilpotent} if the trivial subgroup of $\Gamma$ is the intersection of a (possibly infinite) collection $\caln$ of normal subgroups such that for each $N\in\caln$ the index $|\Gamma:N|$ is a power of $p$.
\EndDefinition

\Proposition\label{rheingold is as good to your taste}
Let $v$ be a valuation of a number field $E$, and suppose that the characteristic $p$ of $ k_v$ is odd. Then the kernel of
$\overline h_v:\piggle(\frako_v)\to\piggle(k_v)$ is residually $p$-nilpotent,
where
\EndProposition

The following proof is similar to the discussion on p. 87 of
\cite{residual}.

\Proof[Proof of Proposition \ref{rheingold is as good to your taste}]
Fix a generator $\pi$ for $\frakm_v$. For each 
$n\ge1$, let $\zeta_{v,n}:\frako_{v}\to\frako_v/(\pi^n)$ denote the quotient homomorphism, and set
$\overline h_{v,n}=\overline h_{\zeta_{v,n}}:\piggle(\frako_{v})\to\piggle(\frako_v/(\pi^n))$ (see \ref{bet you can't eat one}).
Then we have $\zeta_{v,1}=\zeta_v$ and $\overline h_{v,1}=\overline h_v$. For each $n$ set $\Gamma _n=\ker \overline h_{v,n}$. We are required to show that $\Gamma _1$ is residually $p$-nilpotent.

It follows from the definition that for each $n\ge1$ we have $\Gamma _n\triangleleft \Gamma$; in particular $\Gamma_n\triangleleft\Gamma_1$. I claim:
\Claim\label{why don't you shut up}$|\Gamma_1/\Gamma_n|$ is a power of $p$ for each $n$.
\EndClaim
To prove \ref{why don't you shut up}, let $Z_{v,n}$ denote the 
natural homomorphism from
$\frako_v/(\pi^n)$ to $k_{v}$, so that $Z_{v,n}\circ\zeta_v =\zeta_{v,n}$.
If we set $H_{v,n}:=\overline h_{Z_{v,n}}$, we have $H_{v,n}\circ\overline h_v =\overline h_{v,n}$, and hence
$\Gamma _1/\Gamma _n$ is isomorphic to $K_n:=\ker H_{v,n}$.

Set $M^{(n)}=M_2(\frako_v/(\pi^n))$. For each natural number $n$, and for each $A\in M^{(n)}$, we have $\det (I+\pi A)\equiv1\pmod{\pi}$.
Since $\frako_v/(\pi^n)$ is a local ring
with maximal ideal $\pi \frako_v/(\pi^n)$ it follows that $\det (I+\pi A)$ is
invertible in $\frako_v/(\pi^n)$, and hence that $I+\pi A\in
\ggle(\frako_v/(\pi^n))$. We may therefore
define a map of sets $\phi_n:M^{(n)}\to \piggle(\frako_v/(\pi^n))$ by $\phi_n(A)=[I+\pi A]$. 

It is immediate from the definitions that $\phi_n(A)\in K_n$ for every $A\in M^{(n)}$. The definition of $K_n$ also implies that any element of $K_n$ may be written in the form $[\pm I+\pi A]$ for some $A\in M^{(n)}$. Since $[- I+\pi A]=[I+\pi(-A)]$, we may in fact write any element of $K_n$ in the form $[I+\pi A]$ for some $A\in M^{(n)}$. This shows that $\phi_n$ maps $M^{(n)}$ onto $K_n$.

If $A,A'\in M^{(n)}$ and $\phi_n(A)=\phi_n(A')$, then $I+\pi A=\pm(I+\pi A')\in \frako_v/(\pi^n)$. Hence if $A\ne A'$ we must have $I+\pi A=-(I+\pi A')$, so that $\pi(A+A')=-2I$. It follows that $\pi$ divides $2$ in  $\frako_v/(\pi^n)$, which is impossible since $p$ is odd. This shows that $\phi_n$ maps $M^{(n)}$ bijectively onto $K_n$.
Hence $|K_n|=|\pi M^{(n)}|$. 

The local ring $\frako_v/(\pi^n)$ has residue field $k_v$ of
characteristic $p$; hence $|\frako_v/(\pi^n)|$ is a power of $p$. Since
$\pi M^{(n)}$ is a submodule of $M^{(n)}$, which is a rank-four free module over $\frako_v/(\pi^n)$, the order of $\pi M^{(n)}$ is also a power of
$p$. Since $|\Gamma _1/\Gamma _n|=|K_n|=|\pi M^{(n)}|$, this proves (\ref{why don't you shut up}).

Next I claim that
\Claim\label{nuthin left}
$\bigcap_{n=1}^\infty \Gamma _n=\{1\}$.
\EndClaim
To prove this, suppose that $g$ is an element of $\bigcap_{n=1}^\infty \Gamma _n $. Let us write $g=[B]$ for some $B\in\ggle(E)$. Then for every natural number $n$ we have $B\equiv\pm I\pmod{ \pi^n}$. Hence either there are infinitely many $n\in\NN$ for which $B\equiv I\pmod{ \pi^n}$, or
there are infinitely many $n\in\NN$ for which $B\equiv-I\pmod{ \pi^n}$. After changing the sign of $B$ if necessary we may assume that  $B\equiv I\pmod{ \pi^n}$ for every $n\in S$, where $S\subset \NN$ is some infinite set. But if $B\equiv I\pmod{ \pi^n}$ for a given $n$, then $B\equiv I\pmod{ \pi^m}$ for every $m<n$. Hence $B\equiv I\pmod{ \pi^n}$ for every $n\ge1$. Since $\bigcap_{n=1}^\infty(p^n)=\{0\}$, it follows that $B=I$. This proves \ref{nuthin left}.

Since the family of normal subgroups $(\Gamma_n)_{n\ge1}$ has Properties \ref{why don't you shut up} and \ref{nuthin left}, it follows from  the definition that  $\Gamma _1$ is residually $p$-nilpotent.
\EndProof

\Proposition \label{fantabula}
If $p$ is a prime and $X$ is a residually $p$-nilpotent group, then every subgroup of  $X$ is residually $p$-nilpotent.
\EndProposition

\Proof
Let  $\caln$ be a collection of normal subgroups of $X$ such that $\bigcap_{N\in\caln}N=\{1\}$, and such that $|\Gamma:N|$ is a power of $p$ for each $N\in\caln$. Then for any subgroup $Y$ of $X$ wee have in particular that $\bigcap_{N\in\caln}(N\cap Y)=\{1\}$. Furthermore, for each $N\in\caln$, since $N\triangleleft\Gamma$, we have $N\cap Y\triangleleft Y $; and $|Y:Y\cap N|$ divides $|\Gamma:N|$, and is therefore a power of $p$.
\EndProof

\Proposition \label{nabla}
Let $Y$ be a finite-index normal subgroup of a group $X$ and let $p$ be a prime. If $Y$ is residually $p$-nilpotent and $X/Y$ is a $p$-group, then $X$ is residually $p$-nilpotent.
\EndProposition

\Proof
Let  $\caln$ be a collection of normal subgroups of $Y$ such that $\bigcap_{N\in\caln}N=\{1\}$, and such that $|\Gamma:N|$ is a power of $p$ for each $N\in\caln$.  For any given $N\in\caln$ we have $|X:N|<\infty$, and hence $N$ has only finitely many conjugates in $X$, say $N=N_1,\ldots, N_n$. Set  $C_N=N_1\cap\ldots\cap N_n$. It is immediate from the definition of $C_N$ that $C_N\triangleleft X$ and that $C_N\le N$. Since $Y\triangleleft X$ we have $N_i\le Y$ for $i=1,\ldots,n$; furthermore, each $N_i$ is the image of $N$ under an automorphism of $Y$, so that $N_i\triangleleft Y$, and $|Y:N_i|$ is a power of $p$. This implies that $|Y:C_N|$ is also a power of $p$. As $|X:Y|$ is a power of $p$, it follows that $|X:C_N|$ is also a power of $p$. Since $C_N\le N$ and $\bigcap_{N\in\caln}N=\{1\}$, we have $\bigcap_{N\in\caln}C_N=\{1\}$. This shows that $X$ is residually $p$-nilpotent.
\EndProof

\Proposition \label{squash}
Let $p$ a prime, and let $X$ be a finitely generated, residually $p$-nilpotent group which is not cyclic. Let $X_1$ denote the mod $p$ commutator subgroup (\ref{adams spectral nahant}) of $X$. Then $X/X_1$ is an elementary abelian $p$-group of rank at least $2$.
\EndProposition

\Proof
It follows from the definition of $X_1$ that $X/X_1$ is an elementary abelian $p$-group. I will show that it has rank at least $2$.

First consider the case in which the finitely generated group $X$ is abelian. 
If $X/X_1$ has rank at most $1$, then $X$ is a sum of direct product of two subgroups $A$ and $B$, where $A$ is cyclic  and  $|A|$ is either infinite or a power of $p$, while $|B|$ is finite and prime to $p$. Thus $B$ is contained in every subgroup of $X$ whose index is a power of $p$. Since $X$ is residually $p$-nilpotent we have $B=\{1\}$; hence $X$ is cyclic, a contradiction to the hypothesis.

Now consider the case in which $X$ is non-abelian. Let $x$ and $y$ be non-commuting elements of $X$. Since $X$ is residually $p$-nilpotent  and $[x,y]\ne1$,  there is a normal subgroup $N$ of $X$ such that $[x,y]\notin N$, and  $| X:N|$ is a power of $p$. Hence $G:= X/N$ is a finite $p$-group, and in particular it is nilpotent. On the other hand, the images of $x$ and $y$ under the quotient homomorphism do not commute in $G$, and hence $G$ is non-abelian. Thus if $G=G_1\triangleleft\cdots\triangleleft G_n=\{1\}$ denotes the lower central series of $G$, then $G/G_3$ is a non-abelian group having $G_2/G_3$ as a central subgroup, and hence $G/G_2$ is non-cyclic. It follows that $G$ has a quotient $A$ which is an elementary abelian $p$-group of rank at least $2$. In particular, $A$ is isomorphic to a quotient of $ X$, and the conclusion follows in this case.
\EndProof

\Lemma\label{i shoulda stood in bed}
Let $\Gamma$ be a group admitting a reversive system $(x,y)$. Let $p$ be a prime, and let $m$ be a natural number. Suppose that there exists a homomorphism $\ell $ of $\Gamma$ onto a finite abelian group $A$ such that $\ell(x) $ has order $m$ in $A$, and $\ker \ell $ is $p$-residually nilpotent. Then   there exist a finite group $\Delta $ and a homomorphism $\lambda:\Gamma\to \Delta $ such that 
\begin{enumerate}
\item  $\Delta $ is reversive;
\item $\Delta$ has a normal $p$-Sylow subgroup $S$ which is an elementary abelian $p$-group of rank at least $2$
\item $\Delta/S$ is abelian; and
\item the image of $\lambda(x)$ under the quotient homomorphism $\Delta\to\Delta/S$ has order $m$ in $\Delta/S $.
\end{enumerate}
\EndLemma

\Proof
First consider the special case in which the order of $A$ is prime to $p$. In this case, let $n$ denote the order of $\ell(y)$. Since $n$ and $m$ divide $A$ they are prime to $p$. Set $K=\ker\ell$, and let $K_0$ denote the subgroup of $\Gamma$ generated by $x^m$, $y^n$ and the commutator subgroup $[\Gamma,\Gamma]$ of $\Gamma$. Since $[\Gamma,\Gamma]\le K_0$ we have $K_0\triangleleft\Gamma$, and $A_0:=\Gamma/K_0$ is abelian. Furthermore, $A_0$ is generated by the images $\overline x$ and $\overline y$ of $x$ and $y$ under the quotient homomorphism; the definition of $K_0$ implies that $\overline x^m=\overline y^n=1$, and hence the order of $A_0$ is finite and prime to $p$. In particular $K_0$ has finite index in $\Gamma$ and is therefore finitely generated.

 Since $G_1$ is abelian and $\ell(x^m)=\ell(y^n)=1$, we have $K_0\le K$.

Let $L$ denote the mod $p$ commutator subgroup (\ref{adams spectral nahant}) of $K_0$. Since $L$ is a characteristic subgroup of $K_0\triangleleft\Gamma$, we have $L\triangleleft\Gamma$. Define $\Delta$ to be the quotient group $\Gamma/L$, and define $\lambda:\Gamma\to\Delta$ to be the quotient homomorphism. I will show that $\Delta$ is finite and that conclusions (1)--(4) hold.

Since $K_0\le K$, and $K$ is residually $p$-nilpotent by hypothesis, $K_0$ is residually $p$-nilpotent by Proposition \ref{fantabula}. If $K_0$ were cyclic, then since  $A_0 =\Gamma/K_0$ is abelian, $\Gamma$ would be solvable, a contradiction to the hypothesis. Hence $K_0$ is non-cyclic, and it follows from Proposition \ref{squash}  that $S:=K_0/L$ is an elementary abelian $p$-group of rank at least $2$. Since $K_0\triangleleft\Gamma$, we have $S=K_0/L\triangleleft\Gamma/L=\Delta$. Furthermore, we have $\Delta/S\cong\Gamma/K_0=A_0$, so that $\Delta/S$ is abelian and its order is finite and prime to $p$. This shows that $\Delta$ is finite and that $S$ is a $p$-Sylow subgroup. Thus (2) and (3) are established.

To prove (1), let $J$ denote the reversal defined by $(x,y)$. Since $J(x)=x^{-1}$ and $J(y)=y^{-1}$, we have $J(x^n)=x^{-n}$ and $J(y^n)=y^{-n}$. Since $K_0$ is generated by $x^n$ and $y^n$ and the characteristic subgroup $[\Gamma,\Gamma]$ of $\Gamma$, it follows that $K_0$ is invariant under $J$. But $L$ is a characteristic subgroup of $K_0$ and is therefore also invariant under $J$. It now follows from Lemma \ref{and it's gonna stay that way} that $\Delta=\Gamma/L$ is reversive; this is Assertion (1).

To prove (4) we must show that $\pi\lambda(x)$ has order $m$ in $\Delta/S $, where $\pi:\Delta\to\Delta/S$ denotes the quotient homomorphism. We have $\ker(\pi\lambda)=\lambda^{-1}(S)=\lambda^{-1}(K_0/L)=K_0$. Since $x^m\in K_0$ by the definition of $K_0$, the order of $\pi\lambda(x)$ is at most $m$. On the other hand, since $K_0\ge K=\ker\ell$, the order of $\pi\lambda(x)$ is at least the order of $\ell(x)$, which by hypothesis is equal to $m$. Thus (4) is established. \EndProof

\section{Finite groups}\label{Delta section}

The purpose of this section and the next is to prove Corollary \ref{the daughter of rosie o'grady}, which was discussed in the Introduction. This section consists of technical background on finite groups. The main result of the section is Lemma \ref{i shoulda stood in bed}.

\Lemma\label{you talk about teenage werewolves}
Let $p$ be a prime, and let $T$ be a $p$-Sylow subgroup of a finite group $X$. Suppose that $T$ is central in $X$ and that $X/T$ is abelian. Then $X$ is abelian.
\EndLemma

\Proof
Since $T$ is central and $X/T$ is abelian, $X$ is nilpotent, and is therefore the direct product of its Sylow subgroups. Hence it suffices to prove that the Sylow subgroups of $X$ are abelian. Since $T$ is central it is abelian. If $S$ is an $\ell$-Sylow subgroup of $X$ for some prime $\ell\ne p$, then $S\cap T=\{1\}$, and hence the quotient homomorphism $X\to X/T$ maps $S$ injectively into $X/T$. Since $X/T$ is abelian, so is $S$. 
\EndProof

\Lemma\label{wide stance}
Let $Q$ be a finite group, let $p$ be a prime, and let $m$ be a natural number. Assume that $Q$ has a normal $p$-Sylow subgroup $T$ which is an elementary abelian $p$-group. Let $ w$ be an element of $Q$ such that either $C( w)\supset T$ or $C( w)\cap T=\{1\}$. Let $\overline w$ denote the image of $ w$ under the quotient homomorphism $Q\to Q/T$, and suppose that $\overline w$ has  order $m$ in $Q/T$. Assume that the normal closure of $\overline w$ in $Q/T$ is abelian. Then either $ w$ has order $m$ in $Q$, or the normal closure of $\{ w\}$ in $Q$ is abelian.
\EndLemma

\Proof
First suppose that $C( w)\cap T=\{1\}$. Note that since $\overline w$ has order $m$ in $Q/T$, we have $ w ^m\in T$. Since $ w ^m$ commutes with $ w $ it follows that $ w ^m \in C( w)\cap T$ and hence that $ w ^m=1$. Since $\overline w$ has order $m$ in $Q/T$, we have $ w^k\notin T$, and hence $ w^k\ne1$, whenever $0<k<m$. Hence $ w$ has order $m$ in this case.

Now suppose that $C( w)\supset T$. Then $w\in C(T)$. Since $T$ is abelian we have $T\le C(T)$, and hence $\{ w\}\cup T\le C(T)$. Since $T$ is normal in $Q$, its centralizer $C(T)$ is also normal. Hence if $X$ denotes the normal closure of $\{ w\}\cup T$, we have $X\le C(T)$; that is, $T$ is a central subgroup of $X$. On the other hand, the definition of $X$ implies that $X/T$ is the normal closure of $\overline w$ in $Q/T$, which by hypothesis is abelian. Finally, since $T$ is the Sylow subgroup of $Q$, it is also the  Sylow subgroup of $X$. Thus $X$ and $T$ satisfy all the hypotheses of Lemma \ref{you talk about teenage werewolves}, and hence $X$ is abelian. Since the normal closure of $\{ w\}\cup T$ is abelian, in particular, the normal closure of $\{ w\}$ is abelian.
\EndProof

\Lemma\label{betwixt and betwining}Let $z$ be an element of a finite group $H$ such that $\langle z\rangle$ is normal in $H$. Let $K$ be a finite field whose characteristic does not divide $H$. Let $V$ be a vector space over $K$ such that $2\le\dim V<\infty$. Let $\rho:H\to {\rm  GL}(V)$ be a representation such that $\det\rho(z)=1$. Then there exists a direct sum decomposition $V=V_0\oplus V_1$ such that
\begin{itemize}
\item each of the $V_i$ is invariant under $\rho(H)$;
\item $\dim V_1\le 2$; and
\item either $\rho(z)|V_1$ is the identity, or $\rho(z)(x)\ne x$ for every non-zero vector $x\in V_1$.
\end{itemize}
\EndLemma

\Proof
Set $\phi=\rho(z)$. Let $W$ denote the subspace of $V$ consisting of all vectors $x\in V$ such that $\phi(x)=x$. For any $x\in W$ and any $h\in H$, we have $h^{-1}zh\in\langle z\rangle$ since $\langle z\rangle$ is normal; hence $x=\rho(h^{-1}zh)(x)=\rho(h^{-1})\phi \rho(h)(x)$, and so $\rho(h)(x)=\phi \rho(h)(x)$. This shows that $\phi(x)\in W$. Hence $W$ is invariant under $\rho(H)$. 

Now since $H$ is prime to the characteristic $p$ of $K$, it follows from Maschke's Theorem (see \cite[Theorem 10.8]{curtis-reiner}) that the $\rho(H)$-invariant subspace $W$ has a $\rho(H)$ invariant linear complement $W'$. By definition this means that $V=W\oplus W'$ and that $\rho(H)(V)=V$. I claim:
\Claim\label{not small at all} Either $\dim W\ge2$ or $\dim W'\ge2$.
\EndClaim
Since $\dim W+\dim W'=\dim V\ge2$, \ref{not small at all} is immediate unless $\dim W=\dim W'=1$. If the latter situation obtains, let $e$ and $e'$ be non-zero vectors in $W$ and $W'$ respectively. Then $e$ and $e'$ form a basis of $V$. Since $e\in W$ we have $\phi(e)=e$. Since $W'$ is $\rho(H)$-invariant, we have $\phi(e')=\alpha e'$ for some $\alpha\in K$. In the basis $(e,e')$, the matrix of $\phi$ is $\begin{pmatrix}1&0\cr0&\alpha\end{pmatrix}$. In view of the hypothesis we have $\alpha=\det\phi=1$. Hence $\phi$ is the identity, so that $W=V$, a contradiction to the assumption $\dim W=1$. Thus \ref{not small at all} is proved.

Now define $V_0=W'$ and $V_1=W$ if $\dim W\ge2$, and otherwise define $V_0=W$ and $V_1=W'$. It follows from \ref{not small at all} that in either case we have $\dim V_1\ge2$. Since $V=W\oplus W'$, and since $W$ and $W'$ are $\rho(H)$-invariant, we have $V=V_0\oplus V_1$, and $W$ and $W'$ are $\rho(H)$-invariant. If $V_1=W$ then
$\rho(z)|V_1$ is the identity by the definition of $W$. If $V_1=W'$ then $V_1\cap W=\{0\}$,  and hence $\rho(z)(x)\ne x$ for every non-zero vector $x\in V_1$.
\EndProof

\Lemma\label{happenstance}
Let $\Delta$ be a finite group and let $p$ be a prime, and let $m$ be a natural number. Assume that $\Delta$ has a normal $p$-Sylow subgroup $S$ which is an elementary abelian $p$-group of rank at least $2$. Let $\xi$ be an element of $\Delta-S$, and suppose that the automorphism $s\mapsto \xi s\xi^{-1}$ of $S$ is unimodular. Let $\overline\xi$ denote the image of $\xi$ under the quotient homomorphism $\Delta\to \Delta/S$, and suppose that $\overline \xi$ has  order $m$ in $\Delta/S$. Assume that $\langle\overline\xi\rangle$ is normal in $\Delta/S$. Then there exist a group $Q$ and a surjective homomorphism $\sigma:\Delta\to Q$ such that 
\begin{itemize}
\item every cyclic subgroup of $Q$ has index at least $p$, and
\item either $\sigma(\xi)$ has order $m$ in $Q$, or the normal closure of $\sigma(\xi)$ in $Q$ is abelian.
\end{itemize}
\EndLemma

\Proof
I will identify $S$ with the additive group of a vector space over $\FF_p$. Set $H=\Delta/S$, and let $\pi:\Delta\to H$ denote the quotient homomorphism. Since $S$ is abelian, the action of $\Delta$ on $S$ by conjugation factors through an action of $H$. Hence there is a representation $\rho:H\to {\rm  GL}(S)$  such that 
$\rho(\pi(\delta))=\delta s\delta^{-1}$ for all $s\in S$ and $\delta\in\Delta$.
By hypothesis we have $\dim S\ge2$ and $\det\rho(\overline\xi)=1$. 
Note also that $|H|$ is prime to $p$ since $S$ is a $p$-Sylow subgroup of $\Delta$.

Thus all the hypotheses of Lemma \ref{betwixt and betwining} hold if we take $K=\FF_p$, $V=S$ and $z=\overline\xi$, and define $H$ and $\rho$ as above. Hence Lemma \ref{betwixt and betwining}  gives
a direct sum decomposition $S=S_0\oplus S_1$ such that
\begin{enumerate}
\item each of the $S_i$ is invariant under $\rho(H)$;
\item $\dim S_1\le 2$; and
\item either $\rho(\overline\xi)|S_1$ is the identity, or $\rho(\overline\xi)(v)\ne v$ for every non-zero vector $v\in S_1$.
\end{enumerate}
We may regard $S_0$ and $S_1$ as subgroups of the elementary abelian $p$-group $S$, and in group-theoretical notation we have
$S=S_0\times S_1$.

In view of the definition of $\rho$, it follows from (1) that $S_0$ and $S_1$ are normal in $\Delta$. Define $Q=\Delta/S_0$, and define $\sigma:\Delta\to Q$ to be the quotient homomorphism.  I will show that the conclusions of the present lemma hold with these choices of $Q$ and $\sigma$.

Since $S=S_0\times S_1$, the homomorphism $\sigma$ maps $S_1$ isomorphically onto a subgroup $T$ of $Q$, and we have $\sigma^{-1}(T)=S$. Since $T$ is isomorphic to the additive group of the vector space $V$, it follows from (2) above that $T$ is an elementary abelian $p$-group of rank $2$. Since $\sigma$ is surjective and $\sigma^{-1}(T)=S$, the subgroup $T$ is normal in $Q$,  and $Q/T\cong H$. In particular, the order of $Q/T$ is prime to $p$, and so $T$ is a $p$-Sylow subgroup of $Q$.

If $J$ is any cyclic subgroup of $Q$, then $J\cap T$ is a cyclic subgroup of $T$. Since $T$ is an elementary abelian $p$-group, we have $p\le|T:J\cap T|\le|Q:J|$. This proves that every cyclic subgroup of $Q$ has index at least $p$, which is the first conclusion of the lemma.

It remains to prove that either $\sigma(\xi)$ has order $m$ in $Q$, or the normal closure of $\sigma(\xi)$ in $Q$ is abelian. This will follow from Lemma  \ref{wide stance}, provided that the hypotheses of Lemma \ref{wide stance} hold when we set $w=\sigma(\xi)$. I have already shown that  $T$ is an elementary abelian $p$-group of rank $2$, that is a $p$-Sylow subgroup of $Q$, and that it is normal in $Q$. I will complete the proof by verifying that the remaining hypotheses of Lemma \ref{wide stance} hold in the present setting.

The definitions of $\overline\xi$ and $\rho$ imply that the fixed subspace of $\rho(\overline\xi)$ is $C(\xi)\cap S$. It therefore follows from (3) that either $C(\xi)\supset S_1$ or $C(\xi)\cap S_1=\{1\}$. If $C(\xi)\supset S_1$, then since $\sigma$ maps $S_1$ onto $T$, we have $C(w)\supset T$ (where as above I have set $w=\sigma(\xi)$). Now suppose that $C(\xi)\cap S_1=\{1\}$. If $t$ is an element of $C( w)\cap T$, then we may write $t=\sigma(s)$ for some $s\in S_1$. Since $S_1\triangleleft\Delta$, we have $t\xi t^{-1}\xi^{-1}\in S_1$. But $\sigma(t\xi t^{-1}\xi^{-1})=u w u^{-1} w^{-1}=1$, since $u\in C( w)$. Since $\sigma|S_1$ is injective, it follows that $t\xi t^{-1}\xi^{-1}=1$, i.e.  $t$ commutes with $\xi$. Since $C(\xi)\cap S_1=\{1\}$, this implies that $t=1$; it follows that $C( w)\cap T=1$. Thus I have shown that either $C( w)\supset T$ or $C( w)\cap T=1$, as required for Lemma \ref{wide stance}.

Another hypothesis of Lemma \ref{wide stance} is that the image $\overline w$ of $w$ under the quotient map $Q\to Q/T$ has order $m$. To verify this, note that because 
$\sigma:\Delta\to Q$ has kernel $S_0$, and
the quotient homomorphism $\mu:Q\to Q/T$ has kernel $T=\sigma(S_1)$, the kernel of $\mu\sigma$ is $S_0S_1=S$. Hence the order of $\overline w=\mu\sigma(\xi)$ is equal to the order of the image $\overline\xi$ of $\xi$ under the the quotient homomorphism $\Delta\to\Delta/S$, which by hypothesis is $m$.

The last  hypothesis of Lemma \ref{wide stance} is that the normal closure of $\overline w$ in $Q/T$ is abelian. Since I have observed that $\overline w=\mu\sigma(\xi)$, and that the kernel of $\mu\sigma$ is $S$, it suffices to show that the image $\overline\xi$ of $\xi$ under the the quotient homomorphism $\Delta\to\Delta/S$ has abelian normal closure in $\Delta/S$. But according to the hypothesis of the present lemma, $\langle\overline\xi\rangle$ is normal in $\Delta/S$; in particular the normal closure of $\overline\xi$ in $Q/T$ is abelian. 
\EndProof

\section{Two-generator subgroups of linear groups}\label{two-generator section}

This section contains the proof of Corollary \ref{the daughter of rosie o'grady}, which was discussed in the Introduction.

\Lemma\label{basket driving}
Let $k$ be a finite field whose characteristic  $p$ is odd, and let $ G$ be a subgroup of
$\pizzle(k)$. Suppose that $G$ has an element of order strictly greater than $5$.
Then at least one of the following alternatives
holds:
\Alternatives
\item\label{tant pis} $ G$ is isomorphic to $\pizzle(K)$ or to to $\piggle(K)$, where $K$ is a subfield of $k$;
\item\label{puisque je te dis} $ G$ has a (possibly trivial) normal subgroup $T$ which is an elementary abelian $p$-group, and $G/T$ is a (possibly trivial) cyclic group whose order is relatively prime to $p$; or
\item\label{encore un carreau} $G$ is dihedral, and its order is prime to $p$.
\EndAlternatives
\EndLemma

\Proof
This is included in Dickson's classification of subgroups of $\pizzle(k)$ \cite[Hauptsatz 8.27]{huppert}. In \cite{huppert}, eight types of subgroups are listed, of which type (8) satisfies Alternative (i) above; types (1), (2) and (7) satisfy Alternative (ii); type (3) satisfies Alternative (iii); and types (4), (5) and (6) have no elements of order strictly greater than $5$.
\EndProof

\Proposition\label{do monkeys have uncles?}
Let $v$ be a valuation of a number field $E$, and suppose that the characteristic $p$  of $k_v$ is odd. Let $m$ be an odd integer with $7\le m<p$.
Let $x$ and $y$ be elements of $\pizzle(\frako_v)$ such that the group $\Gamma:=\langle x,y\rangle$ is and non-solvable. 
Suppose that 
$\overline h_v(x)$
 has  order $m$ in $\pizzle(k_v)$.
Then there exist a finite group $Q$ and a surjective homomorphism $\lambda:\Gamma\to Q$ such that 
\begin{enumerate}
\item every cyclic subgroup of $Q$ has index at least $p$, and
\item either $\lambda( x)$ has order $m$ in $Q$, or the normal closure of $\lambda( x)$ in $Q$ is abelian.
\end{enumerate}
\EndProposition

\Proof
Set $k=k_v$, $h=\overline h_v|\Gamma$, and $G=h(\Gamma)$. Set $\overline x=h(x)$. The group $G$ is a subgroup of $\pizzle(k)$, and $\overline x\in G$ has order $m\ge7$. Hence $G$ must satisfy one of the alternative conclusions (i)---(iii) of Lemma \ref{basket driving}.

Set $K=\ker h\triangleleft\Gamma$. It follows from Propositions \ref{rheingold is as good to your taste} and \ref{fantabula} that $K$ is residually $p$-nilpotent. 

First suppose that $G$ satisfies (i), so that $ G\cong\pizzle(K)$ or to $ G\cong\piggle(K)$ for some subfield $K$ of $k$. In this case I claim that the conclusions of the proposition are true if we take $Q=G$ and $\lambda=h$. Since $\overline x$ has order $m$, Conclusion (2) holds. To verify Conclusion (1), let $q$ denote the order of $K$, a power of $p$. Any cyclic subgroup $J$  of $\pizzle(K)$ or $\piggle(K)$ has order at most $q+1$. We have $|\piggle(K)|\ge|\pizzle(K)|=(q^3-q)/2$ (since $p$ is odd), and hence
$$|G:J|\ge \frac{(q^3-q)/2}{q+1}=\frac{q^2-q}2\ge q\ge p.$$

Next suppose that $G$ satisfies (ii), so that $ G$ has a normal subgroup $T$ which is an elementary abelian $p$-group, and $G/T$ is cyclic. 
Set $A=G/T$, let $\pi:G\to A$ denote the quotient homomorphism, and set $\ell =\pi\circ h:\Gamma\to A$. If we set $L=\ker \ell =h^{-1}(T)$, then $K\triangleleft L\triangleleft\Gamma$. 
Since $K$ is residually $p$-nilpotent and $L/K\cong T$ is an elementary abelian $p$-group, it follows from Proposition \ref{nabla} that $L$ is residually $p$-nilpotent. The hypothesis implies that the order $m$ of $\overline x$ is prime to $p$ and therefore to the order of $T$; hence $\ell (x)=\pi(\overline x)$ also has order $m$. According to Proposition \ref{it's so diversible it's reversible}, $(x,y)$ is a reversive system for $\Gamma$. Thus $\Gamma$, $x$, $y$, $p$, $A$, $\ell $ and $m$ satisfy all the hypotheses of  Lemma \ref{i shoulda stood in bed}. Hence there exist a finite group $\Delta $ and a homomorphism $\lambda:\Gamma\to \Delta $ such that conditions (1)---(4) of Lemma \ref{i shoulda stood in bed} hold.

Set $\xi=\ell (x)$. Since Conditions (1)--(3) of Lemma \ref{i shoulda stood in bed} hold, it follows from Proposition \ref{unimodular} that
the automorphism $s\mapsto \xi s\xi^{-1}$ of $S$ is semi-unimodular (where $S$ is defined by Condition (2)). On the other hand, since $S$ is abelian, the action of $\Delta$ on $S$ by conjugation factors through an action of $\Delta/S$. By Condition (3) of Lemma \ref{i shoulda stood in bed}, the image of $\overline\xi$ of $\xi$ in $\Delta/S$ has order $m$,  which by hypothesis is odd. It therefore follows from \ref{paris hat} that the automorphism $s\mapsto \xi s\xi^{-1}$ of $S$ is unimodular. Note also that since $\Delta/S$ is abelian by Condition (4) of  Lemma \ref{i shoulda stood in bed}, the subgroup $\langle\overline\xi\rangle$ of $\Delta/S$ is normal. Furthermore, by Condition (2) of Lemma \ref{i shoulda stood in bed}, the elementary abelian $p$-group $S$ has rank at least $2$. Thus the hypotheses of Lemma \ref{happenstance} hold.  Hence there exist a group $Q$ and a surjective homomorphism $\sigma:\Delta\to Q$, such that
(1) every cyclic subgroup of $Q$ has index at least $p$, and
(2) either $\sigma(\xi)$ has order $m$ in $Q$, or the normal closure of $\sigma(\xi)$ in $Q$ is abelian. If we set $\lambda=\sigma\circ \ell :\Gamma\to Q$, the conclusions of the theorem follows in this case.

There remains the case in which $G$ satisfies (iii), so that $ G$ is dihedral and its order is prime to $p$.
Let $K_1$ denote the mod $p$ commutator subgroup (\ref{adams spectral nahant}) of $K$. Since $K_1$ is a characteristic subgroup of $K\triangleleft\Gamma$, we have $K_1\triangleleft\Gamma$. In this case I will define $\Delta$ to be the quotient group $\Gamma/K_1$, and I will let $S$ denote the subgroup $K/K_1$ of $\Delta$. Since $K\triangleleft\Gamma$ we have $S\triangleleft\Delta$. The quotient group $\Delta/S=(\Gamma/K_1)/(K/K_1)$ is isomorphic to $\Gamma/K\cong G$ and is therefore dihedral. 

Let $\xi$ denote the image of $x$ under the quotient homomorphism $\Gamma\to\Delta$. Let $\overline\xi$ denote the image of $\xi$ under the quotient homomorphism $\Delta\to \Delta/S$. The composition of the quotient homomorphisms $\Gamma\to\Delta$ and $\Delta\to \Delta/S$ has kernel $K=\ker h$; hence the order of $\xi$ is equal to the order of $\overline x=h(x)$, namely $m$.

Since $\Gamma/K\cong G$ is solvable, and since $\Gamma$ is non-solvable by hypothesis, $K$ must be non-solvable, and in particular non-cyclic. Hence by Proposition \ref{squash}, $S=K/K_1$ is an elementary abelian $p$-group of rank at least $2$. Since $|\Delta/S|=|G|$ is prime to $p$, the subgroup $S$ is a $p$-Sylow subgroup of $\Delta$.

Since $S$ is abelian, the action of $\Delta$ on $S$ by conjugation factors through an action of $\Delta/S$. Since $\Delta/S$ is dihedral and the order $m$ of $\overline\xi\in\Delta/S$ is odd, $\overline\xi$ belongs to the commutator subgroup of $\Delta/S$. Hence by \ref{harris pat}, the automorphism $s\mapsto \xi s\xi^{-1}$ of $S$ is unimodular. Furthermore, since $\overline\xi$ is an element of odd order in the dihedral group $\Delta/S$, the subgroup $\langle\xi\rangle$ of  $\Delta/S$ is normal. It now follows from Lemma \ref{happenstance} that there exist a group $Q$ and a surjective homomorphism $\sigma:\Delta\to Q$, such that
(1) every cyclic subgroup of $Q$ has index at least $p$, and
(2) either $\sigma(\xi)$ has order $m$ in $Q$, or the normal closure of $\sigma(\xi)$ in $Q$ is abelian. If we define $\lambda$ to be the composition of the quotient homomorphism $\Gamma\to\Delta$ with $\sigma$, the conclusions of the theorem follow in this case.
\EndProof

\Corollary\label{the daughter of rosie o'grady} Let $v$ be a valuation of a number field $E$, and suppose that the characteristic $p$ of $k_v$ is odd. Let $m$ be an odd integer with $7\le m<p$. Let $x$ and $y$ be elements of $\pizzle(\frako_v)$ such that the group $\Gamma:=\langle x,y\rangle$ is non-solvable.  Suppose that $\overline h_v(x)$
 has  order $m$ in $\pizzle(k_v)$. Set 
$$\Theta_1=\langle x^m,y\rangle\le\Gamma$$
and
$$\Theta_2=\langle x, yxy^{-1}xyx^{-1}y^{-1}\rangle\le\Gamma,$$
and set $\theta_i=|\Gamma:\Theta_i|$ for $i=1,2$. Then
$$\max(\theta_1,\theta_2)\ge p.$$
\EndCorollary

\Proof
Note that $E$, $v$, $p$, $x$ and $\Gamma$ satisfy the hypotheses of Proposition \ref{do monkeys have uncles?}.
Let $Q$ be a finite group, and $\lambda:\Gamma\to Q$ a surjective homomorphism, such that Conclusions (1) and (2) of Proposition \ref{do monkeys have uncles?} hold.
Set $\overline x=\lambda(x)$ and $\overline y=\lambda(y)$. For $i=1,2$, set $\overline\Theta_i=h(\Theta_i)$, and note that $\theta_i=|G:\overline\Theta_i|$. 

According to Conclusion (2) of Proposition \ref{do monkeys have uncles?}, either $\overline x$ has order $m$ in $Q$, or the normal closure of $\overline x$ in $Q$ is abelian.
If $\overline x$ has order $m$, we have $\overline\Theta_1=\langle \overline x^m,\overline y\rangle=\langle\overline y\rangle$. 
In particular $\overline\Theta_1$ is cyclic, and by
Conclusion (1) of Proposition \ref{do monkeys have uncles?} we have $\theta_1=|G:\overline\Theta_1|\ge p$.
If the normal closure of $\overline x$ in $Q$ is abelian, then 
$\bar y\bar x\bar y^{-1}$ commutes with $\overline x$, and hence
$$\overline \Theta_2=\langle \bar x, \bar  y \bar  x \bar  y^{-1}\bar  x \bar  y \bar  x^{-1}\bar  y^{-1} \rangle=\langle \bar x, \bar x\rangle=\langle \bar x\rangle.$$
In particular $\overline\Theta_2$ is cyclic, and by
Conclusion (1) of Proposition \ref{do monkeys have uncles?} we have $\theta_2=|G:\overline\Theta_2|\ge p$.
\EndProof

\section{The tree for $\pizzle$}\label{pizzletree}

\Number\label{personal savior}
A {\it simplicial complex} will be understood to be geometric, except when it is specified to be abstract. The underlying space of a simplicial complex $K$, i.e. the union of its simplices, will be denoted $|K|$. 

If $v$ is a vertex in a simplicial complex $K$, I will denote the (open) star of $v$ in $K$ by ${\rm St}_K(v)$.

A simplicial map $f:K\to L$, where $K$ and $L$ are simplicial complexes, will be termed {\it nondegenerate} if $\dim f(\sigma)=\dim\sigma$ for every simplex $\sigma$ of $K$.
\EndNumber

\Number\label{tree stuff-not}
By a {\it tree} I will mean a $1$-connected simplicial complex of dimension $0$ or $1$. If $s$ and $t$ are vertices of a tree $T$, there is a unique arc in $T$ with endpoints $s$ and $t$, which I will denote $[s,t]$. The {\it distance} between $s$ and $t$ is defined to be the length of $[s,t]$, i.e. the number of edges that it contains.
With this definition of distance, $T$ becomes a metric space. 

If $\Gamma$ is a group, I will define a $\Gamma$-tree to be a tree equipped with a simplicial action of $\Gamma$, which has no inversions in the sense that if an edge $e$ is invariant under an element $\Gamma\in \gamma$, then both vertices incident to $e$ are fixed by $\gamma$. A $\Gamma$-tree $T$ will be called {\it essential} if no vertex of $T$ is fixed by the entire group $\Gamma$. (The term ``nontrivial'' has been used for this notion elsewhere, but I am avoiding it in the present paper to prevent confusion.) A $\Gamma$-tree $T$ will be called {\it faithful} if no nontrivial element of $\Gamma$ acts by the identity on $T$. It will be called {\it minimal} if no proper subtree of $T$ is $\Gamma$-invariant. (Note that a minimal $\Gamma$-tree $T$ is essential unless $T$ consists of a single vertex.)
\EndNumber

\Number\label{If it ain't broke}
If $\Gamma$ is a group and $T$ is a $\Gamma$-tree, then for each subgroup $H$ of $\Gamma$ I will denote by $\Fix_T(H)$ the set of all points of $T$ that are fixed by the entire subgroup $H$. Since $\Gamma$ acts simplicially and without inversions  on $T$, the set $\Fix_T(H)$ is a subcomplex of $T$. If $s$ and $t$ are vertices of $\Fix_T(H)$, then any element of $H$ maps $[s,t]$ simplicially onto an arc with endpoints $s$ and $t$, which must be $s$ and $t$ itself; hence $[s,t]\subset\Fix_T(H)$. This shows that if $\Fix_T(H)$ is non-empty then it is connected, and is therefore a subtree of $T$.
\EndNumber

\Number\label{elliptical billiard balls}
Let $T$ be a $\Gamma$-tree. A non-trivial element $\gamma$ of $\Gamma$ is said to be {\it $T$-elliptic} if $\Fix_T(\gamma):=\Fix_T(\langle\gamma\rangle)\ne\emptyset$, and {\it $T$-hyperbolic} if $\Fix_T(\gamma)=\emptyset$. Thus $\gamma$ is $T$-elliptic if and only if
it fixes some vertex of $T$. By the proof of \cite[1.3]{culler-morgan}, $\gamma$ is $T$-hyperbolic if and only if it has an {\it axis}: this is a subcomplex of $T$, isomorphic to the real line triangulated with vertices precisely at the integer points, on which $\gamma$ acts via an integer translation. The proof of  \cite[1.3]{culler-morgan} also shows that the axis of $\gamma$ is contained in every $\gamma_0$-invariant subtree of $\calt$; in particular it is unique.
\EndNumber

\Number\label{periodic billiards}
If $\gamma$ is an element of a group $\Gamma$ and $T$ is a $\Gamma$-tree, I will set $\Per_T(\gamma)=\bigcup_{n=1}^\infty\Fix_T(\gamma^n)$. Rewriting 
$\Per_T(\gamma)$ as $\bigcup_{n=1}^\infty\Fix_T(\gamma^{n!})$ show that
$\Per_T(\gamma)$ is a monotone union of subsets each of which is either empty or a subtree. Hence if
$\Per_T(\gamma)$ is non-empty then it is a subtree of $T$.

Note that if $\gamma$ is $T$-hyperbolic, then its axis is also an axis for any positive power of $\gamma$; hence all powers of $\gamma$ are $T$-hyperbolic, and so $\Per_T(\gamma)=\emptyset$. This shows that
$\Per_T(\gamma)\ne\emptyset$ if and only if $\gamma$ is $T$-elliptic.
If $e$ is an edge in $\Per_T(\gamma)$, I will define the $x$-period of $e$ to be the smallest $d>0$ such that $e\subset\Fix_T(\gamma^n)$.
\EndNumber

\Lemma\label{miminumule}
If $\Gamma$ is a finitely generated group and $\calt$ is a $\Gamma$-tree, some $\Gamma$-invariant subtree of $\calt$ is a minimal $\Gamma$-tree.
\EndLemma

\Proof
If every element of $\Gamma$ has a fixed vertex in $T$, then according to \cite [p. 64, Corollary 2]{Serre}, the $\Gamma$-tree $T$ is inessential. Hence there is a $\Gamma$-invariant subtree of $\calt$ consisting of a single vertex, which is obviously minimal. Now suppose that some element $\gamma_0$ of $\Gamma$ is $T$-hyperbolic. According to \ref{elliptical billiard balls}, $\gamma_0$ has an axis $A_0$, which is contained in every $\gamma_0$-invariant subtree of $\calt$. In particular, every $\Gamma$-invariant subtree of $\calt$ contains $A_0$. Hence if we denote by $\calF$ the set of all $\Gamma$-invariant subtrees of $\calt$ and set $T_0=\bigcap_{T\in\calF}T$, we have $A_0\subset T_0$ and hence $A_0\ne\emptyset$. If $s$ and $s'$ are vertices of $T_0$, each $T\in\calF$ contains $s$ and $s'$, and therefore contains the unique arc $[s,s']$ with endpoints $s$ and $s'$; hence $[s,s']\subset T_0$. This shows that $T_0$ is a subtree of $\calt$; it is immediate that $T_0$ is $\Gamma$-invariant, and that as a $\Gamma$-tree it is minimal.
\EndProof

\Proposition\label{mr fix it}Let $\Gamma$ be a group, let $T$ be a $\Gamma$-tree, and let $u$ and $v$ be $T$-elliptic elements of $\Gamma$ such that $uv$ is also $T$-elliptic. Then $u$ and $v$ have a common fixed vertex in $T$.
\EndProposition

\Proof
I will regard $T$ as a $\ZZ$-tree in the sense of \cite[Chapter II] {ms1}. If $s$ is any vertex of $T$, and $\gamma\in\Gamma$ is $T$-elliptic, then according to
\cite[Lemma II.2.16]{ms1}, 
the length of $[s,\gamma\cdot s]$ is even, and the midpoint of $[s,\gamma\cdot s]$ is fixed by $\gamma$. If we take $s$ to be a fixed vertex of $v$, it follows that the midpoint of $[s,u\cdot s]=[s,uv\cdot s]$ is fixed by both $u$ and $uv$ (and hence by $v$).
\EndProof

\Number\label{wizzlebee}
In the remainder of this section I will be using the tree for $\pizzle(E)$, where $E$ is a field with a valuation. I will be taking the point of view presented in \cite[Section 3]{handbook}, which I will briefly review here; it is a mild variant of the point of view used in \cite{Serre}.

 Let $v$  be a valuation $v$ of a field $E$, and set $\frako=\frako_v$. Let $V=K^2$ denote the
standard $2$-dimensional vector space over $K$, which is in particular an $\frako$-module. We define a {\it lattice} in $V$ to be an
$\frako$-submodule of $V$ which is finitely generated and spans
$V$ as a vector space over $K$. It is pointed out in \cite[Subsection 3.6]{handbook} that any lattice
is a free $\frako$-module of rank $2$. Two lattices
$\Lambda,\Lambda'\subset V$ are said to be
{\it (homothety)-equivalent,} or to represent the same homothety
class, if there is a nonzero element $\alpha$ of $K$ such
that $\Lambda'=\alpha\Lambda$. 

It follows from \cite[Lemma 3.6.8]{handbook} that any two homothety classes of lattices $s_0$ and $s_1$ have respective representatives $\Lambda_0$ and $\Lambda_1$ such that (1) $\Lambda_1\subset\Lambda_0$ and (2) $\Lambda_0/\Lambda_1$ is isomorphic as an $\frako$-module to $\frako/\beta\frako$ for some nonzero element $\beta$ of $\frako$. (There is a typographical error in the statement of \cite[Lemma 3.6.8]{handbook}, where the last two occurrences of $\Lambda_1$ should be replaced by $\Lambda_1'$.) It is shown in the discussion following the proof of \cite[Lemma 3.6.8]{handbook} that the non-negative integer $v(\beta)$ depends only on $s_0$ and $s_1$ and not on the choice of representatives satisfying (1) and (2). I will denote it $\Delta(s_0,s_1)$. According to \cite[Lemma 3.6.11]{handbook}, the set $\calv$ of homothety classes of $\frako$-lattices is a metric space with distance function $\Delta$. Furthermore, if we define an abstract simplicial $1$-complex $\frakT_v$ by defining the vertex set of $\frakT_v$ to be $\calv$ and defining its $1$-simplices to pairs of vertices $\{s,t\}$ such that $\Delta(s,t)=1$, then according to \cite[Theorem 3.6.14]{handbook}, a geometric realization $\calt_v$ of $\frakT_v$ is a tree. I will identify $\calv$ with the vertex set of $\calt_v$. The proof of \cite[Theorem 3.6.14]{handbook} shows that if $s$ and $t$ are vertices of $\calt_v$ then $\Delta(s,t)$ is the distance between $s$ and $t$ in the sense defined in \ref{tree stuff-not}. \EndNumber

\Number\label{kurimu}
According to the discussion in \cite{handbook} following the proof of Lemma 3.6.8, if $s$ and $s'$ are vertices of $\calt_v$ such that $\Delta(s,s')=n$, we may represent $s$ and $s'$ by lattices $\Lambda$ and $\Lambda'$ which respectively have bases of the form $(e,f)$ and $(e,\beta f)$ for some $\beta\in\frako$ with $v(\beta)=n$. Since the $\frako$-module generated by $e$ and $\beta f$ is unchanged when $\beta$ is  multiplied by a unit in $\frako$, we may take $\beta=\pi^n$ where $\pi$ is a prescribed generator for $\frakm_v$.
\EndNumber

\Number\label{act while the acting is good}
%If $\Lambda$ and $\Lambda'$ are homothety-equivalent lattices, then $A(\Lambda)$ and $A(\Lambda')$ are obviously homothety-equivalent for any $A\in\ggle(E)$. Hence we may define an action of $\ggle(E)$ on $\calt_v$ by defining $A\cdot s$ to be the homothety class of 
It is shown in \cite[Subsection 3.7]{handbook} that there is a unique simplicial action of $\zzle(E)$ on $\calt_v$ such that if $A$ is an element of $\ggle(E)$ and $\Lambda$ is a lattice representing a vertex $s$ of $\calt_v$, then $A\cdot s$ is represented by the lattice $A(\Lambda)$. Since $-I(\Lambda)=\Lambda$ for every lattice $\Lambda$, the action of $\ggle(E)$ on $\calt_v$ factors through an action of $\piggle(E)$. By restriction we obtain natural actions of $\zzle(E)$ and $\pizzle(E)$ on $\calt_v$. I will always regard $\calt_v$ as a $\pizzle(E)$-tree via this action.

Since every $\frako$-lattice is generated, as an $\frako$-module, by some basis of $V$, and since $\ggle$ acts transitively on the bases of $V$, the action of $\ggle(E)$ (or of $\piggle(E)$) on $\calt_v$ is also transitive.
\EndNumber

One especially important property of the action is:

\Proposition\label{don't move-a da lattice}
Let $v$ be a valuation of a field $E$, let $s$ be a vertex of $\calt_v$, and let $\Lambda$ be a lattice representing $s$. Then the stabilizer $\pizzle(E)_s$ of $s$ in $\pizzle(E)$ consists of all elements of $\pizzle(E)$ that leave $\Lambda$ invariant. If $\Lambda$ is the lattice $\calo_v^2\subset E^2$, we have $\pizzle(E)_s=\pizzle(\calo_v)$. For any vertex $s$ of $\calt_v$ the stabilizer $\pizzle(E)_s$ is a conjugate of $\pizzle(\calo_v)$ in $\piggle(E)$.
\EndProposition

\Proof
The first assertion is proved in Subsection 3.7 of \cite{handbook}. The second assertion follows from the first because $\pizzle(\calo_v)$ is the stabilizer of the lattice $\calo_v^2$ in $\pizzle(E)$. The third assertion follows from the section because $\piggle(E)$ acts transitively on $\calt$ by \ref{act while the acting is good}.
\EndProof

\Number\label{no tomorrow}
Recall that if $k$ is a field, the {\it projective line} $k{\rm P}^1$ is defined to be the set of all $1$-dimensional subspaces of the standard $2$-dimensional $k$-vector space $k^2$. Under the standard action of $\ggle(k)$ on $k^2$, each element of $\ggle$ carries each $1$-dimensional subspace of $k^2$ onto a $1$-dimensional subspace; hence the action defines an action of of $\ggle(k)$ on $k{\rm P}^1$. Since $-I$ acts obviously lies in the kernel of this action, the action factors through an action of $\piggle(k)$. By restriction we obtain a {\it standard} action of $\pizzle(k)$ on $k{\rm P}^1$.
\EndNumber

\Proposition\label{don't squeeze-a da fruit}
Let $v$ be a valuation of a field $E$, set $\frako=\frako_v$, and let $s_0$ denote the vertex of $\calt_v$ represented by the lattice $\frako^2\subset E^2$. Set $k=k_v$. Then there is a $\pizzle(\frako )$-equivariant bijection between the link $L$ of $s_0$ in $\calt$, equipped with the action of $\pizzle(E)_{s_0}=\pizzle(\frako)$ obtained by restricting the action of $\pizzle(E)_{s_0}$ to $L$, and the projective line $k{\rm P}^1$ equipped with the action of $\pizzle(\frako)$ obtained by pulling back the standard action (\ref{no tomorrow}) of $\pizzle(k)$ on $k{\rm P}^1$ via $\overline h_v:\pizzle(\frako )\to\pizzle(k)$.
\EndProposition

\Proof
Set $V=k^2$ and $\Lambda _0=\frako^2$. According to \cite[Lemma 3.6.8]{handbook} (with the correction pointed out in \ref{wizzlebee} above), each vertex in $L$ is represented by a unique lattice $\Lambda $ such that (1) $\Lambda \subset \Lambda _0$ and (2) $\Lambda_0 /\Lambda\cong\frako/(\pi)$. Conversely, it follows from the definition of $\calt_v$ that every lattice $\Lambda$ satisfying (1) and (2) represents a vertex in $L$. Hence if $\calL$ denotes the set of all lattices $\Lambda$ satisfying (1) and (2), there is a natural bijection $\alpha:\calL\to L$ that maps each lattice to its homothety class.

Note that if a lattice $\Lambda$ satisfies Conditions (1) and (2), then  $\pi\Lambda_0 \subset\Lambda \subset \Lambda _0$. Hence if $\calL'$ denotes the set of all lattices $\Lambda$ such that $\pi\Lambda_0 \subset\Lambda \subset \Lambda _0$, we have $\calL\subset\calL'$.
On the other hand, if we identify $\Lambda_0/\pi\Lambda_0=\frako^2/\pi(\frako^2)$ with $V=(\frako/\pi\frako)^2$, and let $p:\Lambda_0\to V$ denote the quotient map, then there is a bijection $\Lambda\to p(\Lambda)$ from $\calL'$ to the set of subspaces of the vector space $V$. For any $\Lambda\in\calL'$ we have $\Lambda_0/\Lambda\cong V/p(\Lambda)$. In particular, we have $\Lambda_0/\Lambda\cong\frako/\pi\frako$ if and only if $V/p(\Lambda)$ is $1$-dimensional, i.e. if and only if $p(\Lambda)$ is $1$-dimensional. Hence the bijection $\Lambda\to p(\Lambda)$ restricts to a bijection $\beta$ of $\calL$ onto $k{\rm P}^1$. Thus $J:=\beta\circ\alpha^{-1}:L\to k{\rm P}^1$ is a bijection. It remains to show that $J$ is $\pizzle(\frako)$-equivariant in the sense of the statement of the proposition.

Any element of $\pizzle(\frako_v )$ may be written in the form $[A]$ for some $A\in\zzle[ \frako_v]$. If we identify $A$ and $h_v(A)$ with linear automorphisms of the $\frako$-module $L_0=\frako^2$ and the vector space $V$ respectively, we have $p\circ A=h_v(A)\circ p$. Now let $s$ be any vertex in $L$, and set $\Lambda=\alpha^{-1}(s)$. The $\frako$-lattice $A(\Lambda_s)$ is contained in $A(\Lambda_0)=\Lambda_0$, and satisfies $\Lambda_0/A(\Lambda_s)=A(\Lambda_0)/A(\Lambda_s)\cong \Lambda_0/\Lambda_s\cong\frako/(\pi)$. Hence $A(\Lambda)\in\calL$. The definition of the action of $\zzle(E)$ on $\calt_v$ implies that $A(\Lambda)\in\calL$
represents the vertex $[A]\cdot s$, so that $A(\Lambda)=\alpha^{-1}([A]\cdot s)$.
Hence $$J([A]\cdot s)
=p(A(\Lambda_s))=h_v(A)(p(\Lambda_s))=h_v(A)(J(s))=[h_v(A)]\cdot J(s)=\overline h_v([A])\cdot J(s),$$
which proves equivariance. 
\EndProof

\Corollary\label{i'm a stranger here myself}
If $v$ is a valuation of a number field $E$, the tree $\calt_v$ is locally finite.
\EndCorollary

\Proof
If $s_0$ is any vertex of $\calt_v$, we must show that the link of $s_0$ is finite. Since $\piggle(E)$ acts transitively on $\calt$ by \ref{act while the acting is good}, we may assume without loss of generality that $s_0$ is represented by the lattice $\frako^2=\frako_v^2\subset E^2$. Then by Proposition \ref{don't squeeze-a da fruit}, the link of $s_0$ is in bijective correspondence with $k{\rm P}^1$, where $k=k_v$. Since $E$ is a number field, $k_v$ is a finite field, and hence  $k{\rm P}^1$ is finite.
\EndProof

\Corollary\label{who isn't these days} Let $v$ be a valuation of a number field $E$, and let $p$ denote the characteristic of $k_v$. Let $ A$ be an element of $\zzle(k)$ which satisfies $\trace A\in \frako_v$ and $\trace A\equiv\pm2\pmod{\frakm_v}$. Suppose that $s$ is a vertex of $\Fix_{\calt_v}( [A])$. Then $[ A]^p$ fixes each vertex in the link of $s$ in $\calt_v$.  
\EndCorollary

\Proof
Since $\piggle(E)$ acts transitively on $\calt$ by \ref{act while the acting is good}, we may assume without loss of generality that $s$ is represented by the lattice $\frako^2=\frako_v^2\subset E^2$. According to Proposition \ref{don't move-a da lattice} we have $\pizzle(E)_s=\pizzle(\calo_v)$, and by Proposition \ref{don't squeeze-a da fruit} (applied with $s_0=s$), the action of $\pizzle(\calo_v)$ on the link $L$ of $s$ is the pullback of an action of $\pizzle(k_v)$ via the homomorphism $\overline h_v:\pizzle(\calo_v)\to\pizzle(k_v)$. Now since $\trace A\equiv\pm2\pmod{\frakm_v}$, we have $\trace h_v(A)=\pm2$, so that $h_v(A)$ is conjugate in $\ggle(k_v)$ to  $\pm\begin{pmatrix}1&\lambda\cr0&1\end{pmatrix}$ for some $\lambda\in k_v$. Hence $\overline h_v([A]^p)=1$, and so $[A]$ acts trivially on $L$. 
\EndProof

\Proposition\label{it lives}
Let $v$ be a valuation of a number field $E$, and let $p$ denote the characteristic of $k_v$. Let  $\Gamma$ be finitely generated subgroup of $\pizzle(E)$. Suppose that $\Gamma$ has no non-trivial normal abelian subgroup. Suppose also that $\Gamma$ is not conjugate in $\piggle(E)$ to a subgroup of $\pizzle(\frako_v)$. Then:
\begin{enumerate}
\item There exists an essential, faithful, minimal $\Gamma$-tree $T$. 
\item If $x$ is an element of $\Gamma$ which is not conjugate in $\piggle(E)$ to an element of $\pizzle(\frako_v)$, then $T$ may be chosen in such a way that $x$ is $T$-hyperbolic. 
\item If $x\in\Gamma\cap\pizzle(\frako_v)$, and if $m$ denotes the order of $\overline h_v(x)\in\pizzle(k_v),$ then $T$ may be chosen in such a way that $x$ is $T$-elliptic, $\Per_T(x)=T$, and the $x$-period of every edge of $T$ divides $mp^r$ for some $r\ge0$.
\end{enumerate}
\EndProposition

\Proof
Fix a generator $\pi$ for $\frakm_v$.

Set $\calt=\calt_v$ (see \ref{wizzlebee}). Since $\calt$ is a $\pizzle(E)$-tree, it is in particular a $\Gamma$-tree. Since by hypothesis $\Gamma$ is not conjugate in $\piggle(E)$ to a subgroup of $\pizzle(\frako_v)$, it follows from Proposition \ref{don't move-a da lattice} that no vertex of $\calt$ is fixed by $\Gamma$, i.e. that $\calt$ is an essential $\Gamma$-tree. By Corollary \ref{i'm a stranger here myself}, the tree $\calt$ is locally finite.

It follows from Lemma \ref{miminumule} that $\calt$ has a $\Gamma$-invariant subtree $T$ which, regarded as a $\Gamma$-tree, is minimal. Since $\calt$ is an essential, locally finite $\Gamma$-tree, the same is true of $T$.

If $T$ is finite, the action of $\Gamma$ on $T$ factors through an action of a finite quotient of $\Gamma$; it then follows from \cite[Example 6.3.1]{Serre} that $\Gamma$ fixes a point of $T$, a contradiction. Hence $T$ is infinite. Since $T$ is locally finite, we have $\diam T=\infty$. Hence if we fix any base vertex $v_0$, there exists for each $n>0$ a vertex $v_n$ such that $\Delta(v_0,v_n)=n$. Set $Q_n=\Gamma_{v_0}\cap\Gamma_{v_n}$ for each $n>0$. Set $\tQ_n=\Pi_E^{-1}(Q_n)\le\zzle(E)$.

The proof that $T$ is faithful will depend on the following fact:

\Claim\label{mod on}
For every $n>0$, and for every $A\in[\tQ_n,\tQ_n]$, we have $\trace A\equiv2\pmod{\pi^n}$.
\EndClaim

To prove \ref{mod on}, we first note that by \ref{kurimu}, since $\Delta(v_0,v_n)=n$, we may represent $v_0$ and $v_n$ by lattices $\Lambda_0$ and $\Lambda_n$ which respectively have bases of the form $(e,f)$ and $(e,\pi^nf)$. After a conjugation of $\Gamma$ in $\piggle(E)$ we may assume that $e=(1,0)$ and $f=(0,1)$, so that $\pi^n f=(0,\pi^n)$. Hence the element $B:=\begin{pmatrix}1&0\cr0&\pi^n\end{pmatrix}$ of $\zzle(E)$ carries $\Lambda_0$ onto $\Lambda_1$. 

According to Proposition \ref{don't move-a da lattice}, $\Gamma_{v_0}$ and $\Gamma_{v_n}$ respectively leave the lattices $\Lambda_0$ and $\Lambda_n$ invariant. It follows that $\Pi_E^{-1}(\Gamma_{v_0})\le\zzle(\frako_v)$, and that
$$\Pi_E^{-1}(\Gamma_{v_n})\le  B\cdot\zzle(\frako_v)\cdot B^{-1}=\bigg\{\begin{pmatrix}a&\pi^{-n}b\cr\pi^nc&d\end{pmatrix}:a,b,c,d\in\frako_v\bigg\}.$$ 
Hence $\tQ_n=\Pi_E^{-1}(\Gamma_{v_0})\cap \Pi_E^{-1}(\Gamma_{v_n})$ consists of matrices which are upper triangular modulo $\pi^n$. This implies that $[\tQ_n,\tQ_n]$ consists of matrices  which are upper triangular modulo $\pi^n$ and whose diagonal entries are congruent to $1$ modulo $\pi^n$. This gives \ref{mod on}.

To show that $T$ is faithful, let $N\triangleleft\Gamma$ denote the kernel of the action of $\Gamma$ on $T$, consisting of all elements of $\Gamma$ that act on $T$ by the identity. We have $N\le Q_n$ for every $n>0$. Set $\tN=\Pi_E^{-1}( N)$; then every $n>0$ we have $\tN\le \tQ_n$ and hence $[\tN,\tN]\le[\tQ_n, \tQ_n]$. Hence for every $A\in[\tN,\tN]$ and for every $n>0$, it follows from \ref{mod on} that $\trace A\equiv2\pmod{\pi^n}$. Since $\bigcap_{n=1}^\infty(\pi^n)=0$, we have
$\trace A=2$ 
for each $A\in[\tN,\tN]$. It then follows from \cite[Lemma 1.2.1]{varreps} that $\tN$ is a reducible subgroup of $\zzle(E)$. In particular, $[\tN,\tN]$ is abelian. Hence $[N,N]$ is abelian. But since $N$ is normal in $\Gamma$, its characteristic subgroup $[N,N]$ is also normal. By hypothesis, $\Gamma$ has no non-trivial normal abelian subgroup; hence $[N,N]=\{1\}$. This implies that the normal subgroup $N$ of $\Gamma$ is abelian, and a second application of the same hypothesis shows that $N=\{1\}$. Hence $T$ is faithful.

This proves Assertion (1) of the proposition. To complete the proof, it suffices to show that the $\Gamma$-tree $T$ constructed above has the properties stated in Assertions (2) and (3). To prove Assertion (2), note that if $x$ is not in a conjugate of $\piggle(\frako_v)$, then according to Proposition \ref{don't move-a da lattice}, $x$ fixes no vertex of $\calt_v$. In particular, $x$ fixes no vertex of $T$, i.e. $x$ is $T$-hyperbolic.

To prove (3), suppose that $x\in\Gamma\cap\pizzle(\frako_v)$, and let $m$ denote the order of $\overline h_v(x)\in\pizzle(k_v)$. Since $x\in\Gamma\cap\pizzle(\frako_v)$, it follows from  Proposition \ref{don't move-a da lattice} that the vertex of $\calt$ representing $L_0$ is fixed by $x$. In particular $x$ is elliptic. Let $W$ denote the subset of $\Per_\calt(x)$ consisting of all vertices and edges whose $x$-period divides $mp^r$ for some $r\ge0$. To establish (3) it suffices to show that $W=\calt$. Since $x$ is elliptic we have $W\ne\emptyset$. It is therefore enough to show that if a vertex $s$ lies in $W$ then $\overline{\rm St}(s)\subset W$.

Let $\gamma\in\zzle(\frako_v)$ represent $x$.
Since $\overline h_v(x)\in\pizzle(k_v)$ has order $m$, we have $\overline h_v(\gamma^{mk})=\pm1$ for every positive integer $k$; in particular we have $\trace h_v(\gamma^{mk})=\pm 2$. Hence:
\Claim\label{linguine with claim sauce}
If $k$ is any positive integer, $\trace \gamma^{mk}$ lies in $\frako$ and is congruent to $\pm2$ modulo $\pi$. 
\EndClaim

Now let $s$ be any vertex in $W$. Choose an integer $r\ge0$ such that the $x$-period of $s$ divides $mp^r$. Thus if we set $A=\gamma^{mp^r}$, we have $s\in
\Fix_{\calt_v}( [A])$. Applying
\ref{linguine with claim sauce} with $k=p^r$ we find that 
$\trace A$ lies in $\frako$ and is congruent to $\pm2$ modulo $\pi$. It now follows from Corollary \ref{who isn't these days} that $ A^p=\gamma^{mp^{r+1}}$ fixes each vertex in the link of $s$, and therefore fixes every vertex or edge in the closed star of $s$. Hence $\overline{\rm St}(s)\subset W$, as required.
\EndProof

\section{Three-manifolds and trees}\label{tweetwee}

As in \cite{hakenmarg}, I will say that elements $x_1, \ldots, x_n$ of a group
$\Gamma$ are {\it independent} (or that the $n$-tuple $(x_1,\ldots,x_n)$ is independent) if $x_1,\ldots,x_n$ freely generate a free subgroup of $\Gamma$; and I will say that $x_1, \ldots, x_n$  are  {\it semi-independent} (or that the $n$-tuple $(x_1,\ldots,x_n)$ is semi-independent) if $x_1,\ldots,x_n$ freely generate a free semigroup in $\Gamma$. (In other words, $x_1, \ldots, x_n$
are semi-independent if distinct positive words in these elements
represent distinct elements of $\Gamma$.)

\Lemma\label{nielsen was here} If $u$ and $v$ are semi-independent elements of a group $\Gamma$, then $u$ and $uv$ are semi-independent.
\EndLemma

\Proof
We must show that if $W_1=W_1(x,y)$ and $W_2=W_2(x,y)$ are positive abstract words in letters $x$ and $y$, and if $W_1(u,uv)=W_2(u,uv)$, then $W_1=W_2$. For $i=1,2$, expanding the expression $W_i(u,uv)$ gives a positive word $V_i$ in $u$ and $v$. Note that $V_i$ ends with $u$ if $W_i(x,y)$ ends with $x$, and that that $V_i$ ends with $v$ if $W_i(x,y)$ ends with $y$. Since $u$ and $v$ are semi-independent, and since $V_1$ and $V_2$ define the same element of $\Gamma$, the words $V_1$ and $V_2$ are identical; hence $W_1$ and $W_2$ are either both empty, or both end in $x$, or both end in $y$. Now
I will argue by induction on $\length W_1+\length W_2$, the assertion of the lemma being trivial if $\length W_1+\length W_2=0$. If $\length W_1+\length W_2>0$, then $W_1$ and $W_2$ either both end with $x$ or both end with $y$. If they both end with $y$ we may write $W_i=U_iy$ for $i=1,2$; since $V_1=V_2$ we have $U_1(u,uv)uv=U_2(u,uv)uv$. Hence $U_1(u,uv)=U_2(u,uv)$. The induction hypothesis now gives $U_1=U_2$ and hence $W_1=W_2$. If 
$W_1$ and $W_2$ both end with $x$, the argument is similar.
\EndProof

\Number\label{phi edges}
I will be following the conventions of \cite[Section 4]{hakenmarg} concerning surfaces in $3$-manifolds, with one exception: The definition of incompressible surface given in \cite[Definition 4.3]{hakenmarg} includes the requirement that the surface be connected. Here I will define an {\it incompressible surface} in an orientable $3$-manifold $M$ to be a bi-collared surface $F\subset M$ such that every component of $F$ is incompressible in the sense defined in \cite{hakenmarg}. Thus connected incompressible surfaces in the sense of the present paper are the same as
incompressible surfaces in the sense of \cite{hakenmarg}.

If $F$ is an incompressible surface in
an orientable, irreducible $3$-manifold $M$, and if $\Gamma\cong\pi_1(M)$ denotes the group of deck transformations of the universal covering of $M$, I will denote by $T_F$ the dual
$\Gamma$-tree of $F$ in the sense of \cite[Subsection 4.2]{hakenmarg}. (The definition given there applies to any bicollared surface $F\subset M$, not necessarily connected.)

According to the definition given in \cite{hakenmarg}, $T_F$ is a quotient of the universal covering $\tM$ of $M$. If $p:\tM\to M$ denotes the covering projection, the surface $\tF:=p^{-1}(F)$ has a collar neighborhood $\tC\cong F\times[-1,1]$ such that the pre-image of each (open) edge under the quotient map $q:\tM\to T_F$ is a component of the set-theoretic interior of $\tC$; and the pre-image of each vertex under $q$ is a component of $\overline{\tM-\tC}$. Thus there is a natural bijective correspondence between components of $\tF$ and edges of $T_F$: to an edge $e$ there corresponds the component of $\tF$ contained in $q^{-1}(e)$. Likewise, there is a natural bijective correspondence between components of $\tM-\tF$ and vertices of $T_F$: to a vertex $s$ there corresponds the component of $\tM-\tF$ containing $q^{-1}(s)$. Since $q$ is $\Gamma$-equivariant, these bijective correspondences are also $\Gamma$-equivariant. A vertex $s$ is incident to an edge $e$ if and only if the component of $\tF$ corresponding to $e$ is contained in the closure of  the component of $\tM-\tF$ corresponding to $s$.

If $\Phi$ is a component of $F$, I will define a {\it $\Phi$-edge} of $T_F$ to be an edge corresponding to a component of $\tPhi:=p^{-1}(\Phi)\subset\tF$.
\EndNumber

\Proposition\label{obamium}
Let $M$ be
an orientable, irreducible $3$-manifold, and let $\Gamma\cong\pi_1(M)$ denote the group of deck transformations of the universal covering $\tM$ of $M$.
Suppose that $T$ is an essential $\Gamma$-tree. Then there exist an incompressible surface $F\subset M$ and a $\Gamma$-equivariant, nondegenerate simplicial map $f:T_F\to T$. 
\EndProposition

\Proof
Let $E$ denote the set of midpoints of edges of $T$. Let $p:\tM\to M$ denote the covering projection.
According to \cite[Proposition 1.3.8]{cgls} (see also \cite[Subsections 2.2 and 2.4]{handbook}), there is an incompressible surface $F$ which is ``associated'' (or ``weakly dual'') to the $\Gamma$-tree $T$ in the sense that there is a $\Gamma$-equivariant map $\phi:\tM\to |T|$, transverse to $E$, such that $\phi^{-1}(E)=p^{-1}(F)$.

Let $\calv$ and $\calv_F$ denote the vertex sets of $T$ and $T_F$. Let $T'$ denote the first barycentric subdivision of $T$. According to the discussion in \ref{phi edges}, each vertex $s\in\calv_F$ corresponds to a component $X$ of $\tM-\tF$. The set $f(X)$ is a connected subset of $|T|-E$, and is therefore contained in $\St_{T'}(f_\calv(s)$ for a unique vertex $f_\calv(s)$ of $T$. This defines a map $f_\calv:\calv_F\to\calv$, which is $\Gamma$-equivariant since both $\phi$ and the bijective correspondence defined in \ref{phi edges} are $\Gamma$-equivariant. 

Similarly, let  $\cale$ and $\cale_F$ denote the edge sets of $T$ and $T_F$. Each edge $e\in\cale_F$ corresponds to a component $Y$ of $\tF$, and $f(Y)$ is the midpoint of a unique edge $f_\cale(e)$ of $T$. This defines a map $f_\cale:\cale_F\to\cale$, which is $\Gamma$-equivariant since both $\phi$ and the bijective correspondence defined in \ref{phi edges} are. 

I claim that $f_\calv$ extends to a non-degenerate $\gamma$-equivariant simplicial map $f:T_F\to T$, which implies the conclusion of the proposition. It suffices to show that if $s_0$ and $s_1$ are the endpoints of an edge $e\in\cale_T$, then $f_\calv(s_0)$ and  $f_\calv(s_1)$ are the distinct endpoints of $f_\cale(e)$. Let $X_i$ denote the component of $\tM-\tF$ correspnding to $s_i$, and let $Y$ denote the component of $\tF$ correspnding to $e$. For $i=0,1$, since $s_i$ is incident to $e$, we have $Y\subset \overline{X_i}$ by \ref{phi edges}. Hence if $m$ denotes the midpoint of $f_\cale(e)$, we have $\{m\}=\phi(Y)\subset\phi( \overline{X_i})=\overline\St_{T'}(f_\calv(s_i)$, so that $f_\calv(s_i)$ is an endpoint of $f_\cale(e)$. If $f_\calv$ were to map $s_0$ and $s_1$ to the same endpoint $t$ of $f_\cale(e)$, then $\phi$ would map the neighborhood $X_0\cup Y\cup X_1$ of $Y$ into $\overline\St_{T'}(t)$, a subtree of $T'$ having $m$ as an endpoint. This is impossible by transversality.
\EndProof

\Proposition\label{gevaltium}
Let $F$ be an incompressible surface in
an orientable, irreducible $3$-manifold $M$, and let $\Phi$ be a component of $F$. Then there exists a $\Gamma$-equivariant simplicial map $g:T_F\to T_{\Phi}$ such that 
\begin{enumerate}
\item each $\Phi$-edge (\ref{phi edges}) of $T_F$  is mapped by $g$ onto an edge of $T_{\Phi}$;
\item if $e$ and $e'$ are distinct $\Phi$-edges of $T_F$, then $g(e)\ne g(e')$; and
\item each edge of $T_F$ which is not a $\Phi$-edge is mapped by $g$ to a vertex of $T_{\Phi}$.
\end{enumerate}
\EndProposition

\Proof 
Let $\calv_F$ and $\calv_\Phi$ denote the vertex sets of $T_F$ and
$T_\Phi$, and let $\cale_F$ and $\cale_\Phi$ denote their edge
sets. Let $\cale'_\Phi\subset\cale_F$ denote the set of $\Phi$-edges
of $T_F$. To prove the proposition, it suffices to define
$\Gamma$-equivariant maps $f_\calv:\calv_F\to\calv_\Phi$ and
$f_\cale:\cale_\Phi'\to\cale_\Phi$ such that (a) $f_\cale$ is
bijective, (b) for every edge $e$ of $T_F$ which is not a $\Phi$-edge,
$f_\calv$ maps the endpoints of $e$ to the same vertex of $T_\Phi$,
and (c) for every $e\in\cale'_\Phi$, the map
$f_\calv$ takes the endpoints of $e$ to the distinct endpoints of $f_\cale(e)$.

Let $\tM$ and $p:\tM\to M$ denote the universal covering of $M$ and
the covering projection. Set $\tF=p^{-1}(F)$ and $\tPhi=p^{-1}(\Phi)\subset\tF$.
According to the discussion in \ref{phi edges}, each vertex
$s\in\calv_F$ corresponds to a component $X$ of $\tM-\tF$. The
component of $\tM-\tPhi$ containing $X$ corresponds to a vertex of
$T_\Phi$, which I will denote $f_\calv(v)$. This
defines a map $f_\calv:\calv_F\to\calv$, which is $\Gamma$-equivariant
since the correspondences defined in \ref{phi
edges} are $\Gamma$-equivariant.

Similarly, each $\Phi$-edge $e$ of $T_F$ corresponds to a component $Y$ of
$\tPhi$, which also corresponds to a unique edge of $T_\Phi$; I will denote the latter edge by $f_\cale(e)$. This
defines a map $f_\calv:\calv_F\to\calv$, which is $\Gamma$-equivariant
since the bijective correspondence defined in \ref{phi
edges} is $\Gamma$-equivariant. and $f(Y)$ is the midpoint of a unique edge $f_\cale(e)$ of
$T$. This defines a map $f_\cale:\cale_F\to\cale$, which is
$\Gamma$-equivariant since the bijective
correspondences defined in \ref{phi edges} are $\Gamma$-equivariant.

Condition (a) follows from the bijectivity of the correspondences
defined in \ref{phi edges}. To prove Condition (b), let $e$ be an edge
of $T_F$ which is not a $\Phi$-edge, let $s_0$ and
$s_1$  denote its endpoints, let
$X_i$ denote the component of $\tM-\tF$ correspnding to $s_i$, and let
$Y$ denote the component of $\tF$ correspnding to $e$. By \ref{phi
edges} we have $Y\subset\overline{X_i}$ for $i=0,1$, and hence
$X_0\cup Y\cup X_1$ is connected. Since $e$ is not a $\Phi$-edge we
have $X_0\cup Y\cup X_1\subset\tM-\tPhi$. Hence $X_0$ and $X_1$ are contained in the same component of $\tM-\tPhi$, which implies that $f_\calv(s_0)=f_\calv(s_1)$. This proves (b).

To prove Condition (c), let $e$ be a $\Phi$-edge
of $T_F$, let $s_0$ and
$s_1$  denote its endpoints, Let
$X_i$ denote the component of $\tM-\tF$ correspnding to $s_i$, and let
$Y$ denote the component of $\tF$ corresponding to $e$. Since $e$ is a $\Phi$-edge,
$Y$ is a component of $\tPhi$. By \ref{phi
edges} we have $Y\subset\overline{X_i}$ for $i=0,1$, Hence if $Z_i$ denotes the component of $\tM-\tPhi$ containing $X_i$, we have $Y\subset\overline{Z_i}$. Since $f_\calv(s_i)$ is the vertex of $T_\Phi$ corresponding to $Z_i$, and $f_\cale(e)$ is the edge of $T_\Phi$ corresponding to $Z_i$, it follows from \ref{phi 
edges} that $f_\calv(s_i)$ is incident to $f_\cale(e)$ for $i=0,1$. To prove (c) it remains only to show that $f_\calv(s_0)\ne f_\calv(s_1)$, or equivalently that $Z_0\ne Z_1$. But $Y$ separates $\tM$ since $\tM$ is simply connected, and hence $X_0$ and $X_1$ are contained in distinct components of $\tM-Y$. Since $Z_i$ is disjoint from $Y$ and contains $X_i$ it follows that $Z_0\cap Z_1=\emptyset$.
\EndProof

Since a simplicial map takes vertices to vertices, we obtain:

\Corollary\label{jeez}
Let $F$ be an incompressible surface in
an orientable, irreducible $3$-manifold $M$, and let $\Phi$ be a component of $F$. Then every subgroup of $\Gamma$ which fixes a vertex of $T_F$ also fixes a vertex of $T_{\Phi}$.
\EndCorollary
\NoProof

The next three results, Propositions \ref{4.4gen}---\ref{4.9gen}, generalize Propositions 4.4, 4.6, 4.9, and 4.11 of \cite{hakenmarg}. They differ from the latter results only in that the incompressible surface $F$ appearing in each of the results below is permitted to be disconnected. 
%\redcomment{Revise this as needed.}

\Number\label{walk the line}
Following the terminology of \cite{hakenmarg}, I will say that the action of $\Gamma$ on $T$ is {\it linewise
  faithful} if for every line $L$ in $T$, the subgroup of $\Gamma$
that fixes $L$ pointwise is trivial.
\EndNumber

\Proposition\label{4.4gen}
Suppose that $F$ is an
incompressible surface in a compact orientable
$3$-manifold $M$.  Let $\Gamma$ denote the group of deck
transformations of $\tM$. Then
$T_F$ is an essential $\Gamma$-tree.
\EndProposition

\Proof
Choose a component $\Phi$ of $F$. If $\Gamma$ fixes a vertex  of $T_F$ then according to Corollary \ref{jeez}, it fixes a vertex of $T_{\Phi}$. But this contradicts \cite[Proposition 4.4]{hakenmarg}, which asserts that for a connected incompressible surface $\Phi$, the $\Gamma$-tree $T_{\Phi}$ is essential.
\EndProof

\Proposition\label{4.6gen}
Let $M = \HH^3/\Gamma$ be a closed
hyperbolic $3$-manifold containing an incompressible surface $F$, no component of which
is a fiber or a semi-fiber. Then for every
non-trivial $T_F$-elliptic element $\gamma\in\Gamma$, the diameter (as
an integer metric space, cf. \ref{tree stuff-not}) of the set of fixed
vertices $\gamma$ in $T_F$ is at most $14G-12n$, where $n$ is the number of components of $F$ and $G$ is the sum of their genera. In particular, the
action of $\Gamma$ on $T_F$ is linewise faithful.
\EndProposition

\Proof
The second assertion follows from the first, since a
line in a tree has infinite diameter. 

To prove the first assertion, let $F_1,\ldots,F_n$ denote the components of $F$, and let $g_i$ denote the genus of $F_i$ for $i=1,\ldots,n$. It follows from \cite[Corollary
 1.5]{DeB} that for every $i\in\{1,\ldots,n\}$ and every
non-trivial $T_{F_i}$-elliptic element $\gamma\in\Gamma$, the diameter of $\Fix_{T_{F_i}}(\gamma)$ in $T_{F_i}$ is at most $14g_i-12$. 

For each $i\in\{1,\ldots,n\}$ we may apply Lemma \ref{gevaltium} to obtain a $\Gamma$-equivariant simplicial map $g_i:T_{F}\to T_{F_i}$ such that Conditions (1)--(3) of Lemma \ref{gevaltium} hold with $\Phi=F_i$ and $g=g_i$.

Suppose that $\gamma\in\Gamma$ in $T_F$-elliptic. Then Corollary \ref{jeez}, applied to the subgroup $H=\langle\gamma\rangle$ of $\Gamma$ and the equivariant map  $g=g_i$, implies that $\gamma$ is $T_{F_i}$-elliptic for $i=1,\ldots,n$. We are required to prove that if $s$ and $t$ are vertices in $\Fix_{T_F}(\gamma)$ then $\dist(s,t)\le 14G-12n$, where $\dist$ denotes the distance function defined in \ref{tree stuff-not}. For each $i\in\{1,\ldots,n\}$, let $e^{i}_1,\ldots,e^{(i)}_{k_i}$ be the $F_i$-edges in the segment $[s,t]$ (see \ref{phi edges}). Conditions (1) and (2) of Lemma \ref{gevaltium} imply that $g_i$ maps the (open) edges $e^{i}_1,\ldots,e^{(i)}_{k_i}$ onto distinct edges of $T_{F_i}$, and Condition (3) implies that $g_i$ maps each component of $[s,t]-(e^{i}_1\cup\ldots\cup e^{(i)}_{k_i})$ to a vertex. Since $T_i$ is a tree, it follows that $g_i(\overline e^{i}_1\cup\cdots\cup \overline e^{(i)}_{k_i})$ is the segment with endpoints $g_i(s)$ and $g_i(t)$. Hence $k_i=\dist(g_i(s),g_i(t))$. Since 
$s,t\in \Fix_{T_F}(\gamma)$, the equivariance of $g_i$ gives
$g_i(s),g_i(t)\in \Fix_{T_{F_i}}(\gamma)$, and so $k_i\le\diam(\Fix_{T_{F_i}}(\gamma))$. Hence
$$\dist(s,t)=\length[s,t]=\sum_{i=1}^nk_i\le\sum_{i=1}^n\diam({\rm Fix}_{T_{F_i}}(\gamma))\le\sum_{i=1}^n(14g_i-12)=14G-12n.$$
\EndProof

\Proposition\label{4.9gen}
Let $F$ be a (possibly disconnected) incompressible surface in
an orientable $3$-manifold $M$.  Let $\Gamma$ denote the group of deck
transformations of the universal cover of $M$. Suppose that $\gamma$ is an infinite-order
element of $\Gamma$ such that $\Fix_{T_F}(\gamma)$ contains at
least one edge. Then for every integer $n>0$ we have
$\Fix_{T_F}(\gamma^n)=\Fix_{T_F}(\gamma)$.
\EndProposition

\Proof
This is proved in the same way as \cite[Proposition 4.9]{hakenmarg}. In the latter proof, the incompressibility of $F$ is used to show that the surface denoted $A_i$ is incompressible in the $3$-manifold denoted $N$. If $F$ is disconnected, the incompressibility of $A_i$ follows from the fact that the component of $F$ which is covered by the component $\widetilde\Phi_i$ of $\tF$ is incompressible.
\EndProof

\Lemma\label{h-a-double-r-i}Let $F$ be an incompressible surface in
an orientable hyperbolic $3$-manifold $M=\HH^3/\Gamma$.  Let $\gamma$ be an element of $\Gamma$ such that $\Per_{T_F}(\gamma)$ contains at least one edge in $T_F$. Then all edges of $\Per_{T_F}(\gamma)$ have the same $\gamma$-period.
\EndLemma

\Proof
We may assume that $\gamma\ne1$, so that $\gamma$ has infinite order. Let $m$ denote the smallest $\gamma$-period of any edge of $\Per_{T_F}(\gamma)$. Then $\Fix_{T_F}(\gamma^m)$ contains at least one edge of $T_F$. Applying Proposition \ref{4.9gen}, with $\gamma^m$ playing the role of $\gamma$ in that proposition, we find that $\Fix_{T_F}(\gamma^{mn})=\Fix_{T_F}(\gamma^{m})$ for every $n>0$. If $e$ is any edge of $\Per_{T_F}(\gamma)$ we have $\gamma^n\cdot e=e$ for some $n>0$, and hence $e\subset\Fix_{T_F}(\gamma^n)\subset\Fix_{T_F}(\gamma^{mn})=\Fix_{T_F}(\gamma^m)$. By the minimality of $m$ it follows that $e$ has $\gamma$-period exactly $m$.
\EndProof

In view of Lemma \ref{h-a-double-r-i}, it makes sense to introduce the

\Notation\label{beat it to death}
If $F$ is an incompressible surface in
an orientable hyperbolic $3$-manifold $M=\HH^3/\Gamma$, and if $\gamma$ is an element of $\Gamma$ such that $\Per_{T_F}(\gamma)$ contains at least one edge, we will denote by $m_F(\gamma)$ the common $\gamma$-period of all edges of $\Per_{T_F}(\gamma)$.

Note that $m_F(\gamma)$  is defined only if $\Per_{T_F}(\gamma)$ has at least one edge.

\EndNotation

\Proposition\label{ashes to ashes}
Let $M=\HH^3/\Gamma$ be a closed, orientable hyperbolic $3$-manifold, let $x$ and $y$ be non-commuting elements of $\Gamma$, and let $\calp$ be a set of (rational) primes. Suppose that there is a faithful, essential $\Gamma$-tree $T$ such that 
for every edge $\alpha$ of $\Per_{T}(x)$, all the primes dividing the $x$-period of $\alpha$ belong to $\calp$. Then either
\begin{enumerate}[(a)]
\item  at least one of the pairs $(x,yxy^{-1})$ and $(xy^{-1},y^{2})$ is independent; or
\item at least one of the pairs $(x^{-1},y)$ and $(x,y)$ is semi-independent; or
\item there is an incompressible surface $F\subset M$,  no component of which is not a fiber or semifiber, such that $\Per_{T_F}(x)$ has at least one edge, and at least one prime in $\calp$ divides $m_F(x)$.
\end{enumerate}
\EndProposition

\Remark 
The hypothesis that for every edge $e$ of $\Per_{T}(x)$, all the primes dividing the $x$-period of $e$ belong to $\calp$, may sometimes hold vacuously; that is, $\Per_{T}(x)$ may have no edges (and may be empty). The proposition applies in this case, even if $\calp=\emptyset$. 
\EndRemark

\Proof[Proof of Proposition \ref{ashes to ashes}]
According to Proposition \ref{obamium}, there exist an incompressible surface $F\subset M$ and a $\Gamma$-equivariant nondegenerate simplicial map $f:T_F\to T$. 

I claim:
\Claim\label{nattering nabobs} No component of $F$ is a fiber or semifiber. \EndClaim
To prove \ref{nattering nabobs}, consider an arbitrary component $F_0$ of $F$. Let $p:\HH^3\to M$ denote the quotient map, and choose a component $\tF$ of $p^{-1}(F_0)\subset\HH^3$. The definition of $T_F$ implies that $\tF$ is the pre-image under the $\Gamma$-equivariant quotient map $q:\HH^3\to T_F$ of an edge $e_0$ of $T_F$. By the equivariance of $q$, the stabilizer $\Gamma_{\tF}$ of $\tF$ in $\Gamma$ fixes $e_0$. By the equivariance of $f$, the group $\Gamma_{\tF}$ fixes the edge $f(e_0)\subset T$. In particular, $U:=\Fix_T(\Gamma_{\tF})$ is non-empty, and is therefore a subtree  of $T$. On the other hand, if
we identify $\Gamma$ with $\pi_1(M)$, then up to conjugacy $\Gamma_{\tF}$ is identified with the image of the inclusion homomorphism 
$\pi_1(F_0,\star)\to\pi_1(M,\star)$, where $\star\in F_0$ is any base point. Hence if $F_0$ is a fiber or semifiber, $\Gamma_{\tF}$ is a non-trivial normal subgroup of $\Gamma$. The normality of $\Gamma_{\tF}$ implies that $U$ is $\Gamma$-invariant. Since $T$ is minimal it follows that $U=T$. But this implies that the non-trivial subgroup $\Gamma_{\tF}$ acts trivially on $T$, in the sense that each of its elements acts by the identity. This contradicts the hypothesis that $T$ is faithful, and \ref{nattering nabobs} is proved.

In view of Proposition \ref{4.6gen}, it follows from \ref{nattering nabobs}
that the action of $\Gamma$ on $T_F$ is linewise faithful. Note also that $T_F$ is an essential $\Gamma$-tree by Proposition \ref{4.4gen}.

With this background, I will now show that one of the alternatives (a), (b) or (c) holds. I will divide the argument into three cases, each of which has two or more  (possibly overlapping) subcases.

{\bf Case I: $x$ is $T_F$-elliptic but $\Per_{T_F}(x)$ contains no edges.} Then in particular $x$ has a unique fixed vertex $s$ in $T_F$, and we have $\Fix_{T_F} (x^k)=\{s\}$ for every $n>0$. If we set $z=yxy^{-1}$ and $t=y\cdot s$, then $\Fix_{T_F} (z^k)=\{t\}$ for every $n>0$. If $s=t$, then since $T_F$ is essential, it now follows from \cite[Proposition 4.11]{hakenmarg} that $x$ and $z$ are independent in $\Gamma$. If $s\ne t$, it follows from \cite[Proposition 4.10]{hakenmarg} that $x$ and $z$ are independent in $\Gamma$. Thus Alternative (a) of the present proposition holds in both subcases. 

{\bf Case II: $x$ is $T_F$-hyperbolic.} In the subcase where $y$ is also $T_F$-hyperbolic, using the linewise faithfulness of the action of $\Gamma$ on $T_F$, we deduce from \cite[Proposition 3.6]{hakenmarg}
that $x$ and $y$ are semi-independent in $\Gamma$. In the subcase where $xy$ is $T_F$-hyperbolic, we deduce from \cite[Proposition 3.6]{hakenmarg}
that $x^{-1}$ and $xy$ are semi-independent in $\Gamma$. It then follows from Lemma \ref{nielsen was here} that $x^{-1}$ and $y$ are semi-independent in $\Gamma$. The same argument shows that in the subcase where $xy^{-1}$ is $T_F$-hyperbolic, $x^{-1}$ and $y^{-1}$ are semi-independent in $\Gamma$; this implies that $x$ and $y$ are semi-independent in $\Gamma$. Thus in all these subcases, Alternative (b) of the present proposition holds. There remains the subcase in which $y$, $xy$ and $xy^{-1}$ are all $T_F$-elliptic. In particular, the elements $xy^{-1}$ and  $y^2$, and their product $xy$, are all elliptic. Hence by \ref{mr fix it},  $xy^{-1}$ and  $y^2$ have a common fixed vertex. It then follows from \cite[Proposition 4.11]{hakenmarg} 
 that $xy^{-1}$ and  $y^2$ are independent. Thus alternative (a) holds in this subcase.

{\bf Case III: $x$ is $T_F$-elliptic and $\Per_{T_F}(x)$ has at least one edge.} Then according to Lemma \ref{h-a-double-r-i} and the definition in \ref{beat it to death}, the natural number $m:=m_F(x)$ is well defined.
In the subcase where some prime in $\calp$ divides $m$, Alternative (c) holds. 

Now consider the subcase in which no prime in $\calp$ divides $m$. Set $V=\Per_{T_F}(x)$ and $W=\Fix_T(x)$. According to the definition of $m=m_F(x)$, each edge of $V$ has $x$-period $m$. If $e$ is any edge of $V$, and if we set $\alpha=f(e)$, then since $x^m\cdot e=e$, it follows from equivariance that $x^m\cdot \alpha=\alpha$. Hence if $d$ denotes the $x$-period of the edge $\alpha$, we have $d|m$. But by hypothesis each prime dividing $d$ belongs to $\calp$. Since no prime in $\calp$ divides $m$, we must have $d=1$. This shows that $f(V)\subset W$. By equivariance it follows that $f(y\cdot V)\subset y\cdot W$. If we again set $z=yxy^{-1}$, then $y\cdot V=\Per_{T_F}(z)$, and $y\cdot W=\Fix_T(z)$. We now distinguish two sub-subcases, depending on whether $V\cap y\cdot V$ is empty or not. If $V\cap y\cdot V=\emptyset$, so that $\Per_{T_F}(x)\cap \Per_{T_F}(z)=\emptyset$, it follows from \cite[Proposition 4.10]{hakenmarg} that $x$ and $z$ are independent in $\Gamma$, i.e. Alternative (a) holds. If $V\cap y\cdot V\ne\emptyset$, then 
$${\rm Fix}_T(x)\cap{\rm Fix}_T(z)=W\cap y\cdot W\supset f(V)\cap f(y\cdot V)\supset f(V\cap y\cdot V)\ne\emptyset.$$ 
It then follows from 
\cite[Proposition 4.11]{hakenmarg} that $x$ and $z$ are independent in $\Gamma$, and again Alternative (a) holds.
\EndProof

\section{Displacements, volumes and indices of two-generator subgroups}\label{better bounds}
The main results of this section, Proposition \ref{radical chick} and Corollary \ref{ranana}, are refinements of Theorem 4.2 and Corollary 4.3 of \cite{bounds}.

\Proposition\label{hunt it down and kill it}
Let $x$ and $y$ be non-commuting elements of 
$\pizzle(\CC)$ such that $\langle x,y\rangle$ is discrete and torsion-free and has infinite covolume. Then for every $P\in\HH^3$ we have
$$\frac1{1+\exp{d_P(x)}}+\frac1{1+\exp{d_P(y)}}\le\frac12.$$
In particular we have $\max(d_P(x),d_P(y))\ge\log3$.
\end{proposition}
 
\Proof
Set $\Gamma=\langle x,y\rangle$. Proposition 3.14 of \cite{bounds}, applied to the hyperbolic $3$-manifold $M=\HH^3/\Gamma$, shows that $\Gamma$ is free. The conclusion now follows from the case $k=2$ of \cite[Theorem 4.1]{surgery}. \EndProof

\Proposition\label{radical chick}
Let ${\alpha}$ and ${\beta}$ be positive real numbers such that
$$\frac1{1+\exp{\alpha}}+\frac1{1+\exp{\beta}}>\frac12.$$
 Then there is a constant $V_{{\alpha},{\beta}}$ with the following property:
\begin{itemize}
\item Let  $M=\HH^3/\Gamma$ be any orientable hyperbolic $3$-manifold with
$\infty\ge\vol M> V_{{\alpha},{\beta}}$. Let $P$ be any point of $ \HH^3$, and let $x$ and $y$  be elements of $\Gamma$ such that  $d_P(x)\le\alpha$ and $d_P(y)\le\beta$. Then $x$ and $y$ commute in $\Gamma$.
\end{itemize}
\EndProposition

\Proof
Assume that the assertion is false. Then there is a sequence $(M_n)_{n\ge1}$ of 
orientable finite-volume hyperbolic $3$-manifolds, where $M_n=\HH^3/\Gamma_n$, and for each $n$ a pair of non-commuting elements
$x_n,y_n\in\Gamma_n$   and a point $P_n\in\HH^3$,   such that 
$d_{P_n}(x_n)\le\alpha$ and $d_{P_n}(y_n)\le\beta$ for each $n$.

After replacing each $\Gamma_n$ by a suitable conjugate of itself in
$\pizzle(\CC)$, we may assume that the $P_n$ are all the same point of $\HH^3$,    which I will denote by $P$.    Thus for each $n$ we have
\begin{equation}\label{mindgymillion}
d_{P}(x_n)\le\alpha\text{ and }d_{P}(y_n)\le\beta.
\end{equation}

For each $n$, set $\tGamma_n:=\langle x_n,y_n\rangle$. Note that $\tGamma_n$ is discrete and torsion-free since $\Gamma_n$ is, and that $\tGamma_n$ is non-abelian---and hence non-elementary by \cite[Proposition 2.1]{finiteness}---since $x_n$ and $y_n$ do not commute. Set $\tM_n:=\HH^3/\tGamma_n$.

Since 
$$\frac1{1+\exp{\alpha}}+\frac1{1+\exp{\beta}}>\frac12,$$
it follows from (\ref{mindgymillion}) that
\Equation\label{pathetic like a fox}
\frac1{1+\exp{d_P(x_n)}}+\frac1{1+\exp{d_P(y_n)}}>\frac12.
\EndEquation
It therefore follows from Proposition \ref{hunt it down and kill it} that $\vol\tM_n<\infty$. On the other hand, $\tM_n$ covers $M_n$, and hence $\vol \tM_n\ge\vol M_n$.  In particular, $\vol\tM_n\to\infty$.

It follows from (\ref{mindgymillion}) that the $x_n$ and $y_n$ lie in a compact subset of $\pizzle(\CC)$. Hence, after passing to a subsequence, we may assume that the sequences $(x_n)$ and $(y_n)$ converge in $\pizzle(\CC)$ to limits $x_\infty$ and $y_\infty$. 
It then follows, again  from (\ref{mindgymillion}), that
\begin{equation}\label{manx}
\frac1{1+\exp{d_P(x_\infty)}}+\frac1{1+\exp{d_P(y_\infty)}}\ge
\frac1{1+\exp{\alpha}}+\frac1{1+\exp{\beta}}>\frac12.
\end{equation}

For $1\le n\le\infty$ we define a representation $\rho_n$ of the rank-$2$ free group $F_2=\langle\xi,\eta\rangle$ by $\rho_n(\xi)=x_n$, $\rho_n(\eta)=y_n$.  Thus $\rho_n(F_2)=\tGamma_n$ for each $n$. Since $\xi$ and $\eta$ generate $F_2$, and since $\rho_n(\xi)\to\rho_\infty(\xi)$ and $\rho_n(\eta)\to\rho_\infty(\eta)$ as $n\to\infty$, we have $\rho_n(\gamma)\to\rho_\infty(\gamma)$ for every $\gamma\in F_2$. By definition (see, for example, Section 2 of \cite{finiteness}) this means that the sequence $(\rho_n)$ converges algebraically to $\rho_\infty$.   Set $x_\infty=\rho_\infty(\xi)$, $y_\infty=\rho_\infty(\eta)$.

Let $D$ denote the set of representations of $F_2$ in $\pizzle(\CC)$ whose  images are discrete, torsion-free, and non-elementary.
According to \cite[Theorem 2.4]{finiteness} (a theorem essentially due to T. Jorgensen and P. Klein \cite{jk}), the limit of any algebraically convergent sequence of representations in $D$ is again in $D$.  Hence $\rho_\infty\in D$. Thus $\tGamma_\infty:=\rho_\infty(F_2)=\langle x_\infty,y_\infty\rangle$ is a discrete group.

According to \cite[Proposition 3.8]{jm}, since $(\rho_n)$ converges algebraically, the sequence of discrete groups $(\tGamma_n) $ has a geometrically convergent subsequence (in the sense defined in \cite{jm}). Hence without loss of generality we may assume that $(\tGamma_n) $ converges geometrically to some discrete group $\hGamma_\infty$. It then follows, again from \cite[Proposition 3.8]{jm},  that 
$\tGamma_\infty\le \hGamma_\infty$.

According to \cite [Lemma 4.1]{bounds}, there is a neighborhood $W$ of the identity in $\pizzle(\CC)$ such that $\tGamma_n\cap W=\{1\} $ for every $n\in\NN$. Let $E$ denote the set of all torsion-free subgroups $\Delta$ of $\pizzle(\CC)$ such that $\Delta\cap W=\{1\} $. (In particular each group in $E$ is discrete.) According to \cite[Theorem 1.3.1.4]{ceg}, $E$ is compact in the topology of geometric convergence. Since $\tGamma_n\in E$ for each $n\in\NN$, we have $\hGamma_\infty\in E$. In particular $\hGamma_\infty$ is torsion-free. We let $\hM_\infty$ denote the orientable hyperbolic $3$-manifold $\HH^3/\hGamma_\infty$.

Since $(\tGamma_n) $ converges geometrically to $\hGamma_\infty$, the sequence of orientable hyperbolic $3$-manifolds $(\tM_n)$ converges geometrically to $\hM_\infty$ in the sense of \cite[Chapter E]{bp}. If $\vol \hM_\infty$ were finite, it would then follow from \cite[Proposition E.2.5]{bp} that the sequence $(\vol\tM_n)$ had the finite limit $\vol \hM_\infty$, which contradicts 
$\vol\tM_n\to\infty$. Thus $\tM_\infty:=\HH^3/\rho_\infty(F_2)$ is a hyperbolic $3$-manifold of infinite volume. It therefore follows from Proposition \ref{hunt it down and kill it} that
$$
\frac1{1+\exp{d_P(x_\infty)}}+\frac1{1+\exp{d_P(y_\infty)}}\le\frac12.
$$
But this contradicts (\ref{manx}).
\EndProof

\Corollary\label{ranana}
Let ${\alpha}$ and ${\beta}$ be positive real numbers such that
$$\frac1{1+\exp{\alpha}}+\frac1{1+\exp{\beta}}>\frac12.$$
 Then there is a natural number $\delta_{{\alpha},{\beta}}$ with the following property:
\begin{itemize}
\item Let  $\Gamma$ be any cocompact, discrete, torsion-free subgroup of $\pizzle(\CC) $. Let $P$ be any point of
$ \HH^3$, and let $x$ and $y$  be elements of $\Gamma$ such that  $d_P(x)\le\alpha$ and $d_P(y)\le\beta$.
Then either  $x$ and $y$ commute in $\Gamma$, or the subgroup $\langle x,y\rangle$ has index at most $\delta_{{\alpha},{\beta}}$ in $\Gamma$.
\end{itemize}
\EndCorollary

\Proof
Let $v$ denote the infimum of the volumes of all
hyperbolic $3$-manifolds;  we have $v>0$, for example by
\cite[Theorem 1]{meyerhoff}. Let $V_{\alpha,\beta}$ be a positive real number having the property stated in Proposition \ref{radical chick}, and set 
$$\delta_{\alpha,\beta}=\lfloor \frac{V_{\alpha,\beta} }v\rfloor.$$

Suppose that $\Gamma\le\pizzle(\CC) $ is cocompact, discrete and torsion-free. Set
$M=\HH^3/\Gamma$. Let $P$ be any point of
$ \HH^3$, and let $x$ and $y$  be elements of $\Gamma$ such that  $d_P(x)\le\alpha$ and $d_P(y)\le\beta$.
Then either  $x$ and $y$ commute in $\Gamma$, or the subgroup $\langle x,y\rangle$ has index at most $\delta_{{\alpha},{\beta}}$ in $\Gamma$. Set $\tGamma=\langle x,y\rangle$ and $\tM=\HH^3/\tGamma$. Applying 
Proposition \ref{radical chick}, with $\tM$ and $\tGamma$ in place of $M$ and $\Gamma$, we find that either $x$ and $y$ commute in $\tGamma\le\Gamma$, or $\vol \tM\le V_{{\alpha},{\beta}}<\infty$. In the latter case, since $\vol M\ge v$, we find that
$$|\Gamma:\tGamma]=\frac{\vol\tM}{\vol M}\le
\frac{V_{{\alpha},{\beta}}}v.$$
It follows that $|\Gamma:\tGamma|\le \delta_{{\alpha},{\beta}} $.
\EndProof

\section{Haken manifolds, measures and displacements}\label{hakenmeasure}

This section is devoted to some refinements of results from \cite{hakenmarg} which will be needed in the next section. for the proof of Theorem \ref{N lives} (Theorem B of the Introduction).

The following result, Proposition \ref{gen6.2}, generalizes Proposition 6.2 of \cite{hakenmarg}. As in the latter result, if $\gamma$ is an element of $\pizzle(\CC)$ I will denote by $\gamma_\infty$ the M\"obius transformation of the sphere at infinity $S_\infty$ defined by $\gamma$. I will use
the same notation as in \cite{paradoxical}.  In particular, the
conformal expansion factor (\cite [2.4]{paradoxical}) of the M\"obius transformation
$\gamma_\infty$ associated to the point $z\in\HH^3$ will be denoted
$\lambda_{\gamma,z}$; the pull-back (\cite[3.1]{paradoxical}) of a
measure $\mu$ under $\gamma_\infty$ will be denoted
$\gamma^*\mu$; and $\mathcal A = (A_z)$ will denote the area density on
$S_\infty$ (see \cite [3.3]{paradoxical}).

\Proposition\label{gen6.2}
 Let $M = \HH^3/\Gamma$ be a
closed hyperbolic $3$-manifold and let $P$ be a point in $\HH^3$.
Let $T$ be a $\Gamma$-tree and let $x$ and $y$  be $T$-elliptic
elements of $\Gamma$ such that $\Fix_T(x) \cap \Fix_T(y) =
\emptyset$.  Let $n$ be a positive integer, and suppose that there is no integer $k$ with $0<k\le2n$ such that either $x^k$ or $y^k$ has a fixed edge in $T$. Then
there exist Borel measures $\sigma_i$ and $\tau_i$ on $S_\infty$, for $i=-n,\ldots,0,\ldots,n$, such that 
\begin{enumerate}
\item $\sigma_0 + \tau_0 \le A_P$;
\item $\sum_{i=1}^n(
\sigma_i + \sigma_{-i} )\le \tau_0$ and  $\sum_{i=1}^n(\tau_i + \tau_{-i}) \le \sigma_0$;  and

\item $\int_{S_\infty}\lambda^2_{x^{\epsilon}, P}d\sigma_{\epsilon i}=\sigma_{\epsilon (i-1)}(S_\infty)$
  and 
$\int_{S_\infty}\lambda^2_{y^{\epsilon}, P}d\tau_{\epsilon i}=\tau_{\epsilon (i-1)}(S_\infty)$
for $i=1,\ldots,n$ and $\epsilon=\pm1$.
\end{enumerate}

Furthermore, in the case that $T=T_F$ for some
incompressible surface $F$ in $M$,
we have $\sigma_0 + \tau_0 = A_P$.
\EndProposition

\Proof
According to \cite[Proposition 6.1]{hakenmarg}, there exist an edge $e$ of $T$ and disjoint subsets $X$ and $Y$ of $\Gamma$ such that
\begin{itemize}
\item $\Gamma$ is the disjoint union of $X$, $Y$ and $\Gamma_e$;
\item $x^{\pm k}  Y \subset X$ and $y^{\pm k}  X
  \subset Y$ for $0 < k \le n$; and
\item $x^i  Y \cap x^j  Y = \emptyset$ and
$y^iX \cap y^jX = \emptyset$
for any pair of  distinct integers $i$ and $j$
 with $-n\le i \le n$ and $-n\le j\le n$.
\end{itemize}

Since $\Gamma$ is
discrete, the set $\Gamma\cdot P$ is uniformly discrete, in the sense
of \cite [Subsection 4.1]{paradoxical}.

Set 
$$X' = X -\bigcup_{i=1}^n (xY\cup x^{-1}Y)\quad\text{ and }\quad Y' = Y- \bigcup_{i=1}^n (yX\cup y^{-1}X).$$
Define a subset $\cals$ of the power set of $\Gamma\cdot P$ by
$$\cals=\{X'\cdot P,Y'\cdot P\, \Gamma_e\cdot P\}\cup\bigcup_{i=1}^n\{x^iY\cdot P,\; x^{-i}Y\cdot P,\;
y^iX\cdot P,\; y^{-i}X\cdot P\}.$$
Define $\mathfrak V$ to be the collection of
all unions of sets in $\cals$.

We apply \cite[Proposition 4.2]{paradoxical} with $W = \Gamma\cdot
P$, to construct a family $(\calm_{V})_{V\in \mathfrak V}$, of
$D$-dimensional conformal densities, for some $D\in [0,2]$, such that
conditions (i)-(iv) of \cite[Proposition 4.2]{paradoxical} hold.  
For $-n\le i\le n$, set 

$$\calm_{x^{i}Y\cdot P} = (\sigma^{(i)}_z) _{z\in\HH^3}, \quad{\rm and}\quad
\calm_{y^iX\cdot P} = (\tau^{(i)}_z) _{z\in\HH^3}.$$

It follows
from conditions (i) and (ii) of \cite[Proposition 4.2]{paradoxical}
that $\calm_\Gamma = \calm_{X\cdot P} + \calm_{Y\cdot P} + \mathcal
M_{\Gamma_e\cdot P}$ is a $\Gamma$-invariant conformal density.  Since
$M = \HH^3/\Gamma$ is a closed manifold, every $\Gamma$-invariant
superharmonic function on $M$ is constant.  Thus by
\cite [Proposition 3.9]{paradoxical}, $D=2$ and $\calm_{X\cdot P} +
\mathcal M_{Y\cdot P} + \calm_{\Gamma_e\cdot P} = k\mathcal A$ for
some constant $k$.  Condition (i) of \cite[Proposition
  4.2]{paradoxical} guarantees that $k > 0$.  Thus by normalizing the family
$(\mathcal M_{V})_{V\in \mathfrak V}$ appropriately we may assume that
$k=1$.

Now define $\sigma_i = \sigma^{(i)}_P$ and $\tau_i = \tau^{(i)}_P$ for $-n\le i\le n$. Then Conclusion (1) follows from the equality
$\calm_{X\cdot P} + \mathcal M_{Y\cdot P} + \calm_{\Gamma_e\cdot P} =
\mathcal A$ by specializing to $z=P$.

According to \cite[Proposition 4.2 (ii)]{paradoxical}, we have 
$$\calm_{X\cdot P} = \calm_{X'\cdot P} +\sum_{i=1}^n \calm_{x^iY\cdot P}
+\sum_{i=1}^n M_{x^{-i}Y\cdot P}$$
and 
$$\calm_{Y\cdot P} = \calm_{Y'\cdot P} +
\sum_{i=1}^n\calm_{y^iX\cdot P} +\sum_{i=1}^n \calm_{y^{-i}X\cdot P}.$$
Specializing to $z=P$, we obtain
Conclusion (2).

The proofs of Conclusions (3) and (4) are based on
\cite[Proposition 4.2 (iii)]{paradoxical}. Let $\epsilon\in\{1,-1\}$ be given. Applying \cite[Proposition 4.2 (iii)]{paradoxical} with $\gamma=x^\epsilon$ and with $V=x^iY$ for any $i\in\{1,\ldots,n\}$, we find that
$$(x^\epsilon)^*(\calm_{x^{\epsilon i}Y}) = \calm_{x^{\epsilon(i-1)}Y}.$$  Since $\calm_{x^{\epsilon i}Y\cdot P}
 = (\sigma^{(\epsilon i)}_z) _{z\in\HH^3}$
and that $\calm_{x^{\epsilon (i-1)}Y\cdot P}
 = (\sigma^{(\epsilon (i-1) )}_z) _{z\in\HH^3}$,
the
definition of the pull-back (see \cite[3.4.1]{paradoxical}) gives
that
$$(x^\epsilon)^*_\infty \sigma^{(\epsilon i)}_z=
\sigma^{ (\epsilon (i-1)) }_z$$ 
which implies that
$$d\sigma^{ (\epsilon (i-1)) }_z= \lambda^2_{x^\epsilon,z}\,d\sigma^{(\epsilon i)}_z.$$

Taking $z=P$
and integrating over $S_\infty$ we obtain
$$\int_{S_\infty}
\lambda^2_{x^\epsilon,P}\,d\sigma_{\epsilon i}
=
 \int_{S_\infty}d\sigma_{ \epsilon (i-1) }=\sigma_{ \epsilon (i-1) }
(S^\infty)$$
which is the first part of Conclusion (3). The second part is proved in exactly the same way, using $X$ and $y$ in place of $Y$ and $x$.

The final assertion of the proposition, that $\sigma_0 + \tau_0 = A_P$ in the case where $T=T_F$ for some
incompressible surface $F$ in $M$, is proved in exactly the same way as 
the final assertion of [Proposition 6.2] of \cite{hakenmarg}. Note that the latter proof depends on the fact that neither $x^2$ nor $y^2$ has a fixed edge in $T_F$; in the context of the present proof, this is a consequence of the hypothesis that 
there is no integer $k$ with $0<k\le2n$ such that either $x^k$ or $y^k$ has a fixed edge in $T$.
\EndProof

\Proposition\label{toil and trouble} Let $M = \HH^3/\Gamma$ be a closed
hyperbolic $3$-manifold and let $P$ be a point of $\HH^3$. Let $F$ be an incompressible surface in
$M$, and let $x$ and $y$ be $T_F$-elliptic elements of $\Gamma$ such
that $\Fix_{T_F}(x) \cap \Fix_{T_F}(y) = \emptyset$. Let $n$ be a positive integer, and suppose that there is no integer $k$ with $0<k\le2n$ such that either $x^k$ or $y^k$ has a fixed edge in $T_F$.  Set $D_x=\exp(2 d_P(x))$ and
$D_y=\exp(2 d_P( y))$. Then there exist strictly positive real numbers $\alpha_i$ and $\beta_i$, for $i=-n,\ldots,0,\ldots,n$, such that the following conditions hold:
\begin{enumerate}
\item $\alpha_0 +\beta_0=1$;
\vspace{6pt plus 3pt}
\item $\sum_{i=1}^n(
\alpha_i + \alpha_{-i} )\le \beta_0$ and  $\sum_{i=1}^n(\beta_i + \beta_{-i}) \le \alpha_0$; and
\vspace{6pt plus 3pt}
\item$\displaystyle\frac{\alpha_{\epsilon (i-1)}  (1-\alpha_{\epsilon i})}{\alpha_{\epsilon i}(1-
\alpha_{\epsilon (i-1)} )} \le D_x$\   and 
$\displaystyle\frac{\beta_{\epsilon (i-1)}  (1-\beta_{\epsilon i})}{\beta_{\epsilon i}(1-
\beta_{\epsilon (i-1)} )} \le D_y$
for $i=1,\ldots,n$ and $\epsilon=\pm1$. 
\end{enumerate} 
\EndProposition

\Proof
Set $T=T_F$, so that the hypotheses of Proposition \ref{gen6.2} hold. For $i=-n,\ldots,0,\ldots,n$
let $\sigma_i$ and $\tau_i$ denote the Borel measures on $S_\infty$ given by Proposition \ref{gen6.2}. Let $\alpha_i$ and $\beta_i$ denote the total masses of $\sigma_i$ and $\tau_i$ respectively.

%\redcomment{The statement of Prop. \ref{gen6.2} has changed---$\tau_0$ and $\nu_0$ are reversed and $\nu$ is replaced by $\sigma$. Fix defs. of $\alpha_i$ and $\beta_i$ and/or the statement of Prop. \ref{toil and trouble} accordingly, and make the discussion below fit.}

Since $T=T_F$, the last sentence of Proposition \ref{gen6.2} gives that $\sigma_0+\tau_0=A_P$, which implies Conclusion (1) of the present proposition. 
Conclusion (2) of the present proposition follows from Conclusion (2) of Proposition \ref{gen6.2}. 

As a preliminary to proving Conclusion (3)  of the present proposition, I will show that 
\Equation\label{the ninth circle}
0<\alpha_i<1\text{ and }0<\beta_i<1\text{ for }i=-n,\ldots,0,\ldots,n.
\EndEquation
To prove (\ref{the ninth circle}), first note that since the function $\lambda^2_{x^\epsilon,P}$ are continuous and positive-valued on $S^2$, it follows from  Conclusion (3) of Proposition \ref{gen6.2} that for a given  $i\in\{-(n-1),\ldots,0,\ldots,n\}$ we have 
\Equation\label{the eighth circle}
\alpha_i=0\text{ if and only if }\alpha_{i-1}=0.
\EndEquation
The same argument shows that 
\Equation\label{the seventh circle}
\beta_i=0\text{ if and only if }\beta_{i-1}=0.
\EndEquation
Now suppose that for some $i_0\in\{-n,\ldots,0,\ldots,n\}$
we have $\alpha_{i_0}=0$. Then by (\ref{the eighth circle}) we have $\alpha_{i}=0$ for every $i\in\{-n,\ldots,0,\ldots,n\}$. In particular $\alpha_0=0$, and hence by Conclusion (2) of Proposition \ref{gen6.2}, we have $\beta_i=\beta_{-i}=0$ for $i=1,\ldots,n$. In particular $\beta_1=0$, and hence by (\ref{the seventh circle}) we have $\beta_0=0$. We now have $0+0=\alpha_0+\beta_0=1$, a contradiction. 

This shows that $\alpha_i>0$ for $i=-n,\ldots,0,\ldots,n$. The same argument shows that $\beta_i>0$ for $i=-n,\ldots,0,\ldots,n$. Since $\alpha_0+\beta_0=1$ it follows that $\alpha_0<1$ and $\beta_0<1$. It then follows from Conclusion (2) of Proposition \ref{gen6.2} that $\alpha_i<1$ and $\beta_i<1$ for $i=-n,\ldots,0,\ldots,n$. Thus (\ref{the ninth circle}) is established.

Conclusion (3)  of the present proposition will be deduced from Conclusions (1)---(3) of Proposition \ref{gen6.2} via \cite[Lemma 2.1]{hakenmarg}. Let $i\in\{1,\ldots,n\}$ and $\epsilon\in\{+1,-1\}$ be given. To facilitate the transition between the notation of the present argument and that of \cite[Lemma 2.1]{hakenmarg}, let $\nu$ denote the measure $\sigma_{\epsilon i}$, let $\gamma$ denote the loxodromic isometry $x^\epsilon$ of $\HH^3$, and set $a=\alpha_{\epsilon i}$ and $b=\alpha_{\epsilon(i-1)}$.
By (\ref{the ninth circle}) we have $0<a<1$ and $0<b<1$.
It follows from
Conclusions (1) and (2) of Proposition \ref{gen6.2} that $\sigma_{\epsilon i}\le\tau_0\le A_P$, so that
\Equation\label{A1A}\nu\le A_P.\EndEquation
The definition of $\alpha_{\epsilon i}$ gives $\sigma_{\epsilon i}(S_\infty)=\alpha_{\epsilon i}$, i.e.
\Equation\label{A2A}\nu(S_\infty)=a.\EndEquation
 Conclusion (3) of Proposition \ref{gen6.2} gives $\int_{S_\infty}\lambda^2_{x^{\epsilon}, P}d\sigma_{\epsilon i}=\sigma_{\epsilon (i-1)}(S_\infty)=\alpha_{\epsilon (i-1)}$, i.e.
\Equation\label{A3A}
\int_{S_\infty}\lambda^2_{\gamma, P}d\nu=b.
\EndEquation
 According to \cite[Lemma 2.1]{hakenmarg}, if  (\ref{A1A}), (\ref{A2A}) and (\ref{A3A}) for numbers $a,b\in(0,1)$, then we have
$$\dist(z,\gamma \cdot z) \ge \frac12\log\frac{b(1-a)}{a(1-b)}.$$   
Hence this inequality holds with the choices of $a$ and $b$ made above; it immediately implies the first inequality of Conclusion (3)  of the present proposition. To prove the second inequality we take $\nu=\tau_{\epsilon i}$, $\gamma=y^\epsilon$, $a=\beta_{\epsilon i}$ and $b=\beta_{\epsilon(i-1)}$ and reason in exactly the same way.
\EndProof

\Lemma\label{it's da bomb}
Let $M = \HH^3/\Gamma$ be a closed
hyperbolic $3$-manifold and let $P$ be a point of $\HH^3$.  Let $F$ be an incompressible surface in
$M$, and suppose that $x$ and $y$ are $T_F$-elliptic elements of $\Gamma$ such
that $\Fix_{T_F}(x) \cap \Fix_{T_F}(y) = \emptyset$. Suppose that there is no integer $k$ with $0<k<5$ such that either $x^k$ or $y^k$ has a fixed edge in $T_F$.  Then either $ d_P(x)>0.34$ or $ d_P(y)>1.08$.
\EndLemma

\Proof
I will apply Proposition \ref{toil and trouble} with $n=2$. This gives numbers $\alpha_j$ and $\beta_j$ for $j=-2,-1,0,1,2$ such that Conditions (1)--(4) of Proposition \ref{toil and trouble} hold. I will distinguish several cases.

First
consider the case in which $\alpha_0>0.35$ and $\alpha_\epsilon\le0.214$ for some $\epsilon\in\{1,-1\}$. In this case we have
$$D_x\ge\frac{\alpha_0 (1-\alpha_\epsilon)}{\alpha_\epsilon(1-\alpha_0 )} 
\ge\frac{0.35(1-0.214)}{0.214(1-0.35 )}=1.9777 \ldots,$$
which gives $ d_P(x)\ge0.3409\ldots$, implying the conclusion in this case. Next consider the case in which $\alpha_0>0.35$ and $\min(\alpha_1,\alpha_{-1})>0.214$. In this case we have $$\alpha_2+\alpha_{-2}\le\beta_0-(\alpha_1+\alpha_{-1})=1-\alpha_0-(\alpha_1+\alpha_{-1})<1-0.35-2\cdot0.214=0.222$$
and hence $\alpha_\epsilon<0.111$ for some $\epsilon\in\{1,-1\}$. It follows that 
$$D_x\ge\frac{\alpha_{\epsilon}(1-\alpha_{2\epsilon})}{\alpha_{2\epsilon}(1-\alpha_{\epsilon} )} 
\ge\frac{0.214(1-0.111)}{0.111(1-0.214 )}=2.1805 \ldots.$$ This gives $ d_P(y)\ge0.389\ldots$, implying the conclusion in this case.

There remains the case in which $\alpha_0\le0.35$ and hence $\beta_0=1-\alpha_0\ge0.65$. In this case we have
$\beta_1+\beta_{-1}\le\alpha_0\le0.35$, and hence $\beta_\epsilon\le0.175$ for some $\epsilon\in\{1,-1\}$. Hence
$$D_y\ge\frac{\beta_0 (1-\beta_\epsilon)}{\beta_\epsilon(1-\beta_0 )} 
\ge\frac{0.65(1-0.175)}{0.175(1-0.65 )}=8.755 \ldots.$$
This gives $ d_P(x)\ge
1.084 \ldots$,
and the conclusion is established in all cases.
\EndProof

\section{Displacements and finitistic orders}\label{bee section}

This section contains the proof of one of the main results of the paper, Theorem \ref{N lives}, which was stated in the Introduction as Theorem B.

\Proposition\label{hempel says}Let $\Gamma$ be a torsion-free group whose abelian subgroups are all cyclic. Then every solvable subgroup of $\Gamma$ is cyclic.
\EndProposition

\Proof
Let $\Delta$ be a solvable subgroup of $\Gamma$, and let $\Delta=\Delta_0\triangleright\Delta_1\triangleright\cdots\triangleright\Delta_n=\{1\}$ be its derived series. If $n\le1$ then $\Delta$ is abelian and hence cyclic. Now suppose that $n>1$. Set $X=\Delta_{n-2}$ and $Y=\Delta_{n-1}$.  Then $ Y\le X$ is abelian, while $X$ is not. By Zorn's lemma, $ Y$ is contained in a maximal abelian subgroup $A$ of $ X$. Since $A$ contains the commutator subgroup $ Y$ of $ X$, it is normal in $ X$. But $A\le\Gamma$ is abelian and therefore infinite cyclic; hence if $C=C(A)\cap X$ denotes the centralizer of $A$ relative to $X$, we have $| X:C|\le2$. For any $c\in C$, the subgroup  $\langle A\cup\{c\}\rangle$ of $X$ is abelian. By maximality we must have $c\in A$. Hence $C=A$, so that $|X:A|\le2$. Since $X$ contains an infinite cyclic subgroup with index $2$, it is either dihedral or infinite cyclic. The former alternative is ruled out by the hypothesis that $\Gamma$ is torsion-free, while the latter alternative is ruled out by the fact that $X$ is non-abelian.
\EndProof

\Theorem\label{N lives} 
There exists a natural number $N$ with the following property. 
Let  $\Gamma$ be any cocompact, discrete, torsion-free subgroup of $\pizzle(\CC) $.
Suppose that $\Gamma\le\pizzle(E)$, where $E$ is a number field.
Let $v$ be a valuation of $E$. Let $x$ and $y$ be non-commuting elements of $\Gamma$. Suppose that either
\Alternatives
\item $x$ does not lie in a $\piggle(E)$-conjugate of $\pizzle(\frako_v)$, or
\item $x\in\pizzle(\frako_v)$, the characteristic of $k_v$ is greater than $N$, and $\overline h_v(x)$ has order $7$ in $\pizzle(k_v)$. 
\EndAlternatives
Then for every point $P\in\HH^3$ we have
$$\max( d_P(x), d_P(y))>0.34.$$
\EndTheorem

\Proof
We have 
$$\frac1{1+\exp{\alpha}}+\frac1{1+\exp(7\alpha)}=0.5005\ldots>\frac12.$$
Hence we may define a natural number $\delta_{\alpha,7\alpha}$ as in Corollary \ref{ranana}. I will show that the assertions of the theorem hold with $N=\max(7,\delta_{\alpha,7\alpha})$.

Suppose that $\Gamma$, $E$, $v$, $x$ and $y$ satisfy the hypotheses. Set $M=\HH^3/\Gamma$.

It suffices to prove that the conclusions are true in each of the following three cases:

{\bf Case A:} Alternative (ii) of the hypothesis holds and $\Gamma\le\pizzle(\frako_v)$. 

\smallskip

{\bf Case B:} Alternative (ii) of the hypothesis holds and $\Gamma$ is not conjugate in $\piggle(E)$ to a subgroup of $\pizzle(\frako_v)$. 

\smallskip

{\bf Case C:} Alternative (i) of the hypothesis holds.

I will consider Case A first. Since $M$ is a closed, orientable hyperbolic $3$-manifold, all abelian subgroups of $\Gamma\cong\pi_1(M)$  are cyclic; hence by Proposition \ref{hempel says}, all solvable subgroups of $\Gamma$  are cyclic. Since $x$ and $y$ do not commute, it follows that the subgroup $\Gamma_1:=\langle x,y\rangle$ of $\Gamma$ is non-solvable. Since $\Gamma$ is torsion-free, so is $\Gamma_1$.

Since we are in Case A, we have 
$\Gamma\le\pizzle(\frako_v)$; furthermore, $\overline h_v(x)$ has order $7$ in $\pizzle(k_v)$. Let $p$ denotes the characteristic of $k_v$. Since Alternative (ii) holds in this case, we have $p>N\ge7$. Thus the hypotheses of Corollary \ref{the daughter of rosie o'grady} now hold with $m=7$, and with $\Gamma_1$ playing the role of $\Gamma$. Hence if, as in the statement of Corollary \ref{the daughter of rosie o'grady}, we set 
$$\Theta_1=\langle x ^7,y \rangle\le\Gamma_1$$
and
$$\Theta_2=\langle x , y x y ^{-1}x y x ^{-1}y ^{-1}\rangle\le\Gamma_1,$$
and set $\theta_i=|\Gamma_1:\Theta_i|$ for $i=1,2$, then we have
\Equation\label{hello}\max(\theta_1,\theta_2)\ge p>N.\EndEquation

Now suppose that $\max( d_P(x), d_P(y))\le\alpha$. Then $ d_P (x ^7)\le7 d_P (x)\le7\alpha$. Since $x$ does not commute with $y$, it follows from \cite[Corollary 8.2]{hakenmarg} that $x^7$ does not commute with $y$. Hence the property of $\delta_{\alpha,7\alpha}$ stated in Corollary \ref{ranana} implies that 
\Equation\label{goodbye}
N\ge\delta_{\alpha,7\alpha}\ge|\Gamma:\langle x ^7,y \rangle|\ge |\Gamma_1:\langle x ^7,y \rangle|= \theta_1 .\EndEquation

Likewise, we have $ d_P(y x y ^{-1}x y x ^{-1}y ^{-1})\le3 d_P(x)+4 d_P(y)\le7\alpha$.
Since $y$ does not commute with $x$, it follows from two applications of \cite[Lemma 8.1]{hakenmarg} that 
$y x y ^{-1}x y x ^{-1}y ^{-1}$ does not commute with $x$. Hence the property of $\delta_{\alpha,7\alpha}$ stated in Corollary \ref{ranana} implies that 
\Equation\label{that's what you say}
N\ge\delta_{\alpha,7\alpha}\ge|\Gamma:\langle x ,
y x y ^{-1}x y x ^{-1}y ^{-1} \rangle|\ge |\Gamma_1:\langle 
x ,
y x y ^{-1}x y x ^{-1}y ^{-1} \rangle|= \theta_2.\EndEquation
Now (\ref{hello}), (\ref{goodbye}) and (\ref{that's what you say}) give a contradiction, and the proof in Case A is complete.

I will now turn to the proof in Cases B and C. In these cases $\Gamma$ is not in a $\piggle(E)$-conjugate of $\pizzle(\frako)$. Furthermore, since $\Gamma\le\pizzle(\CC)$ is discrete, torsion-free and cocompact, it has no non-trivial normal abelian subgroup.  Thus the hypotheses of Proposition \ref{it lives} hold. In Case B it follows from Assertion (3) of 
Proposition \ref{it lives} that there is an essential, faithful $\Gamma$-tree $T$ such that $x$ is $T$-elliptic, $\Per_T(x)=T$, and the $x$-period of every edge of $T$ has the form $p^r$ or $7p^r$ for some $r\ge0$. In Case C it follows from Assertion (2) of 
Proposition \ref{it lives} that there is an essential, faithful $\Gamma$-tree $T$ such that $x$ is $T$-hyperbolic; in this case, $\Per_T(x)=\emptyset$, and it is vacuously true that the $x$-period of every edge of $\Per_T(x)$  has the form $p^r$ or $7p^r$ for some $r\ge0$.  Thus the following assertion holds both in Case B and in Case C:

\Claim \label{vacuum for dust}
There exists an essential, faithful $\Gamma$-tree $T$ such that the $x$-period of every edge of $\Per_T(x)$ has the form $p^r$ or $7p^r$ for some $r\ge0$.  
\EndClaim

I will complete the proof of the theorem by showing that the conclusion follows from \ref{vacuum for dust}.

It follows from \ref{vacuum for dust} that the hypothesis of Proposition
Proposition \ref{ashes to ashes} holds with 
$\calp=\{7,p\}$. Hence one of the alternative conclusions (a), (b) or (c) of Proposition \ref{ashes to ashes} must hold. Consider the subcase in which (a) holds, i.e. at least one of the pairs $(x,yxy^{-1})$ and $(xy^{-1},y^{2})$ is independent. Then according to the case $k=2$ of \cite[Theorem 4.1]{surgery}, we have {\it either}
\Equation\label{stickin' to the onion}\frac1{1+\exp  d_P(x)}+\frac1{1+\exp  d_P(yxy^{-1})}\le\frac12\EndEquation
{\it or}
\Equation\label{till the day i fry}\frac1{1+\exp  d_P(xy^{-1})}+\frac1{1+\exp  d_P(y^2)}\le\frac12.\EndEquation
Now assume that $\max( d_P(x), d_P(y))\le\alpha$. Then we have 
$$ d_P(y x y ^{-1})\le  d_P(x)+2 d_P(y)\le3\alpha,$$ $$ d_P(x y ^{-1})\le d_P(x)+ d_P(y)\le2\alpha$$
 and $$ d_P(y ^2)\le2  d_P(y)\le2\alpha.$$ Hence
\Equation\label{the thirty-third of may}\frac1{1+\exp  d_P(x)}+\frac1{1+\exp  d_P(yxy^{-1})}\ge \frac1{1+\exp \alpha}+\frac1{1+\exp (3\alpha)}=0.68\ldots\EndEquation
{\it and}
\Equation\label{face like a ghost}
\frac1{1+\exp  d_P(xy^{-1})}+\frac1{1+\exp  d_P(y^2)}\ge \frac2{1+\exp (2\alpha)}=0.67\ldots.\EndEquation
But (\ref{the thirty-third of may}) contradicts (\ref{stickin' to the onion}), and (\ref{face like a ghost}) contradicts (\ref{till the day i fry}).
This establishes the conclusion in the subcase where (a) holds.

If (b) holds, i.e. if at least one of the pairs $(x^{-1},y)$ and $(x,y)$ is semi-independent, then it follows from \cite[Corollary 5.3]{hakenmarg} that
$$\max( d_P(x), d_P(y))=\max( d_P(x^{-1}), d_P(y))\ge\frac{\log2}2=0.346\ldots>\alpha,$$
so that the conclusion holds in this subcase as well.

Finally suppose that (c) holds, i.e. that there is an incompressible surface  $F\subset M$,  no component of which is not a fiber or semifiber, such that $\Per_{T_F}(x)$ has at least one edge, and $m_F(x)$ is divisible by either $7$ or $p$.

Since $\Per_{T_F}(x)$ has an edge, in particular $x$ is $T_F$-elliptic by \ref{periodic billiards}, and hence $yxy^{-1}$ is also $T_F$-elliptic. If 
$\Fix_{T_F}(x) \cap \Fix_{T_F}(yxy^{-1}) \ne \emptyset$, it follows from \cite[Proposition 4.11]{hakenmarg} that $x$ and $yxy^{-1}$ are independent in $\Gamma$; thus (a) holds, and by the subcase already proved, the conclusion of the theorem is true. We may therefore assume that 
$\Fix_{T_F}(x) \cap \Fix_{T_F}(yxy^{-1}) = \emptyset$. Note also that since $m_F(x)$ is divisible by either $7$ or $p$, and since $p>N\ge7$, we have $7\le m_F(x)=m_F(yxy^{-1})$. In view of the definition of $m_F$, this implies that
there is no integer $k$ with $0<k<7$ such that either $x^k$ or $(yxy^{-1})^k$ has a fixed edge in $T_F$. In particular, the hypotheses of Lemma \ref{it's da bomb} hold with $yxy^{-1}$ playing the role of $y$ in that lemma. (The condition that $x$ and $yxy^{-1}$ do not commute follows from the fact that $x$ and $y$ do not commute, in view of 
\cite[Lemma 8.1]{hakenmarg}.) Hence according to Lemma \ref{it's da bomb}, we have either $ d_P(x)>0.34=\alpha$ or $ d_P(yxy^{-1})>1.08$. If $ d_P(x)\le\alpha$ then $ d_P(yxy^{-1})\le  d_P(x)+2 d_P(y)\le3\alpha<1.08$. Hence we must have $ d_P(x)>\alpha$, and the conclusion of the theorem is established in all cases. \EndProof

\Corollary\label{but who needs it}
Let  $\Gamma$ be any cocompact, discrete, torsion-free subgroup of $\pizzle(\CC) $. Let $x$ and $y$ be non-commuting elements of $\Gamma$. Suppose that $x=[A]$, where $A$ is an element of $\zzle(\CC)$ such that $\trace A$ is not an algebraic integer. 
Then for every point $P\in\HH^3$ we have
$$\max( d_P(x), d_P(y))>0.34.$$
\EndCorollary

\Proof
It follows from \ref{just conjugate} that $\rho$ is conjugate in $\ggle(\CC)$ to a representation of $\Gamma$ in $\zzle(E)$ for some number field $E$. Hence we may assume without loss of generality that $\rho(\Gamma)\subset\zzle(E)$. Set $\tau=\trace A$. Since $\tau$ is not an algebraic integer, it follows from statement ($\alpha$) on p. 264 of \cite{rib} that there is a valuation $v$ of $E$ such that $v(\tau)<0$. 
Thus $\tau\notin\frako_v$, and hence $A$ does not lie in a $\ggle(E)$-conjugate of $\zzle(\frako_v)$. Hence
$x$ does not lie in a $\piggle(E)$-conjugate of $\pizzle(\frako_v)$ The assertion now follows from Theorem \ref{N lives}.
\EndProof

\section{Character varieties}\label{character variety section}

\Number\label{erence}
This section is devoted to a little background on the variety of $\zzle(\CC)$-characters of a finitely generated group $\Gamma$, which will be needed for Section \ref{curve and degree section}. I will be taking a point of view close to that of \cite[Section 1]{varreps} and \cite[Section 4]{handbook}, and I will briefly review the relevant material here. 

The set of all representations of $\Gamma$ in $\zzle(\CC)$ will be denoted $R(\Gamma)$. An  {\it$\zzle(\CC)$-character}, or more briefly a {\it character}, of $\Gamma$ is a complex-valued function on $\Gamma$ of the form $\chi=\chi_\rho:\gamma\mapsto\trace\rho(\gamma)$ for some representation $\rho\in R(\Gamma)$. The set of all $\zzle(\CC)$-characters of $\Gamma$ will be denoted $X(\Gamma)$. The map $\rho\mapsto\chi_\rho$ from $R(\Gamma)$ to $X(\Gamma)$ will be denoted by $t$.

If $(\xi_1,\ldots,\xi_n)$ is a finite system of generators for $\Gamma$, there is a bijective correspondence
$\rho\leftrightarrow(\rho(\xi_1),\cdots,\rho(\xi_n))$ between
$R(\Gamma)$ and a (closed) algebraic subset of the complex affine space
$\calm_2(\CC)^n=\CC^{4n}$. 
As in \cite{handbook} I will identify $R(\Gamma)$ with this algebraic set via this correspondence, once a generating system has been specified.
\EndNumber

\Number\label{more on characters} 
For each $\gamma\in\Gamma$, one can define a function
$\tau_\gamma:R(\Gamma)\to\CC$ by
setting $\tau_\gamma(\rho)=\chi_{\rho}(\gamma)=\trace\rho(\gamma)$ for every
representation $\rho\in R(\Gamma)$. Then $\tau_\gamma$ is defined by polynomial functions in the ambient coordinates of $\CC^{4n}$, i.e. it belongs to the coordinate ring $\CC(R(\Gamma)$. Let $T(\Gamma)$ denote the sub-ring of $\CC(R(\Gamma)$ generated by all functions $\tau_\gamma$ for $\gamma\in\Gamma$. According to \cite[Proposition 4.4.2]{handbook}, 
 If we set $N=2^n-1$, and we index the words of the
form $\xi_{i_1}\dots \xi_{i_k}$, with $1\le k\le n$ and
$1\le i_1<\dots<i_k\le n$, in some order as
$V_1,\dots,V_N$, then
$\tau_{V_1},\dots,\tau_{V_N}$ generate $T(\Gamma)$. This implies in particular that a character of $\Gamma)$ is determined by its values at $V_1,\ldots,V_N$. Hence the map $\chi\mapsto(\chi(V_1),\ldots,\chi(V_N))$ is a bijection of $X(\Gamma)$ to some subset of $\CC^N$, which I will henceforth identify with $X(\Gamma)$.

In terms of the identifications described above, $X(\Gamma)$ is the image of the map $t:R(\Gamma)\to\CC^N$ defined by
$t(\rho)=(\tau_{V_1}(\rho),\dots,\tau_{V_N}(\rho))$.
As is stated in \cite{handbook} and proved as \cite[Proposition 1.4.4]{varreps}, the set $X(\Gamma):=t(R(\Gamma)$ is an algebraic subset of $\CC^N$. For this reason it is called the ($\zzle(\CC)$)-{\it character variety} of $\Gamma$. It is immediate that conjugate representations in $R(\Gamma)$ have the same image under $\rho$. In the converse direction, one crucial property of $X(\Gamma)$ which will be used below is that if $\rho,\rho'\in R(\Gamma)$ satisfy $t(\rho)=t(\rho')$, and if $\rho$ is irreducible, then $\rho$ and $\rho'$ are conjugate representations; this is Proposition 1.5.2 of \cite{varreps}.

Let $\gamma$ be any element of $\Gamma$. Since $\tau_{V_1},\dots,\tau_{V_N}$ generate $T(\Gamma)$, there is an $N$-variable integer polynomial $f$ such that $\tau_\gamma=f (\tau_{V_1},\dots,\tau_{V_N})$. The polynomial function $f$ on $\CC^N$ restricts to a function $I:X(\Gamma)\to\CC$ with the property that $I\circ t=\tau_\gamma$; furthermore, this property characterizes $I$ since $t$ is surjective. In particular $I$ is uniquely determined by the element $\gamma\in\Gamma$, and will be denoted $I_\gamma$. Since $I_\gamma$ is defined by an integer polynomial in the ambient coordinates, it is in particular an element of the coordinate ring $\CC[X(\Gamma)]$.

The definitions given above involve a specific choice of a generating system for $\Gamma$. It is not hard to show that, up to isomorphism of algebraic sets, $R(\Gamma)$ and $X(\Gamma)$ are independent of the choice of generators. However, in this paper I will always be working in terms of a particular system of generators, so that $R(\Gamma)$ and $X(\Gamma)$ will be concretely defined as subsets of affine spaces.
\EndNumber

\Number\label{weinermandias}
Note that with the definitions of $R(\Gamma)$ and $X(\Gamma)$ given above, if $F$ denotes the free group on the generators $\xi_1,\cdots,\xi_n$ and $N\triangleleft F$ is the group of defining relations for $\Gamma$, then $R(\Gamma)$ is a(n algebraic) subset of $R(F)$, and hence $X(\Gamma)=tR(\Gamma)$ is a(n algebraic) subset of $X(F)=t(R(F))$.
\EndNumber

\Number\label{you and who else}
Now suppose that $H$ is a finitely generated subgroup of a finitely group $\Gamma$. If $\chi$ is an $\zzle(\CC)$-character of $\Gamma$ then $\chi|H$ is an $\zzle(\CC)$-character of $H$; indeed, if $\chi=\chi_\rho$ for some representation $\rho:\Gamma\to\zzle(\CC)$, then $\chi|H=\chi_{\rho|H}$. Hence we may define a {\it restriction map} $r:\chi\mapsto\chi|H$ from 
$r:X(\Gamma)\to X(H)$. 

%If we regard $r:X(\Gamma)$ and $ X(H)$ as algebraic sets, 

Now suppose that
%by $is given by polys in the coordinates, for application in proof of Proposition \ref{long reid type prop}. Also make sure the non-constancy assertion there is obvious. Use this stuff as needed:
$(\xi_1,\ldots,\xi_n)$ and $(h_1,\ldots,h_m)$ are generating systems for $\Gamma$ and $H$ respectively Let us set $N=2^n-1$ and $M=2^m-1$, and identify $X(\Gamma)$ and $X(H)$ with algebraic subsets of $\CC^N$ and $\CC^M$ respectively, as in \ref{more on characters}. Then 
%for suitable words $V_1,\ldots,V_N$ in $(\xi_1,\ldots,\xi_n)$, each character $\chi\in X(\Gamma)$ is identifed with $(I_{V_1}( \chi),\ldots,I_{V_N}( \chi))$; while 
for suitable words $W_1,\ldots,W_M$ in $(h_1,\ldots,h_m)$, each character $\psi\in X(H)$ is identifed with  $(I_{W_1}( \psi),\ldots,I_{W_M}( \psi))$. Since $H\le\Gamma$, we may regard the $W_i$ as elements of $\Gamma$, and for any $\chi\in X(\Gamma)$ we have $r(\chi)=(I_{W_1}( \chi),\ldots,I_{W_M}( \chi))$. According to \ref{more on characters} the functions $I_{W_i}$ are given by integer polynomials in the coordinates of $\CC^N$; hence the map $r$ is defined by integer polynomials in the ambient coordinates. 
\EndNumber

The following result generalizes Proposition 1.1.1 of \cite{varreps}.

\Proposition\label{when i'm 65}
Let $\Gamma$ be a finitely generated group, and let $V$ be an affine algebraic subset of $R(\Gamma)$ which is invariant under conjugation (that is, if a representation $\rho$ belongs to $V$, so does the representation $\rho^A:\gamma\mapsto A\rho(\gamma)A^{-1}$ for every $A\in\zzle(\CC)$). Then each irreducible component of $V$ is also invariant under conjugation.
\EndProposition

\Proof
Let $Z$ be any component of $V$. Consider the map of algebraic sets $F:Z\times\zzle(\CC)\to R(\Gamma)$ defined by $F(\rho,A)=\rho^A$. Since $V$ is  invariant under conjugation, we have $F(Z\times\zzle(\CC))\subset V$. Since the product $Z\times\zzle(\CC)$ of irreducible varieties is irreducible, and since $F$ is defined by polynomials in the ambient coordinates, $F(Z\times\zzle(\CC))$ must be contained in a single component of $V$. Since $F(Z\times\zzle(\CC))\supset F(Z\times\{I\})=Z$, we must have $F(Z\times\zzle(\CC))=Z$, so that $Z$ is invariant under conjugation as claimed.
\EndProof

\Proposition\label{spilious becs} Let $\Gamma$ be a finitely generated group, and let $V$ be an irreducible affine algebraic subset of $R(\Gamma)$ which is invariant under conjugation and contains of an irreducible representation. Then $t(V)\subset X(\Gamma)$ is a (closed) affine algebraic set.
\EndProposition

\Proof
The special case of this result in which $V$ is an irreducible component of $R(\Gamma)$ containing an irreducible representation is proved as Proposition 1.4.1 of \cite{varreps}. An examination of the proof of \cite[Proposition 1.4.1]{varreps} reveals that the assumption that $V$ is an irreducible component of $R(\Gamma)$ is used only to guarantee that it is invariant under conjugation.
\EndProof

\Proposition\label{bilious specs} Let $\Gamma$ be a finitely generated group, and let $C$ be an irreducible affine algebraic subset of $X(\Gamma)$ which contains the character of an irreducible representation. Then there is an irreducible algebraic subset $Z$ of $R(\Gamma)$ such that $t(Z)=C$.
\EndProposition

\Proof
The set $t^{-1}(C)$ is an affine algebraic set and therefore has finitely many components, say $Z_1,\ldots,Z_n$. Since conjugate representations have the same character, $t^{-1}(C)$ is invariant under conjugation. Hence by Proposition \ref{when i'm 65}, each $Z_i$ is invariant under conjugation.

For each $i$, let $W_i\subset C$ denote the Zariski closure of $t(Z_i)$. Since $t:R(\Gamma)\to X(\Gamma)$ is surjective, we have $C=W_1\cup\cdots\cup W_n$. Since $C$ is irreducible we have $C=W_i$ for some $i$, and after re-indexing we may assume that $C=W_1=\overline{t(Z_1)}$. 

Let $Y$ denote the subset of $R(\Gamma)$ consisting of all reducible representations. It follows from \cite[Corollary 1.2.2]{varreps} that $Y=t^{-1}(S)$ for some Zariski-closed set $S\subset X(\Gamma)$. By hypothesis we have $C\not\subset S$, and hence $C\setminus S$ is a non-empty subset of $C$, Zariski-open in $C$. It therefore meets the Zariski-dense subset $t(Z_1)$ of $C$. This means that $Z_1$ contains at least one irreducible representation. As $Z_1$ is invariant under conjugation, it now follows from Proposition \ref{spilious becs} that $t(Z_1)\subset X(\Gamma)$ is Zariski-closed, and is therefore equal to $C$.
\EndProof

\Number\label{truncation}
I will be using some of the conventions of \cite{finitistic} in discussing finite-volume hyperbolic $3$-manifolds. Suppose that $M$ is an orientable hyperbolic $3$-manifold of  finite volume, and let $N$ denote a truncation of $M$ in the sense of \cite[Section 3]{finitistic}. According to \cite[Proposition 3.7]{finitistic}, $N$ is compact, each component of $\partial N$ is a torus, and each component of $\overline{M-N}$ is diffeomorphic to $T^2\times[0,\infty)$. Hence $N$ is a strong deformation retract of $M$. 
\EndNumber

I will need the following well-known fact:

\Lemma\label{friend matrix of discord}Let $M$ be an orientable hyperbolic $3$-manifold of  finite volume, let $N$ be a truncation of $M$, and let $\gamma\ne1$ be an element of $\pi_1(M)$ that is mapped to a parabolic element under some discrete faithful representation of $\pi_1(M)$ in $\pizzle(\CC)$. Then $\gamma$ lies in a conjugate of  the image of the inclusion homomorphism $\pi_1(T)\to\pi_1(M)$ for some component $T$ of $\partial N$.
\EndLemma

\Proof
Let us write $M=\HH^3/\Gamma$, where $\Gamma\le\pizzle(\CC)$ is discrete and torsion-free. Up to conjugacy there is a canonical isomorphism $J:\pi_1(M)\to\Gamma$. If $\rho:\pi_1(M)\to\pizzle(\CC)$ is the given discrete faithful representation, then Mostow rigidity implies that $\rho\circ J^{-1}$ extends to an automorphism of $\zzle(\CC)$; hence we may assume after a conjugation that $\rho$, up to complex conjugation, is the composition of $J$ with the inclusion $\Gamma\to\pizzle(\CC)$. Thus the hypothesis implies that $J(\gamma)$ is parabolic. The centralizer $C$ of $J(\gamma)$ then consists entirely of parabolic elements of $\Gamma$, and is free abelian of rank $1$ or $2$. If $C$ has rank $1$, then applying the Margulis lemma as on p. 64 of \cite{morgan} we deduce that $M$ has a ``$\ZZ$ cusp'' and therefore has infinite volume, a contradiction. Hence $C\cong\ZZ\times\ZZ$. It then follows from \cite[Proposition 3.5]{finitistic} that $C$ is conjugate of  the image of the inclusion homomorphism $\pi_1(T)\to\pi_1(M)$ for some component $T$ of $\partial N$. Since $\gamma\in C$, the conclusion follows.
\EndProof

The following result is similar to \cite[Proposition 1.1.1]{sepsurfs}, but I am supplying a proof here because the argument given in \cite{sepsurfs} contained a reference to L. Lok's unpublished thesis.

\Proposition\label{loop de loo}
Let $M$ be an orientable hyperbolic $3$-manifold of finite volume, let $N$ denote a truncation of $M$, and let $T_1,\ldots,T_k$ denote the components of $\partial N$. For $i=1,\ldots,k$ let $\gamma_i$ be an element of $\pi_1(M)$ representing the conjugacy class determined by some homotopically non-trivial closed curve in $T_i$. Suppose that $\chi_0\in X(\pi_1(M))$ is the character of a discrete, faithful representation of $\pi_1(M)$ in $\zzle(\CC)$.  Then $x_i:=I_{\gamma_i}(\chi_0)=\pm2$ for $i=1,\ldots,k$, and $\chi_0$ is an isolated point of the set $\bigcap_{i=1}^kI_{\gamma_i}^{-1}(\{x_i\})$.
\EndProposition

\Proof
Let $\rho_0:\pi_1(M)\to\zzle(\CC)$ be a discrete, faithful representation 
with character $\chi_0$. 
For $i=1,\ldots,k$ let $H_i\le\pi_1(M)$ denote the subgroup defined, up to conjugacy, as the image of the inclusion homomorphism $\pi_1(T_i)\to\pi_1(M)$. It follows from \cite[Proposition 3.5]{finitistic} that each $T_i$ is a free abelian group of rank $2$, and hence $\rho_0(T_i)$
is a discrete, rank-$2$ free abelian subgroup of $\zzle(\CC)$. Hence $\rho_0(T_i)$ is conjugate to a group of matrices of the form $\pm\begin{pmatrix}1&\lambda\cr0&1\end{pmatrix}$; in particular, since $\gamma_i$ lies in a conjugate of $H_i$, we have
$x_i=\trace\rho_0(\gamma_i)=\pm2$. This is the first assertion of the proposition.

Assume that the second conclusion is false, and let $C\subset X(\pi_1(M))$ be an affine algebraic curve with $\chi_0\subset C\subset\bigcap_{i=1}^kI_{\gamma_i}^{-1}(\{x_i\})$. Note that since $\rho_0$ is faithful and $M$ has finite volume, $\rho_0(\pi_1(M))$ is non-solvable and hence $\rho_0$ is irreducible. It therefore follows from Proposition \ref{bilious specs} that there is an irreducible algebraic subset $Z$ of $R(\pi_1(M))$ such that $t(Z)=C$. Then for every $\rho\in Z$, and for $i=1,\ldots,k$, we have $\trace\rho(\gamma_i)=x_i=\pm2$.

On the other hand, since $\rho_0$ is faithful, there is a Zariski-dense subset $U$ of $Z$ such that $\trace\rho(\gamma_i)=x_i=\pm2$ for every $\rho\in U$ and $i=1,\ldots,k$. Hence for each $i$ there exists $A_i\in\ggle(\CC)$ such that $A_i\rho_i(\gamma_i)A_i^{-1}=\pm\begin{pmatrix}1&1\cr0&1\end{pmatrix}$. We may take the subgroup $H_i$ to be chosen within its conjugacy class so that $\gamma_i\in H_i$; and since $H_i$ is abelian, for any $\gamma\in H_i$ we then have $A_i\rho_i(\gamma_i)A_i^{-1}=\pm\begin{pmatrix}1&c_\gamma\cr0&1\end{pmatrix}$ for some $c_\gamma\in\CC$. In particular:

\Claim\label{duck}Given any $\gamma\in H_i$, we have $\trace\rho(\gamma)=\pm2$ for every $\rho\in U$ and hence for every $\rho\in Z$.
\EndClaim

According to Lemma \ref{friend matrix of discord}, the only parabolic elements of $\Pi_\CC\circ\rho_0(\pi_1(M))$ are the conjugates of $\Pi_\CC\circ\rho_0(H_i)$ for $i=1,\ldots,k$. It therefore follows from \ref{duck} that for every $\gamma\in\pi_1(M)$ such that $\Pi_\CC\circ\rho_0(\gamma)$ is parabolic, $\Pi_\CC\circ\rho(\gamma)$ is parabolic
for every $\rho\in Z$. In view of \cite[Definition 9.1]{marden}, this implies that for every $\epsilon>0$ there is a (classical) open neighborhood $W$ of $\rho_0$ in $Z$ such that $\Pi_\CC\circ\rho$ is an $\epsilon$-deformation of the Kleinian group $\Pi_\CC\circ\rho_0$ for every $\rho\in W$. But the proof of \cite[Lemma 9.2]{marden} shows that $\epsilon>0$ may be chosen in such a way that every $\epsilon$-deformation of $\Pi_\CC\circ\rho_0$ is a discrete and faithful representation. It now follows from Mostow rigidity that for every representation in $\rho\in W$, the representation $\Pi_\CC\circ\rho:\pi_1(M)\to\pizzle(\CC)$ is conjugate to $\Pi_\CC\circ\rho_0$. But since $\pi_1(M)$ is finitely generated, the representation $\Pi_\CC\circ\rho_0$ admits only finitely many lifts to $\zzle(\CC)$. Hence the set $W\subset R(\pi_1(M))$ meets only finitely many conjugacy classes of representations, and thus $t(Z)\subset C$ is finite. This is impossible, because the map $t$ from the irreducible algebraic set $Z$ to $C$, defined by polynomials in the ambient coordinates, is surjective, and must therefore map the non-empty open set $W$ onto an infinite set. Thus the second conclusion of the proposition is established.

(The definition of a Kleinian group used in  \cite{marden} includes the condition that the group have a non-empty set of discontinuity on the sphere at infinity. Therefore the {\it statement} of  \cite[Lemma 9.2]{marden} does not directly apply to the Kleinian group $\Pi_\CC\circ\rho_0(\pi_1(M))$ in the argument above, as its set of discontinuity is empty. However, the condition of having empty set of discontinuity does not appear to be used in the proof of  \cite[Lemma 9.2]{marden}.)
\EndProof

The following result will be quoted in Section \ref{masterpiece section}:

\Proposition\label{i never sausage a bad pun} Let $F$ be a free group on two generators $\xi$ and $\eta$. Then the character variety $X(F)$, defined in terms of the generating system $(\xi,\eta)$, is $\CC^3$, and the map $t:R(F)\to R(\eta)$ is given by $\rho\mapsto(\trace\rho(\xi),\trace\rho(\eta),\trace\rho(\xi\eta))$.
\EndProposition

\Proof
In the notation of \ref{erence} we have $N=2^2-1=3$, and the $V_i$ may be indexed so that $V_1=\xi$, $V_2=\eta$ and $V_3=\xi\eta$. Hence $X(F)$ is a (closed) algebraic subset of $\CC^3$, and we have $t(\rho)=(\trace\rho(\xi),\trace\rho(\eta),\trace\rho(\xi\eta))$ for each $\rho\in R(F)$. It remains to show that $X(F)$ is all of $\CC^3$.

Note that $R(F)=\zzle(\CC)\times\zzle(\CC)\subset M_2\times M_2$, so that $X$ is irreducible and $\dim R(F)=6$. Since $F$ admits a faithful representation in $\zzle(\CC)$, in particular $R(F)$ contains an irreducible representation. Hence by \cite[Corollary 1.5.3]{varreps}, we have $\dim X(\Gamma)=\dim R(\Gamma)-3=3$. Thus $X(\Gamma)$ is a $3$-dimensional closed algebraic subset of $\CC^3$ and is therefore all of $\CC^3$.
\EndProof

\section{A little algebraic geometry}

The algebro-geometric observations made in this section will be needed in Section \ref{curve and degree section}. 

Recall that a complex affine algebraic subset of an affine space $F^N$, where $F$ is an algebraically closed field, is said to be {\it defined over} a subfield $K$ of $F$ if it is the locus of zeros of a set of polynomials whose coefficients lie in $K$.

\Proposition\label{cubar}Let $K$ be a subfield of an algebraically closed field $F$, let $N$ be a positive integer, and let $S$ be a subset of $K^N\subset F^N$. Then the Zariski closure of $S$ in $F^N$ is defined over $K$.
\EndProposition

\Proof
Let $I\subset F[X_1,\ldots,X_N]$ denote the ideal consisting of all polynomials that vanish on $S$. The Zariski closure of $S$ in $F^N$. Then $Y$ is the locus of zeros of $I$ in $F^N$; hence we need only show that $I$ is generated by polynomials with coefficients in $K$. 

By the Hilbert basis theorem, $I$ is finitely generated. Hence we may fix an integer $k>0$ such that $I$ is generated by polynomials of degree at most $k$. Let $M$ denote the set of all monomials in $F[X_1,\ldots,X_N]$ having degree at most $k$, and let $V$ denote the linear span of $M$ in
$F[X_1,\ldots,X_N]$. Then $I$ is generated by $I\cap V$.

For each $s\in S$, let $I_s\subset F[X_1,\ldots,X_N]$ denote the ideal consisting of all polynomials that vanish at $s$. Then $I=\bigcap_{s\in S}I_s$ and hence $I\cap V=\bigcap_{s\in S}(I_s\cap V)$. Since $V$ is a finite-dimensional vector space and the
$I_s\cap V$ are subspaces, there is a finite set $T\subset S$ such that
$I\cap V=\bigcap_{s\in T}(I_s\cap V)$. Hence if we define a linear map $A:V\to F^T$ of finite-dimensional vector spaces by $A(f)=(f(s))_{s\in T}$, we have $I\cap V=\ker A$.
But since $S\subset K^N$, the matrix of $A$ with respect to the basis $M$ of $V$ and the standard basis for $F^T$ has entries in $K$. It follows that $\ker A$ has a basis $B$ consisting of vectors whose coefficients in the basis $M$ have coefficients in $K$. This means that $B\subset F[X_1,\ldots,X_N]$. As $I\cap V$ generates the ideal $I$, the basis $B$ also generates $I$, and thus $I$ is indeed generated by polynomials in $F[X_1,\ldots,X_N]$.
\EndProof

\Proposition\label{slider}Let $K$ be an algebraically closed subfield of an algebraically closed field $F$, let $V$ and $Y$ be (closed) algebraic subsets of affine spaces over $F$. Let $f:V\to Y$ be a map defined by polynomials in the ambient coordinates with coefficients in $K$.  Suppose that $V$ is defined over $K$ and that $f(V)$ is Zariski-dense in $Y$. Then $Y$ is defined over $K$.
\EndProposition

\Proof
Let $F^P$ and $F^N$ be the affine spaces containing $V$ and $Y$.
Since $K$ is algebraically closed and $V$ is defined over $K$, it follows from \cite[Theorem 30.2]{isaacs} that $V\cap K^P$ is Zariski-dense in $V$. Since $f$ is continuous in the Zariski topology, $S:=f(V\cap K^P)$ is Zariski-dense in $f(V)$, and hence in $Y$. But since $f$ is defined by polynomials with coefficients in $K$, we have 
$S\subset f(K^P)\subset K^N$. Hence by Proposition \ref{cubar}, the Zariski closure $Y$ of $S$ in $F^N$ is defined over $K$.
\EndProof

An algebraic set in $\CC^2$ will be termed {\it purely one-dimensional} if all its irreducible components are curves.

\Lemma\label{harald}
Let $V$ be an algebraic curve in $\CC^2$ which is defined over the field $\overline\QQ$ of algebraic numbers in $\CC$. Then there is a purely one-dimensional algebraic set $Z$ in $\CC^2$ such that (1) $Z\supset V$, (2) $Z$ is defined over $\QQ$, and (3) $Z$ has no proper, non-empty purely one-dimensional algebraic subset which is defined over $\QQ$.
\EndLemma

%\Lemma\label{harald}
%Let $C$ be a curve in $\CC^2$ which is defined over the field $\overline\QQ$ of algebraic numbers in $\CC$. Then there is a non-zero two-variable integer polynomial $A\in\ZZ[X,Y]$ which vanishes everywhere on $C$.
%\EndLemma

\Proof
Since the curve $V$ is defined over $\overline\QQ$, it is the locus of zeros of a non-zero polynomial $f\in \overline\QQ[X,Y]$. Let $K$ be a finite normal extension of $\QQ$ containing the coefficients of $f$. Write $f=\sum a_{ij}X^iY^j$, let $G$ denote the Galois group of $K$ over $\QQ$, and for each $\sigma\in G$ set  $f^\sigma=\sum \sigma(a_{ij})X^iY^j$. Set $B_0=\prod_{\sigma\in G}f^\sigma$. Then $B_0$ is non-zero and vanishes on $C$; furthermore, its coefficients are fixed by every element of $G$, and hence $B_0\in\QQ[X,Y]$. Hence if we define  $Z_0$ to be the locus of zeros of $B_0$, then (1) and (2) hold with $Z_0$ in place of $Z$. Now among all purely one-dimensional algebraic sets satisfying (1) and (2), let $Z$ be one that has the smallest possible number of irreducible components. Since (1) and (2) hold, a defining polynomial $B$ for $Z$ lies in $\QQ[X,Y]$ and is divisible by $f$. If (3) does not hold then a defining polynomial $B$ for $Z$ may be written as as product of two non-constant polynomials $B_1$ and $B_2$ where $B_1\in\QQ[X,Y]$. Since $B\in\QQ[X,Y]$ it then follows that $B_2\in\QQ[X,Y]$ as well. Thus if $Z_i$ denotes the locus of zeros of $B_i$ then $Z_1$ and $Z_2$ are defined over $\QQ$. Since $f|B$, we have either $f|B_1$ or $f|B_2$, and hence one of the $Z_i$ contains $V$, in contradiction to the minimality of $Z$.
\EndProof

\section{Character curves and degrees of number fields}\label{curve and degree section}

The first main result of this section, Proposition \ref{long reid type prop} below, is a partial generalization of a result due to Long and Reid \cite[Theorem 3.2]{long-reid}, and the proof closely parallels the proof of their result. Long and Reid describe their result as ``a strong form of an observation due to Hodgson.''  I will define the {\it length} of an integer polynomial $f$ (in an arbitrary number of variables) to be the sum of the absolute values of the coefficients of $f$. Note that $\length f=0$ if and only if $f=0$.

Long and Reid's proof of their Theorem 3.2 implicitly involves the following fact.

\begin{lemma}\label{long reid type lemma}Let $\cald$ and $L$ be non-negative integers. Let $W\subset\CC$ denote the set of all algebraic numbers $w$ such that (i) $w$ has degree at most $\cald$ and (ii) $w$ is a root of some (possibly reducible) integer polynomial $f$ with $0<\length f\le L$. Then $W$ is finite.
\EndLemma

\Proof
Recall from \cite[Section 3]{long-reid} that the {\it Mahler measure} $\meas(f)$ of a non-zero one-variable integer polynomial 
$f(X)=a_nX^n+\cdots+a_0=a_n(X-r_1)\cdots(X-r_n)$  is defined by
$$\meas(f)=|a_n|\prod_{i=1}^n\max(|r_i|,1).$$
According to \cite[Lemma 3.3]{long-reid}, we have $\meas(f)\le\length(f)$ for any $f$. 

Every $w\in W$ is the root of a non-zero integer polynomial $f$ with $\meas(f)\le\length(f)\le L$. In particular every root of $f$ has absolute value at most $L$. It follows from Gauss's lemma that we may write $f=g_1\cdots g_r$, where $g_1,\ldots,g_m$ are $\ZZ$-polynomials which are $\QQ$-irreducible. We may take the $g_i$ to be indexed so that $w$ is a root of $g_1$, and we may write $g_1(X)=b(X-s_1)\cdots(X-s_d)$, where $d=\deg g=\deg w\le\cald$, and the $s_j$ are among the roots of $f$ and hence have absolute value at most $L$. Furthermore, the leading coefficient $b$ of $g_1$ divides $a_n$, and hence $|b|\le L$. Since the coefficients of $g_1/b$ are elementary symmetric functions of $s_1,\ldots,s_d$, they are all of absolute value at most $2^dL^d$. Hence the coefficients of $g_1$ are all of absolute value at most $2^DL^{D+1}$. As there are only finitely many integer polynomials of degree at most $D$ whose coefficients are all of absolute value at most $2^DL^{D+1}$, the conclusion follows. 
\EndProof

\begin{proposition}\label{long reid type prop}
 Let $M$ be an orientable hyperbolic 3-manifold of finite volume.  Let $C$ be an algebraic curve contained in $X(\pi_1(M))$ which is defined over the algebraic closure $\overline{\QQ}\subset\CC$ of $\QQ$ and contains the character of a discrete, faithful representation. 
Let $D$ be a positive integer, and let $\cals$ denote the set of points of $C$ defined by 
representations $\rho$ such that (a) $\rho(\pi_1(M))$ is discrete and torsion-free, and (b) $\trace(\rho(\pi_1(M))\subset L$ for some number field $L$ of degree at most $D$.  Then $\cals$ is a closed and discrete subset of $C$ in the complex topology; that is, the intersection of $\cals$ with any compact subset of $C$ is finite. 
\end{proposition}

\Proof
Let $N$ be a truncation of $M$ (cf. \ref{truncation}).
Let $T_1,\ldots,T_k$ denote the components of $\partial N$. I claim: 

\Claim\label{eat your peas}There is an index $i_0\in\{1,\ldots,k\}$ such that for every
element $\gamma$ of $\pi_1(M)$ which represents the conjugacy class determined by some homotopically non-trivial closed curve in $T_{i_0}$,
the restriction of the function $I_\gamma$ to $C$ (see \ref{erence}) is non-constant. 
\EndClaim

%(Note that since the function $I_\gamma$ depends only on the conjugacy class of $\gamma$, the truth of \ref{eat your peas} is independent of the choices of the subgroups $H_i$ within their conjugacy classes.)

To prove \ref{eat your peas}, assume that for every $i\in\{1,\ldots,k\}$ there is a non-trivial element $\gamma_i\in H_{i}$ such that $I_{\gamma_i}|C$  is constant. 
By hypothesis, $C$ contains the character $\chi_0$ of some discrete, faithful representation of $\pi_1(M)$ in $\zzle(\CC)$. According to Proposition \ref{loop de loo} we have $x_i:=I_{\gamma_i}(\chi_0)=\pm2$ for $i=1,\ldots,k$, and $\chi_0$ is an isolated point of the set $\bigcap_{i=1}^kI_{\gamma_i}^{-1}(\{x_i\})$. But since
$I_{\gamma_i}|C$  is constant 
for $i\in\{1,\ldots,k\}$, we have $C\subset\bigcap_{i=1}^kI_{\gamma_i}^{-1}(\{x_i\})$. This is a contradiction, and thus \ref{eat your peas} is proved.

After re-indexing the $T_i$ if necessary, we may assume that the index $i_0$ given by \ref{eat your peas} is equal to $1$. Now fix a basepoint $\star\in T_1$ and let $H_1$ denote the image of the inclusion homomorphism from $\pi_1(T_1,\star)$ to $\pi_1(M,\star)=\pi_1(M)$. Then \ref{eat your peas} implies that

\Claim\label{off with her band-aid}
$I_\gamma|C$  is non-constant for every non-trivial element $\gamma\in H_1$.
\EndClaim

As in \ref{you and who else}, we may define a restriction map $r:X(\Gamma)\to X(H_1)$.
It follows in particular from \ref{off with her band-aid} that $r$ is non-constant. Since the discussion in \ref{you and who else} shows that $r$ is defined by polynomials in the ambient coordinates, we deduce that

\Claim\label{corson porson}
$Y:=\overline{r(C)}\subset X(H_1)$ is an (affine algebraic) curve, and the map $r|C:C\to Y$ is finite-to-one.
\EndClaim

The discussion in \ref{you and who else} also shows that the polynomials defining $r$ have coefficients in $\ZZ\subset\overline\QQ$. Furthermore, $r(C)$ is dense in the complex topology of $Y$ and therefore in its Zariski topology. We may therefore apply Proposition \ref{slider}, letting $r$, $\CC$ and $\overline\QQ$ play the respective roles of
$V$, $f$, $F$ and $K$, to deduce that
\Claim\label{porlock sherlock}
$Y$ is defined over $\overline\QQ$.
\EndClaim

Let $\Delta$ denote the group homomorphism from $\CC^\times=\CC-\{0\}$ to $\zzle(\CC)$ defined by $\Delta(x)=\begin{pmatrix}x&0\cr 0& x^{-1}\end{pmatrix}$.
Define a map $p:(\CC^\times)\times (\CC^\times)\to X(H_1)$ by $p(x,y)=t(\phi)$ where $\phi(\lambda)=\Delta(x)$ and $\phi(\mu)=\Delta(y)$. 

It follows from Proposition \ref{i never sausage a bad pun}, and from the discussion in \ref{weinermandias}, that the functions $I_\lambda$, $I_\mu$ and $I_{\lambda\mu}$ may be taken as coordinate functions on $X(H_1)$. In terms of these coordinates we have
\Equation\label{no ducks were injured} p(x,y)=(x+x^{-1},y+y^{-1},xy+x^{-1}y^{-1})
\EndEquation
 for any $(x,y)\in (\CC^\times)\times (\CC^\times)$.

It is clear from \ref{no ducks were injured} that $p$ is finite-to-one. 
Since $H_1$ is abelian, each representation $\rho\in R(H_1)$ is conjugate to an upper triangular representation, and hence has the same character as a diagonal representation. It follows that $p$ is surjective.

Set $\calv$ denote the closure of $p^{-1}(Y)$, in $\CC^2$, so that $\calv\cap(\CC^\times)^2=p^{-1}(Y)$. Since $Y$ is an algebraic curve defined over $\QQ$, and since $p:(\CC^\times)\times (\CC^\times)\to X(H_1)$ is a finite-to-one surjection defined by rational functions with rational coefficients in the ambient coordinates (cf. \ref{no ducks were injured}),
$\calv\subset\CC^2$ is a purely one-dimensional algebraic set defined over $\overline\QQ$. Choose an irreducible component $V$ of $\calv$.
Since $\calv$ is defined over the algebraically closed field $\overline\QQ$, it follows from \cite[Theorem 3.10.9]{nagata} that its irreducible component $V$ is also defined over $\overline\QQ$.
Thus the curve $V$ satisfies the hypotheses of Lemma \ref{harald}, and there therefore exists an algebraic set $Z\subset\CC^2$ for which Conditions (1)--(3) of Lemma \ref{harald} hold. It follows from Condition (2) of Lemma \ref{harald} that $Z$ is the locus of zeros of a non-zero two-variable polynomial $A\in\ZZ[X,Y]$.  According to Condition (1) of Lemma \ref{harald} we have $V\subset Z$.

For any two integers $m_1$ and $m_2$ that are not both zero, set $g_{m_1,m_2}(X)=A(X^{m_1},X^{m_2})$, so that $g_{m_1,m_2}$ is a Laurent polynomial in $X$ with integer coefficients. I claim:

\Claim\label{that ain't hay}
For any two integers $m_1$ and $m_2$ that are not both zero, we have
$g_{m_1,m_2}\ne0$.
\EndClaim

To prove \ref{that ain't hay}, assume that $g_{m_1,m_2}=0$, so that $A(x^{m_1},x^{m_2})=0$ for every $x\in\CC$. Then $(x^{m_1},x^{m_2})\in Z$ for every $x\in\CC$, so that 
$$Z\supset Z_0:=\{(x^{m_1},x^{m_2}):x\in\CC\}.$$
But $Z_0\subset\CC$ is an irreducible affine curve defined over $\QQ$. (Indeed, if
$d$ denotes the greatest common divisor of $|m_1|$ and $|m_2|$, then $Z_0$ is the locus of zeros of the polynomial $X^{|m_2|/d}-Y^{|m_1|/d}$ if $m_1m_2\ge0$, and it  is the locus of zeros of the polynomial $X^{|m_2|/d}Y^{|m_1|/d}-1$ if $m_1m_2<0$.) Since $Z_0$ is defined over $\QQ$, and since $Z$ satisfies Condition (3) of Lemma \ref{harald}, we must have $Z=Z_0$. Since $Z_0$ is irreducible, and since $V\subset Z=Z_0$ is a curve, it follows that $V=Z_0$, and hence
$(x^{m_1},x^{m_2})\in V$ for every $x\in\CC$.

 For every $x\in\CC^\times$, define a diagonal representation $\phi_x:H_1\to\zzle(\CC)$ by $\phi_x(\gamma_1)=\Delta(x)$ and $\phi_x(\gamma_2)=I$. Then we have
$\phi_x(\lambda)=\Delta(x^{m_1})$ and 
$\rho(\lambda)=\Delta(x^{m_2})$, 
so that $t(\phi_x)=p(x^{m_1},x^{m_2})$. Since $(x^{m_1},x^{m_2})\in V\subset
\calv=\overline{p^{-1}(Y)}$, it follows that $t(\phi_x)\in Y$ for every $x\in\CC^\times$. The map of sets $x\mapsto t(\phi_x)$ from $\CC^\times$ to $Y$ is at most two-to-one since $\trace \phi_x(\gamma_1)=x+x^{-1}$ for each $x\in\CC^\times$; in particular, this map has infinite image. On the other hand, we have $I_{\gamma_2}(t(\phi_x))=\trace\phi_x(\gamma_2)=2$ for each $x\in\CC^\times$. Thus the function $I_{\gamma_2}$ takes the value $2$ infinitely many times on the curve $Y$, and is therefore constant on $Y$. This implies that $I_{\gamma_2}$ is constant on $C$. But this contradicts \ref{off with her band-aid}. This proves \ref{that ain't hay}.

Set $L=\length(A)$. Let $W\subset\CC$ denote the set of all algebraic numbers $w$ such that (i) $w$ has degree at most $2D$ and (ii) $w$ is a root of some integer polynomial $f$ with $0<\length f\le L$. According to Lemma \ref{long reid type lemma}, the set $W$ is finite. In particular:

\Claim\label{prime rib} The set
$W':=\{w\in W:|w|\ne1\}$
is finite.
\EndClaim

I now claim:

\Claim\label{kickeroo}
 For every point $\chi\in r(\cals)\subset Y$, either (i) $\chi(\gamma)=\pm2$ for every $\gamma\in H_1$, or (ii) $\chi=p(w^{m_1},w^{m_2})$ for some $w\in W'$ and some integers $m_1$ and $m_2$, not both zero.
\EndClaim

To prove \ref{kickeroo}, suppose that $\chi\in r(\cals)$ is given and that (i) does not hold. Fix an element $\gamma_0\in H_1$ such that $\chi(\gamma_0)\ne\pm2$. Since $\chi\in r(\cals)$, there is a representation $\rho\in R(\pi_1(M))$, satisfying Conditions (a) and (b) of the statement of the proposition, such that $\chi=t(\rho|H_1)$. We have $\trace\rho(\gamma_0)=\chi(\gamma_0)\ne\pm2$, and hence after modifying $\rho$ by a conjugation in $\ggle(\CC)$ we may assume that $\rho(\gamma_0)$ is a diagonal matrix distinct from $\pm I$. Since $H_1$ is abelian, $\rho(\gamma)$ commutes with $\rho(\gamma_0)$ for every $\gamma\in H_1$, and hence $\rho(H_1)$ consists entirely of diagonal matrices. Condition (a) implies that $\rho(H_1)$ is discrete and torsion-free, and it is non-trivial since $\rho(\gamma_0)\ne I$. But any non-trivial, discrete, torsion-free group of diagonal matrices in $\zzle(\CC)$ is infinite cyclic. Hence there is a pair of generators $\gamma_1,\gamma_2$ of $H_1$ such that $\rho(\gamma_1)$ has infinite order and $\rho(\gamma_2)=I$. We may write $\rho(\gamma_1)=\Delta(w)$ for some $w\in\CC^\times$ which is not a root of unity. The discreteness of $\rho(H_1)$ implies that $|w|\ne1$. 

Since $\gamma_1$ and $\gamma_2$ generate $H_1$, there are integers $m_1$, $m_2$, $n_1$ and $n_2$ such that $\lambda=\gamma_1^{m_1}\gamma_2^{n_1}$ and $\mu=\gamma_1^{m_2}\gamma_2^{n_2}$. We cannot have $m_1=m_2=0$, since $\lambda$ and $\mu$ generate $H_1$. We now have
$\rho(\lambda)=\Delta(w^{m_1})$ and 
$\rho(\mu)=\Delta(w^{m_2})$, 
so that $\chi=t(\rho|H_1)=p(w^{m_1},w^{m_2})$. Hence condition (ii) of \ref{kickeroo} will hold provided that $w\in W'$. As I have already shown that $|w|\ne1$, it suffices to show that $w\in W$.

Since $\rho(\gamma_1)=\Delta(w)$, we have $w+w^{-1}=\trace\rho(\gamma_1)\in \trace(\rho(\pi_1(M))$. Since $\rho$ satisfies Condition (b) of the statement of the proposition, it follows that $w+w^{-1}$ is an algebraic number of degree at most $D$. Hence $w$ is an algebraic number of degree at most $2D$. Thus $w$ satisfies Condition (i) in the definition of $W$.

Since $p(w^{m_1},w^{m_2})=\chi\in r(C)$, we have $(w^{m_1},w^{m_2})\in p^{-1}(r(C))\subset V\subset Z$. Hence $g_{m_1,m_2} =A(w^{m_1},w^{m_2})=0$. Thus if $q$ denotes the least non-negative integer such that $f_{m_1,m_2}(X):=X^qg_{m_1,m_2}(X)$ is a polynomial in $X$, then $w$ is a zero of $f_{m_1,m_2}\in\ZZ[X]$. It is immediate from the definition of length and the definition of that $\length f_{m_1,m_2}\le\length A=L$. On the other hand, it follows from \ref{that ain't hay} that $\length f_{m_1,m_2} >0$. Thus $w$ satisfies Condition (ii) in the definition of $W$, taking $f=f_{m_1,m_2}$ in that condition. This proves that $w\in W$, and \ref{kickeroo} is therefore proved.

% Some of this stuff may be useful:
%For each $\rho\in S$ and for $i=1,\ldots,k$, the representation $\rho|H_i$ is well-defined up to conjugacy; and since $\rho(\pi_1(M))$ is torsion-free, either $\rho|H_i$ is injective or trivial, or $\rho(H_i)$ is infinite cyclic. Hence we may write $S$ as a finite union
%$S=Y_1\cup\cdots\cup Y_k\cup Z$,
%where $Y_i$ denotes the set of all representations $\rho\in S$ such that $H_i$ is infinite cyclic; and $Z$ denotes the set of all representations $\rho\in S$ such that each of the representations $\rho_1|H_1,\ldots,\rho_1|H_k$ is either injective or trivial. If we set $\calz=t(Z)\subset X(\pi_1(M)$, and $\caly_i=t(Y_i)\subset X(\pi_1(M)$ for $i=1,\ldots,k$, it follows that $\cals=\caly_1\cup\cdots\cup \caly_k\cup \calz$. Hence if suffices to prove that (1) $\calz$ is finite, and (2) for $i=1,\ldots,k$, the intersection of $\caly_i$ with any compact subset of $C$ is finite. }

Now let $K$ be an arbitrary compact subset of $C$. I claim that

\Claim\label{loude sing goddam}
For every $w\in W'$, the set $\{(m_1,m_2)\in\NN\times\NN:p(w^{m_1},w^{m_2})\in r(K)
\}$ is finite.
\EndClaim

To prove \ref{loude sing goddam}, note that since $K$ is compact, the functions $I_\lambda$ and $I_\mu$ are bounded on $K$. Fix a positive constant $B$ such that $|I_\lambda(\chi)|\le B$ and $|I_\mu(\chi)|\le B$ for every $\chi\in K$. Now let $w\in W$ be given, and set $R=\max(|w|,|w|^{-1})$. The definition of $W'$ implies that $|w|\ne1$ and hence $R>1$. Consider any $\{(m_1,m_2)\in\NN\times\NN$ such that $p(w^{m_1},w^{m_2})\in r(K)$, and choose $\chi\in K$ so that $r(\chi)=p(w^{m_1},w^{m_2})$. The definitions give 
$$I_\lambda(\chi)=I_\lambda(r(\chi))=\trace \Delta(w^{m_1})=w^{m_1}+ w^{-m_1},$$
so that  
$$ B \ge|I_\lambda(\chi)|\ge R^{|m_1|}-R^{-|m_1|}\ge R^{|m_1|}-1$$
and hence $|m_1|\le\log(B+1)/(\log R)$. The same argument, with $\mu$ in place of $\lambda$, shows that $|m_2|\le\log(B+1)/(\log R)$. This establishes the finiteness assertion \ref{loude sing goddam}.

I will now complete the proof of the proposition by showing that $\cals\cap K$ is finite. In view of \ref{corson porson}, it suffices to show that $r(\cals\cap K)$ is finite. According to \ref{kickeroo}, we have $r(\cals)\subset \Omega_1\cup \Omega_2$, where 
$$\Omega_1=\{\chi\in X(H_1):\chi(\gamma)=\pm2\text{ for every }\gamma\in H_1\}$$
and
$$\Omega_2=\{p(w^{m_1},w^{m_2}): w\in W',m_1,m_2\in\ZZ\}.$$
According to the discussion in \ref{erence}, $\Omega_1$ consists of all points of the algebraic set $X(H_1)$ whose coordinates in a suitable coordinate system are all $\pm2$, and is therefore a finite set. Since 
$$r(\cals\cap K)\subset r(\cals)\cap r(K)\subset \Omega_1\cup(\Omega_2\cap r(K)),$$
it suffices to show that $\Omega_2\cap r(K)$ is finite. We may write
$\Omega_2=\bigcup_{w\in W'}\Omega_2^w$, where
$$\Omega_2^w=\{p(w^{m_1},w^{m_2}):m_1,m_2\in\ZZ\}$$
for each $w\in W'$. Hence 
$$\Omega_2\cap r(K)=\bigcup_{w\in W'}(\Omega_2^w\cap r(K)).$$
Here $W'$ is finite by \ref{kickeroo}, and it follows from \ref{loude sing goddam} that $\Omega_2^w\cap r(K)$ is finite for each $w\in W'$. Hence $\Omega_2\cap r(K)$ is finite as required.
\EndProof 
%YZf_

\begin{proposition}\label{free analogue}
 Let $F$ denote a free group of rank $2$ on generators $\xi$ and $\eta$.  Let $C$ be a curve contained in the $\zzle(\CC)$-character variety of $F$ and defined over $\overline\QQ$.
Let $D$ be a positive integer, let $\alpha$ be a positive number less than $\log3$, and let $\calt$ denote the set of points of $C$ of the form $t(\rho)$, where $\rho$ is a
representation such that (1) $\rho(\pi_1(M))$ is discrete, non-elementary and torsion-free (2) $\max(d_P([\rho(\xi)]),d_P([\rho(\eta)]))\le\alpha$ for some $P\in\HH^3$, and (3) $\trace(\rho(F))\subset L$ for some number field $L$ of degree at most $D$.  Then $\calt$ is a finite set. 
\end{proposition}

\Remark\label{more watson}It follows from \ref{watson} that Condition (1) is equivalent to the condition that $\Pi_\CC(\rho(\pi_1(M)))\le\pizzle(\CC)$ is discrete, non-elementary and torsion-free. 
\EndRemark

\Proof[Proof of Proposition \ref{free analogue}]

Assume that $\calt$ is infinite, and choose a sequence $(\chi_j)$ of distinct points in $\calt$. For each $j$ let $\rho_j:F\to\zzle(\CC)$  be a representation with $t(\rho_j)=\chi_j$.
According to the definition of $\calt$, for each $j$ there is a
a point $P_j\in\HH^3$   such that 
$\max(d_{P_j}([\rho_j(\xi)]),d_{P_j}([\rho_j(\eta)]))\le\alpha$.
After modifying each $\rho_j$ by a suitable conjugation in
$\zzle(\CC)$, we may assume that the $P_j$ are all the same point of $\HH^3$,    which I will denote by $P$.    Thus for each $j$ we have
\begin{equation}\label{saskatoon}\max(d_{P}([\rho_j(\xi)]),d_{P}([\rho_j(\eta)]))\le\alpha.\end{equation}

It follows from (\ref{saskatoon}) that the $\rho_j$ all lie
in a compact subset of $R(F)$. Hence after passing to a
subsequence we may assume that $\rho_j\to \rho_\infty$
for some $\rho_\infty\in R(F)$. Then the sequence $[\rho_j]_{j\in\NN}$ converges algebraically (see, for example, Section 2 of \cite{finiteness}) to $[\rho_\infty]$.

For $1\le j\le\infty$, set
$[\rho_j]=\Pi_\CC\circ\rho_j:\Gamma\to\pizzle(\CC)$, and set $\Gamma_j=[\rho_j](F)$.
For each $j\in\NN$, since $\rho_j\in\calt$, the group $\rho_j(F)$ is discrete, non-elementary and torsion-free. Hence by \ref{watson}, $\Gamma_j$ is discrete, non-elementary and torsion-free, and $\Pi_\CC$ maps $\rho_j(F)$ isomorphically onto $\Gamma_j$. It then follows from Assertion (1) of \cite[Theorem 2.4]{finiteness} that  $[\rho_\infty](F)$ is discrete and non-elementary. Furthermore, it follows from Assertion (2) of \cite[Theorem 2.4]{finiteness}  that for sufficiently large $j<\infty$ there is a homomorphism $\psi_j:\Gamma _\infty\to\Gamma_j$ such that $\psi_j\circ[\rho_\infty]=[\rho_j]$. After passing to a subsequence we may assume that such a $\psi_j$ exists for every $j\in\NN$.

We regard $\psi_j$ as a representation of $\Gamma_\infty$ in $\pizzle(\CC)$, and we let $\psi_\infty:\Gamma_\infty\to \pizzle(\CC)$ denote the inclusion homomorphism. Since the sequence $([\rho_j])_{i\in\NN}$ converges algebraically to $[\rho_\infty]$, the
sequence $(\psi_j)_{i\in\NN}$ converges algebraically to $\psi_\infty$.

If $\gamma\in \Gamma_\infty$ has finite order, then for each $j\in\NN$ we have $\psi_j(\gamma)=1$, since $\Gamma_j$ is torsion-free. By algebraic convergence we have 
$$\gamma=\psi_\infty(\gamma)=\lim_{j\to\infty}\psi_j(\gamma)=1.$$
This shows that $\Gamma_\infty$ is
torsion-free. Hence $M_\infty=\HH^3/\Gamma_\infty$ acquires the structure of a hyperbolic $3$-manifold in a natural way.

We now distinguish two cases; in each case we shall obtain a contradiction, thus completing the proof.

{\bf Case I: $\vol M_\infty=\infty$.} 

Note that $[\rho_\infty(\xi)]$ and $[\rho_\infty(\eta)]$ do not commute, since $\Gamma_\infty$ is non-elementary. It therefore follows from Proposition \ref{hunt it down and kill it} that in Case I we have $$\max(d_P([\rho_\infty(\xi)] ),d_P([\rho_\infty(\eta)] ))\ge\log3.$$ On the other hand, since $[\rho_j(\xi)]\to[\rho_\infty(\xi)]$ and $[\rho_j(\eta)]\to[\rho_\infty(\eta)]$, it follows from (\ref{saskatoon}) that $$\max(d_{P}([\rho_j(\xi)]),d_{P}([\rho_j(\eta)]))\le\alpha<\log3.$$
This gives the required contradiction in this case.

{\bf Case II: $\vol M_\infty<\infty$.}

Set $N=\ker[\rho_\infty]\triangleleft F$.  For every $j\in\NN$, since $\psi_j\circ[\rho_\infty]=[\rho_j]$,
we have $[\rho_j](N)=\{1\}$. Since $\rho_j(F)$ is torsion-free, we in fact have $\rho_j(N)=\{1\}$ for every $j\in\NN$. Since $\rho_j\to\rho_\infty$, we have have $\rho_\infty(N)=\{1\}$. Thus $N\le\ker\rho_\infty\le\ker[\rho_\infty]=N$, and hence $\ker\rho_\infty=N$.

According to \ref{weinermandias}, the character variety $X(\Gamma_\infty)$ may be regarded as the subset of $X(F)$ consisting of all points of the form $t(\rho)$ where $\rho:F\to\zzle(\CC)$ is a representation such that $\rho(N)=\{1\}$. Hence $\chi_j\in X(\Gamma_\infty)$ for every $j\in\NN$. Thus $C\cap X(\Gamma_\infty)$ is infinite. Since $C$ is a curve and $X(\Gamma_\infty)$ is a (closed) algebraic set, it follows that $C\subset X(\Gamma_\infty)$. In addition we have $\chi_\infty\in C$; and $\chi_\infty$ is the character of a discrete, faithful representation of $\Gamma_\infty$, since $\rho_\infty$ has kernel $N$ and has discrete image. By hyporhesis $\CC$ is defined over $\overline\QQ$.
Thus $C$ satisfies the hypotheses of Proposition \ref{long reid type prop}. Since the $\chi_i$ belong to $\calt$, they belong to the set $\cals$ defined in the statement of  Proposition \ref{long reid type prop}. However, the infinite set $\{\chi_j:j\in\NN\}$ is contained in a compact subset of $C$, and  Proposition \ref{long reid type prop} asserts that $\cals$ has finite intersection with every compact subset of $C$.
This gives the required contradiction in this case.
\EndProof

\section{Trace fields and Margulis numgers}\label{masterpiece section}

This section contains the proof of the motivating result of the paper, Theorem \ref{motivational research}, which was stated in the Introduction as Theorem A.

\Lemma\label{uneeda biscotte}
Let  $K$ be a number field, and let $m$ and $N$ be positive integers with $m>2$. Then the ring of integers $\ok$ of $K$ has a finite subset 
$W$ with the following property. 
Let $\gamma$ be an element of
$\zzle(\ok)$ such that $\trace\gamma\notin W$. Then there exists a valuation $w$ of $K$ such that (1) the characteristic of $k_w$ is greater than $N$, and (2)
$\overline h_w([\gamma])\in \pizzle(k_w)$ has
order $m$.
\EndLemma

\Proof
Recall from \cite[Chapter III, Section 1]{neukirch} that a {\it place} of a number field is an equivalence class of absolute values.

If $\frakp$ is a place represented by a valuation $v$, I will set $\op=\frako_v$, $k_\frakp=k_v$, $h_\frakp=h_v$, and $\overline h_\frakp=\overline h_v$, where $\frako_v$, $k_v$ and $h_v$ are defined as in \ref{bet you can't eat one}. I will also set 
$\kappa_\frakp=\Pi_{k_\frakp}\circ h_\frakp:\zzle(\op)\to\pizzle(k_\frakp)$, where $\Pi_{k_\frakp}$ is defined as in \ref{lions and tigers and bears}. 
%the homomorphism defined by $\kappa_\frakp(\gamma )=[h_\frakp(\gamma )]$.
This definition of $\kappa_\frakp$ is a paraphrase of the definition given in \cite[Subsection 2.4]{finitistic}. Comparing it with the definition of $\overline h_v$ given in \ref{bet you can't eat one}, we find that $\kappa_\frakp=\overline h_\frakp\circ\Pi_{\op}$.
%, using the conventions introduced in \ref{bet you can't eat one} and \ref{lions and tigers and bears} above.''

%\redcomment{Everything needs fixing up.}

%In this proof I will be using the same conventions regarding algebraic number theory as in \cite[Section 2]{finitistic}. Thus if $\frakp$ be a nonarchimedean place of a number field $E$, I will
%denote by $\op$ the valuation ring defined by $\frakp$, by  $k_\frakp$  the
%residue field $\op/\frakp$, and by $\eta_\frakp:\op\to k_\frakp$ the quotient
%homomorphism. I will denote by
%$h_\frakp$ the natural homomorphism $\zzle(\op)\to\zzle(k_\frakp)$, defined by
%$$\begin{pmatrix}a&b\\c&d\end{pmatrix}\mapsto\begin{pmatrix}\eta_\frakp(
%a)&\eta_\frakp( b)\\\eta_\frakp( c)&\eta_\frakp( d)\end{pmatrix}.$$
%I will denote by
%$\kappa_\frakp:\zzle(\op)\to\pizzle(k_\frakp)$
%the homomorphism defined by $\kappa_\frakp(\gamma )=[h_\frakp(\gamma )]$. \redcomment{Fix the above as needed. Revise references to valuations and places throughout the paper as appropriate. Consider referring more to \cite{neukirch}. Also use whatever version of this is appropriate: ``These definitions of $h_\frakp$ and $\kappa_\frakp$ are paraphrases of the definitions given in \cite[Subsection 2.4]{finitistic}, using the conventions introduced in \ref{bet you can't eat one} and \ref{lions and tigers and bears} above.''
%}

For a given prime $p$, there are at most finitely many finite places of $K$ whose residue fields have characteristic $p$. Indeed, it follows from \cite[Proposition 3.7]{neukirch} that any (exponential) valuation of $K$ whose residue fields has characteristic $p$ restricts to a multiple of the $p$-adic valuation on $\QQ$, and the latter admits only finitely many extensions to $K$, for example by \cite[Proposition 8.2]{neukirch}.

Hence we may fix a finite set $S$  of places of $K$ which includes all infinite places and all finite places whose residue fields have characteristic at most $N$. As $S$ is finite and contains all infinite places, it is by definition \cite[Subsection 2.3]{finitistic} an {\it admissible} set of places. 
Let $\ooks$ denote the ring of $S$-integers of $K$.

Since $S$ is an admissible set of places of the number field $K$,
Proposition 2.7 of \cite{finitistic} asserts that for every integer $m>2$, there is a finite set 
$W\subset\ooks$ with the following property. 
\begin{claim}\label{gadzooks}
Let $\gamma$ be an element of
$\zzle(\ooks)$ such that $\trace\gamma\notin W$. Then there exists a
place $\frakp$ of $K$, with $\frakp\notin S$, such that $\kappa_\frakp(\gamma)\in \pizzle(k_\frakp)$ has
order $m$.
\end{claim}

%\redcomment{Fix from here and
%Finish.
%Recall that in \ref{bet you can't eat one} I defined $\overline h_\zeta: %\piggle(R)\to\piggle(S)$ by $\overline h_\zeta\circ\Pi_R= \Pi_S\circ h_\zeta$.
%Relate $\kappa_\frakp$ to $\overline h_v$. 
%Some notation about places, valuations, ideals, etc. has to be reconciled.
%

Now let $\gamma\in\zzle(\ok)$ be given with $\trace\gamma\in W$. In particular we have $\gamma\in\zzle(\ooks)$, and it follows from \ref{gadzooks} that there
is a
place $\frakp$ of $K$, with $\frakp\notin S$, such that $\kappa_\frakp(\gamma)\in \pizzle(k_\frakp)$ has
order $m$. Let $v$ be a valuation representing $\frakp$. Since $\frakp\notin S$, the characteristic of $k_w=k_\frakp$ is greater than $N$. Furthermore, the order of
$\overline h_w([\gamma])=\overline h_\frakp\circ\Pi_{\op}(\gamma])=\kappa_\frakp(\gamma)$ in $\pizzle(k_w)=\pizzle(k_\frakp)$ has
order $m$. 
\EndProof

\Lemma\label{just about} Let $K$ be a number field, and let $F$ denote a free group on two generators $\xi $ and $\eta$. Then up to conjugacy there are at most finitely many representations of $\rho:F\to\zzle(\CC)$ such that 
\begin{enumerate}
\item $\rho(F)$ is discrete, torsion-free and cocompact:
\item $\max(d_P([\rho(\xi )]),d_P([\rho(\eta)]))\le0.34$ for some $P\in\HH^3$; and
\item $\trace\rho(F)\subset K$.
\end{enumerate}
\EndLemma

\Proof
Fix a natural number $N\ge7$ having the property stated in Theorem \ref{N lives}. 
With this choice of $N$, and with $m=7$, fix a set $W\subset\ok$ having the property stated in Lemma \ref{uneeda biscotte}.

I now claim:
\Claim\label{diddle i} If a representation $\rho:F\to\zzle(\CC)$ satisfies Conditions (1)--(3) of the statement of the present lemma, then $\trace\rho(\xi )$ and $\trace\rho(\eta)$ belong to $W$.
\EndClaim

To prove \ref{diddle i}, assume that $\rho:F\to\zzle(\CC)$ satisfies Conditions (1)--(3), and that either $\trace\rho(\xi )$ or $\trace\rho(\eta)$ lies outside $W$; by symmetry we may assume that $t:=\trace\rho(\xi )\notin W$. Since $\rho$ satisfies Condition (1), it follows from Proposition \ref{just conjugate} that $\rho$ is conjugate in $\ggle(\CC)$ to a representation of $\Gamma$ in $\zzle(E)$ for some number field $E$. Hence we may assume without loss of generality that $\rho(F)\subset\zzle(E)$. After enlarging $E$ if necessary we may also assume that  $E\supset K$. Set $\Gamma=\Pi_\CC(\rho(F))$. Condition (1) implies that $\Gamma$  is discrete, torsion-free and cocompact. Let $M$ denote the closed, orientable hyperbolic $3$-manifold $\HH^3/\Gamma$. Set 
$x=[\rho(\xi)]\in\Gamma$ and $y=[\rho(\eta)]\in\Gamma$, so that
$\Gamma=\langle x,y\rangle$. Since $\Gamma$ is cocompact, $x$ and $y$ do not commute. 

If $t$ is not an algebraic integer, it follows from Corollary \ref{but who needs it} that $\max( d_P(x), d_P(y))>0.34$. This contradicts Condition (3). Hence $t$ is an algebraic integer. 
In view of Condition (3) we therefore have $t\in\ok$.

Since $x$ and $y$ do not commute, we have $\rho(\xi)\ne\pm1$. Hence $\rho(\xi)$ is conjugate in $\zzle(E)$ to a matrix of the form 
$\begin{pmatrix}0&b\cr1&d\end{pmatrix}$. Since $\det\rho(\xi)=1$ and $\trace\rho(\xi)=t$ we must have $b=-1$ and $d=t$. After modifying $\rho$ by a conjugation in $\zzle(E)$ we may therefore assume that $$\rho(\xi)=\begin{pmatrix}0&-1\cr1&t\end{pmatrix}\in\zzle(\ok).$$

Since $\rho(\xi) \in\zzle(\ok)$ and $t\notin W$, the property of $W$ stated in Lemma \ref{uneeda biscotte}  implies that there is a valuation $w$ of $K$ such that the characteristic of $k_w$ is greater than $N$, and 
$\overline h_w(x)\in \pizzle(k_w)$ has order $7$. By the extension theorem for valuations (see for example \cite[Proposition 8.2]{neukirch}, there exists a (not necessarily normalized) valuation $v$ of $E$ such that $v|(K-\{0\})=\cdot w$. We then have $\frako_w\subset\frako_v$ and $\frakm_v\cap\frako_w=\frakm_w$. It follows that $k_w$ is a subfield of $k_v$, so that in particular $k_v$ has characteristic greater than $N$. It also follows that $\ker \overline h_w=\ker \overline h_v\cap \pizzle(K)$; hence $\overline h_v(x)\in \pizzle(k_v)$ has order $7$. The property of $N$ stated in Theorem \ref{N lives} now implies that
$\max( d_P(x), d_P(y))>0.34$. This contradicts Condition (2), and \ref{diddle i} is proved.

Now, according to Proposition \ref{i never sausage a bad pun}, the character variety $X(F)$ (defined in terms of the generating system $(\xi,\eta)$ is equal to $\CC^3$, and the map $t:R(F)\to X(F)$ is given by $\rho\mapsto(\trace\rho(\xi),\trace\rho(\eta),\trace\rho(\xi\eta))$. If $\calr$ denotes the set of all $\rho\in R(F)$ that satisfy Conditions (1)---(3) of the statement of the present lemma, it follows from \ref{diddle i} that 
\Equation\label{pack of cards}t(\calr)\subset W\times W\times\CC.
\EndEquation

Set $D=[K:\QQ]$, and let $\calr'$ denote the set of all $\rho\in R(F)$ that satisfy Conditions (1), (2) and

\quad(3$'$)\quad $\trace(\rho(\pi_1(M))\subset L$ for some number field $L$ of degree at most $D$.

We obviously have $\calr\subset\calr'$. According to Proposition \ref{free analogue}, applied with $\alpha=0.34<\log3$, the set $t(\calr')\cap C$ is finite for every curve $C\subset X(F)$ which is defined over $\overline\QQ$. In particular this holds for $C= \{w_1\}\times \{w_2\}\times\CC$, since $W\subset\ok\subset\overline\QQ$. Thus the set $t(\calr')\cap( \{w_1\}\times \{w_2\}\times\CC)$ is finite for all $w_1,w_2\in W$. But (\ref{pack of cards}) gives
$$t(\calr)=\bigcup_{w_1,w_2\in W}(t(\calr)\cap( \{w_1\}\times \{w_2\}\times\CC))\subset
\bigcup_{w_1,w_2\in W}(t(\calr')\cap( \{w_1\}\times \{w_2\}\times\CC));$$
thus $t(\calr)$ is contained in a finite union of finite sets and is therefore finite. 

Finally, every $\rho\in\calr$ satisfies Condition (1) and is therefore irreducible. Hence by \cite[Proposition 1.5.2]{varreps}, representations in $\calr$ having the same image under $t$ are conjugate. Since $t(\calr)$ is finite, $\calr$ is a finite union of conjugacy classes of representations. This is the conclusion of the lemma.
\EndProof

\Theorem\label{motivational research} Let $K$ be a number field. Then up to isometry, among the closed, orientable hyperbolic $3$-manifolds that have trace field $K$, all but (at most) finitely many admit $0.34$ as a Margulis number.
\EndTheorem

\Proof
Let $F$ denote a free group on two generators $\xi $ and $\eta$.
According to Lemma \ref{just about}, for some $n\ge0$ there are representations $\rho_1,\ldots,\rho_n$ of $F$ in $\zzle(\CC)$, satisfying Conditions (1)---(3) of the lemma, such that every representation $\rho:F\to\pizzle(\CC)$ satisfying Conditions (1)---(3) is conjugate to one of the $\rho_i$. Condition (1) implies that for $1=1,\ldots,n$, the quotient $M_i=\HH^3/\rho_i(F)$ is a closed hyperbolic $3$-manifold. Set $\cald=\max_{1\le i\le n}\diam M_i$.

Suppose that $M=\HH^3/\Gamma$ is a closed hyperbolic $3$-manifold that has trace field $K$ and does not admit $0.34$ as a Margulis number. According to Definitions \ref{Margdef} this means that there exist non-commuting elements $x $ and $y$ of $\Gamma$ and a point $P\in\HH^3$ such that $\max(d_P(x),d _P(y))\le0.34$ for some $P\in\HH^3$. The subgroup $\langle x ,y\rangle$ must have finite covolume, as otherwise we would have $\max(d_P(x),d _P(y))\ge\log3$ by Proposition \ref{hunt it down and kill it}. It follows that $\langle x ,y\rangle$ has finite index in $\Gamma$. 

According to \cite[Proposition 3.1.1]{varreps}, there is a subgroup of $\zzle(\CC)$ which maps isomorphically onto $\langle x ,y\rangle$ under $\Pi_\CC$. Let $\tx$ and $\ty$ denote the generators of this group that map to $x$ and $y$.

It now follows that if and if we define a representation $\rho:F\to\pizzle(\CC)$  by $\rho(\xi)=\tx$, $\rho(\eta)=\ty$, then Conditions (1)---(3) of Lemma \ref{just about} hold. Hence $\rho$ is conjugate to $\rho_i$ for some $i$, and hence $M$ admits $M_i$ as a finite-sheeted covering. In particular $\diam M\le\diam M_i\le\cald$. 

Let $v$ denote the infimum of the volumes of all
closed hyperbolic $3$-manifolds;  we have $v>0$, for example by
\cite[Theorem 1]{meyerhoff}. Let $\calm$ denote the class of closed hyperbolic $3$-manifolds that have trace field $K$, and do not admit $0.34$ as a Margulis number. The manifolds in $\calm$ have volume at least $v$,
constant curvature $-1$ and dimension $3$, and I have shown that they have diameter at most $\Delta$. It then
follows from the main theorem of \cite{peters} that the manifolds in $\calm$
represent only finitely many diffeomorphism types. By the Mostow
rigidity theorem, they represent only finitely many isometry
types.
\EndProof

\bibliographystyle{plain}
\bibliography{cubic-not}

\end{document}